\documentclass[12pt,fleqn]{article}
\usepackage{latexsym}
\usepackage{amsmath}
\usepackage{amssymb}
\usepackage{epsfig}
\usepackage{amsfonts}

\begin{document}

\newtheorem{definition}{Definition}[section]
\newtheorem{theorem}[definition]{Theorem}
\newtheorem{lemma}[definition]{Lemma}
\newtheorem{remark}[definition]{Remark}
\newtheorem{remarks}[definition]{Remarks}
\newtheorem{examples}[definition]{Examples}
\newtheorem{leeg}[definition]{}
\newtheorem{definitions}[definition]{Definitions}
\newtheorem{proposition}[definition]{Proposition}
\newtheorem{example}[definition]{Example}
\newtheorem{comments}[definition]{Some comments}
\newtheorem{corollary}[definition]{Corollary}
\def\square{\Box}
\newtheorem{observation}[definition]{Observation}
\newtheorem{observations}[definition]{Observations}
\newtheorem{defsobs}[definition]{Definitions and Observations}
\newenvironment{prf}[1]{ \trivlist
\item[\hskip \labelsep{\it
#1.\hspace*{.3em}}]}{~\hspace{\fill}~$\square$\endtrivlist}
\newenvironment{proof}{\begin{prf}{Proof}}{\end{prf}}

\title{Galois theory of $q$-difference equations}
\author{Marius van der Put $\ $ and $\ $ Marc Reversat \\
\footnotesize Department of Mathematics, University of Groningen,
 P.O.Box 800,\\
\footnotesize 9700 AV Groningen, 
The Netherlands, mvdput@math.rug.nl,$\;\;\;$ and\\
\footnotesize Laboratoire de Math\'ematiques Emile Picard, Universit\'e Paul 
Sabatier, U.F.R. M.I.G.\\
\footnotesize 118, route de Narbonne, 31602 Toulouse, France,  
marc.reversat@math.ups-tlse.fr }
\maketitle
\noindent

\section*{Introduction}
Choose $q\in {\mathbb C}$ with $0<|q|<1$. The main theme of this paper is the
study of linear $q$-difference equations over the field $K$ of germs of 
meromorphic functions at $0$. A more detailed and systematic treatment of
classification and moduli is developed as a continuation of [vdP-S1] (Chapter
12), [vdP-R] and [vdP]. It turns out that a difference module $M$ 
over $K$ induces in a functorial way a vector bundle $v(M)$ on the
Tate curve $E_q:={\mathbb C}^*/q^{\mathbb Z}$ 
(this is done here for all slopes, the case of integral slopes has been 
treated in [Sau1] and [Sau2]).
As a corollary one rediscovers Atiyah's classification ([At]) of the
indecomposable vector bundles on the  
complex Tate curve. Linear $q$-difference equations are also
studied in positive characteristic in order to derive Atiyah's results for
elliptic curves for which the $j$-invariant is not algebraic over 
${\mathbb F}_p$. 

A universal difference ring and a universal formal difference Galois group is 
introduced. For pure difference
modules this ring provides an explicit expression of the difference Galois
group. If the difference module has more than one slope, then part of the
difference Galois group has an interpretation as `Stokes matrices', related
to a summation method for divergent solutions. 
We do not make any hypothesis on the slopes, when they are integers see [R-S-Z, Sau2]. The above moduli space is the
algebraic tool to compute this part of the difference Galois group.

It is possible to provide the vector bundle $v(M)$ on $E_q$, corresponding to a
difference module $M$ over $K$, with a connection $\nabla _M$. If $M$ is 
regular singular, then $\nabla _M$ is essentially determined by the
absense of singularities and `unit circle monodromy'. More precisely, the
monodromy of the connection $(v(M),\nabla _M)$ coincides with the action of
two topological generators of the universal regular singular difference Galois
group 
([vdP-S1, Sau1]). 
For irregular difference modules, $\nabla _M$ will have singularities
and there are various Tannakian choices for $M\mapsto (v(M),\nabla _M )$. 
Explicit computations are difficult, especially for the case of non integer 
slopes.

The case of modules with integer slopes, 
has been studied in [R-S-Z]. This answers a question of  G.D.~Birkhoff and 
follows ideas of G.D.~Birkhoff, P.E.~Guenther, C.R.~Adams (see [Bir]). 

\section{Classification of $q$-difference equations}

\subsection{Some notation and formulas}
A difference ring is a commutative $R$ with a given automorphism $\phi$.
The skew ring of difference operators $R[\Phi ,\Phi ^{-1}]$, consists of the 
finite formal sums $\sum _{n\in {\mathbb Z}} a_n\Phi ^n$ with all $a_n\in R$.
The multiplication is defined by $\Phi r=\phi (r)\Phi$. A difference
module $M$ is a left $R[\Phi ,\Phi ^{-1}]$-module which is free and finitely
generated as $R$-module. The action $\Phi _M$ of $\Phi$ on $M$ is an additive
bijection satisfying $\Phi _M(rm)=\phi (r)\Phi _M(m)$. Thus we may describe
a difference module as a pair $(M,F)$, with $F$ an additive bijective map and
such that $F(fm)=\phi (m)F(m)$ holds. 

\smallskip

As before we choose $q\in {\mathbb C}$ with $0<|q|<1$. Further we fix $\tau$
in the complex upper half plane with $e^{2\pi i \tau}=q$. The fields 
$K={\mathbb C}(\{ z \})$ and $\widehat{K}={\mathbb C}((z))$ are made into 
difference fields by the automorphisms $\phi$ given by $\phi (z)=qz$. These
automorphisms are extended to the finite extensions $K_n$ and $\widehat{K}_n$
of degrees $n$ and to the algebraic closures $K_\infty$ and
$\widehat{K}_\infty$ of $K$ and $\widehat{K}$, by 
$\phi (z^\lambda )=q^\lambda z^\lambda$ where 
$q^\lambda :=e^{2\pi i \tau \lambda}$ for all $\lambda \in {\mathbb Q}$. 

\smallskip

The formula $\phi (z^\lambda )=q^\lambda z^\lambda$ makes
${\mathbb C}[z,z^{-1}]$ into a difference ring. A difference module over this
ring will be called a {\it global difference module}. As we will prove later
on, any difference module $M$ over $K$ is obtained from a unique global module
$M_{global}$ as a tensor product, i.e., 
$M\cong K\otimes _{{\mathbb C}[z,z^{-1}]}M_{global}$.

Other difference rings that we will use are ${\mathbb C}[z^{1/n},z^{-1/n}]$
and  $O=O({\mathbb C}^*)$, the ring of the holomorphic functions on 
${\mathbb C}^*$.

\smallskip

Closely related to $q$-difference equations is the {\it complex Tate curve}
$E_q:={\mathbb C}^*/q^{\mathbb Z}$. We write $pr:{\mathbb C}^*\rightarrow E_q$
for the natural map. {\it Theta functions} are related to both $E_q$ and
$q$-difference equations. Put
$\Theta :=\sum _{n\in {\mathbb Z}}q^{n(n-1)/2}(-z)^n$. Then
\[\Theta =d\prod _{n>0}(1-q^nz^{-1})\cdot \prod _{n\geq 0}(1-q^nz)
\mbox{ for some constant } d\neq 0.\]
The divisor of $\Theta$ on ${\mathbb C}^*$ is $\sum _{n\in {\mathbb Z}}[q^n]$.
Further $-z\Theta (qz)=\Theta (z)$ or $-z\phi (\Theta )=\Theta$. The latter 
implies $\frac{dz}{z}+\phi (\frac{d\Theta}{\Theta})=\frac{d\Theta}{\Theta}$.
Moreover, the poles of $\frac{d\Theta}{\Theta}$ form the set $q^{\mathbb Z}$.
Each pole is simple and has residue 1.

\smallskip

For $c\in {\mathbb C}^*$ one defines 
$\theta _c:=\frac{\Theta (cz)}{\Theta (z)}$. This function has the property:
$c\cdot \theta _c(qz)=\theta _c(z)$. Moreover, the differential
form $\omega _c:=\frac{d\theta _c}{\theta _c}$ is $\phi$-invariant and
defines a differential form on $E_q$. If $c\not \in q^{\mathbb Z}$, then
$\omega _c$ has simple poles in $pr (c^{-1})$ and $1=pr(1)$ with residues
$1$ and $-1$. Further $\omega _{qc}=\omega _c-\frac{dz}{z}$.

\subsection{Regular singular difference modules}
We recall some classical results (a modern proof is given in [vdP-S]). 
The classification of regular singular modules over $K$ and $\widehat{K}$ are 
similar and we restrict our attention to $q$-difference modules $M$ over $K$. A
difference module over $K$ is called {\it regular singular}
if there exists a lattice $M^0\subset M$
over ${\mathbb C}\{z\}$ (i.e., $M^0={\mathbb C}\{z\}e_1\oplus \cdots \oplus 
{\mathbb C}\{z\}e_m$ for some $K$-basis $\{e_1,\dots ,e_m\}$ of $M$)
which is invariant under $\Phi$ and $\Phi ^{-1}$
(or in later terminology, $M$ is pure of slope 0). Then $M$ can uniquely be 
written as $K\otimes _{\mathbb C}V$ where $V$ is a finite dimensional vector
space over $\mathbb C$ provided with a linear map $A:V\rightarrow V$ such that
all its eigenvalues $\alpha$ satisfy $|q| <|\alpha |\leq 1$. The action of
$\Phi$ on this tensor product is given by 
$\Phi (a\otimes v)=\phi (a)\otimes A(v)$, for $a\in K$ and $v\in V$.

For a regular singular $M$ we define $M_{global}\subset M$ as the set of
elements $m$ such that the $\mathbb C$-vector space generated by 
 $\{\Phi ^nm|\ n\geq 0\}$ has finite dimension. Equivalently,
 $m\in M_{global}$, if and only if there exists
a non zero $L\in {\mathbb C}[\Phi ]$ such that $L(m)=0$. Clearly
$M_{global}$ is a ${\mathbb C}[z,z^{-1}]$-submodule. More precisely,

\begin{lemma}
$M_{global}={\mathbb C}[z,z^{-1}]\otimes _{\mathbb C}V$ and consequently
the natural morphism $K\otimes _{{\mathbb C}[z,z^{-1}]}M_{global}\rightarrow
M$ is an isomorphism.
\end{lemma}
\begin{proof} Suppose that $m\in M$ (or even $m\in \widehat{K}\otimes _KM$)
satisfies $L(m)=0$ with $L=\Phi ^d+c_{d-1}\Phi ^{d-1}+\cdots +c_0\in
{\mathbb C}[\Phi ]$ and $c_0\neq 0$. Write 
$m=\sum _{n>>-\infty} z^n\otimes v_n$.
Then $L(m)=\sum _nz^n\otimes 
(q^{nd}A^d+c_{d-1}q^{n(d-1)}A^{d-1}+\cdots +c_0)v_n$.
One provides $V$ with some norm. For large $|n|$, the linear map 
$q^{nd}A^d+c_{d-1}q^{n(d-1)}A^{d-1}+\cdots +c_0$ is invertible since the
norm of either $c_0$ or $q^{nd}A^d$ is larger than the norm of the remaining 
part of the linear map. Thus $M_{global}\subset {\mathbb C}[z,z^{-1}]\otimes
V$. The other inclusion is obvious. \end{proof}

\begin{remarks}$\;$\\ {\rm
(1) The {\it unipotent difference module} $U_m$ over $K$ (or over
$\widehat{K}$) is $U_m:=K\otimes _{\mathbb C}{\mathbb C}^m$
with $\Phi (f\otimes v)=\phi (f)\otimes A(v)$, where 
$A:{\mathbb C}^m\rightarrow {\mathbb C}^m$ is the unipotent map
which has a unique Jordan block. Any 1-dimensional regular singular
difference module has the form $E(c):=Ke$ with $\Phi (e)=ce,\ c\in {\mathbb
C}^* $ and one may
normalize $c$ such that $|q|<|c|\leq 1$. From the modules 
$\{E(c) \}_{|q|<|c|\leq 1}$ and $U_2$ one constructs every regular singular
module by taking tensor products and direct sums.

\smallskip

\noindent (2) Let $M$ be a regular singular module.
For any $c\in {\mathbb C}^*$,  $Eigen(\Phi ,c)\subset M$ denotes the 
generalized eigenspace for the eigenvalue $c$. In other words, 
$Eigen(\Phi ,c)$ consists of the elements $m\in M$ such that there exists an
integer $N>0$ with $(\Phi -c)^N(m)=0$. From the above one easily concludes
that each $Eigen(\Phi ,c)$ has finite dimension. Further 
$V=\oplus _{|q|<|c|\leq 1}Eigen(\Phi ,c)$ and
$M_{global}=\oplus _{c\neq 0}Eigen(\Phi ,c)$. }\end{remarks}

\subsection{The slope filtration}
We describe here the slope filtration and give references for more details
and proofs.
It is well known that any difference module $M$ over $K$ contains a
cyclic vector. This means that there exists an element $e\in M$ such that the 
homomorphism $K[\Phi ,\Phi ^{-1}]\rightarrow M$, given by 
$\sum a_n\Phi ^n\mapsto \sum a_n\Phi ^n(e)$, is surjective (compare 
Lemma 4.1). Thus $M$ is isomorphic to
$K[\Phi ,\Phi ^{-1}]/K[\Phi ,\Phi ^{-1}]L$ for some $L$ of the form
$\Phi ^d+a_{d-1}\Phi ^{d-1}+\cdots +a_0$, with all $a_i\in K$ and $a_0\neq 0$. 
The difference operator $L$ has a Newton polygon. For completeness we recall
its definition. Let $ord:\ \widehat{K}\rightarrow {\mathbb Z}\cup \{+\infty
\}$ denote the order function on $\widehat{K}^*$ extended by $ord(0)=+\infty$. 
In ${\mathbb R}^2$ one considers the convex hull of
$\bigcup _{i=0}^d \{(i,-ord(a_i)+x_2)|\ x_2\leq 0\}$. The finite part of the
boundary of this convex set is the Newton polygon of $L$.

 The module $M$ (over $K$ or over $\widehat{K}$) is called
{\it pure} if this Newton polygon has only one slope. As in the case of
differential operators, one can factorize $L$, viewed as an element of
$\widehat{K}[\Phi ,\Phi ^{-1}]$, according to the slopes in any order that one
chooses. This gives a unique decomposition of $\widehat{K}\otimes _KM$ as
direct sum $N_1\oplus N_2\oplus \cdots \oplus N_r$ of pure difference
modules over $\widehat{K}$ with slopes 
$\lambda _1<\lambda _2<\cdots <\lambda _r$.
The rule $\Phi z^n=q^nz^n\Phi$ and $|q|<1$ imply that the slope
factorization $L_1\cdot L_2\cdots L_r$ of $L$ where $L_i$ has slope 
$\lambda _i$ for $i=1,\dots ,r$ is convergent, i.e., all $L_i$ are in 
$K[\Phi ,\Phi ^{-1}]$. 

One deduces from this the ascending {\it slope filtration} of $M$ by submodules
$0=M_0\subset M_1\subset M_2\subset \cdots \subset M_r=M$ such that each 
$M_i/M_{i-1}$ is pure of slope $\lambda _i$ and moreover 
$\widehat{K}\otimes M_i/M_{i-1}\cong N_i$. The slope filtration is unique.
The {\it graded module} $gr(M)$ associated to $M$ is
$\oplus _{i=1}^r M_i/M_{i-1}$. 
We note that the above facts on
slope filtration are already present in the work of G.D.~Birkhoff, 
P.E.~Guenther and C.R.~Adams, see [Bir]. A modern proof is provided in [Sau3].
The difference module $M$ over $K$ is called {\it split} if
$M$ is isomorphic to $gr(M)$ (in other words, $M$ is a direct sum of
pure modules). Fix a direct sum $A=\oplus _{i=1}^rA_i$ of pure modules with
slopes $\lambda _1<\cdots <\lambda _r$.
In section 3 we will construct a fine moduli space for the
equivalence classes of the pairs $(M,f)$, consisting of a difference module
over $K$ and an isomorphism $f:gr(M)\rightarrow A$.

\subsection{Classification of pure modules over $\widehat{K}$ and $K$}
Let $F\subset G$ be a finite extension of difference fields. Let $M$ be a 
difference module over $G$. Then $Res(M)$ (the restriction of $M$ to $F$)
denotes $M$, considered as a  difference module over $F$. One observes that\\
$\dim _FRes(M)=[G:F]\cdot \dim _GM$.

Put $\widehat{K}_n=\widehat{K}(z^{1/n})$ for any integer $n\geq 1$. We apply
the above restriction to the extension $\widehat{K}\subset \widehat{K}_n$ in
order to construct all irreducible modules over $\widehat{K}$. Consider 
integers $t,n$ with $n\geq 1$ and $g.c.d.(t,n)=1$ and $c\in {\mathbb C}^*$ 
with $|q|^{1/n}<|c|\leq 1$. Let $E(cz^{t/n}):=\widehat{K}_ne$ denote the
difference module over $\widehat{K}_n$ given by $\Phi (e)=cz^{t/n}e$. 
Put $E:=Res(E(cz^{t/n}))$.
  
\begin{proposition}[The irreducible modules over $\widehat{K}$]$\;$\\ 
{\rm (1)} $E$ depends only on $t,\ n,\ c^n$.\\
{\rm (2)} $E$ is irreducible of dimension $n$ and has slope $t/n$. The algebra
of the $\widehat{K}$-linear endomorphisms of $E$, commuting with $\Phi$, is 
${\mathbb C}$.\\
{\rm (3)} For any irreducible difference module $I$ there are unique 
$t,n$ and $c^n$ with $n\geq 1$, $g.c.d.(t,n)=1$ and $|q|<|c^n|\leq 1$, such 
that $I\cong Res(E(cz^{t/n}))$.
\end{proposition}
\begin{proof} (1) and (2). 
$E$ has basis $e,\Phi e,\cdots ,\Phi ^{n-1}e$ over $\widehat{K}$ and
thus $e$ is a cyclic vector for $E$. The minimal monic polynomial 
$L\in \widehat{K}[\Phi ]$ with $Le=0$ is $L=\Phi ^n- q^{t(n-1)/2}c^nz^{t}$. 
Thus $E\cong \widehat{K}[\Phi ,\Phi ^{-1}]/\widehat{K}[\Phi ,\Phi ^{-1}]L$
and depends only on $t,\ n,\ c^n$. The operator $L$ has slope $t/n$ and degree
$n$. If $L$ has a non trivial decomposition $L_1L_2$, then the Newton polygon
of $L$ is the sum of the Newton polygons of $L_1$ and $L_2$. In particular
the Newton polygon of $L$ contains (at least) three points with integral
coordinates. Since $g.c.d.(t,n)=1$, this is not the case and hence $L$ 
and  $E$ are irreducible. We note that every non zero
endomorphism ($\widehat{K}$-linear and commuting with $\Phi$) of $E$
is bijective. Since $\mathbb C$ is algebraically closed,
this implies that the algebra of the endomorphisms ($\widehat{K}$-linear 
and commuting with $\Phi$) of $E$ is $\mathbb C$.\\
(3) Let $I$ be an irreducible difference module, then $I$ is pure and has
a slope $t/n$ with $n\geq 1$ and $g.c.d.(t,n)=1$. Take a cyclic vector and
let $L=\Phi ^d+a_1\Phi ^{d-1}+\cdots +a_{d-1}\Phi +a_d$ be its minimal
polynomial (with $d=\dim _{\widehat{K}}I$). Then 
$\frac{ord(a_i)}{i}\geq t/n$ for all
$i$ and $\frac{ord(a_d)}{d}=t/n$. It follows that $d$ is a multiple of $n$.
Now $L\in \widehat{K}[\Phi ,\Phi ^{-1}]\subset \widehat{K}_n[\Phi ,\Phi ^{-1}]
=\widehat{K}_n[\Psi, \Psi ^{-1}]$ with $\Psi =\Phi z^{t/n}$.\\
Then $L$, as
operator in $\Psi$, has slope 0. Hence $L$ has a right hand factor of degree 1
in $\Psi$ (or in $\Phi$). This means that we have a morphism of
$q$-difference modules over $\widehat{K}$
\[I=\widehat{K}[\Phi ,\Phi ^{-1}]/\widehat{K}[\Phi ,\Phi ^{-1}]L\rightarrow 
\widehat{K}_n[\Phi ,\Phi ^{-1}]/\widehat{K}_n[\Phi ,\Phi ^{-1}](\Phi -a)\ ,\]
for a suitable $a\in \widehat{K}_n$. This morphism is injective since
$I$ is irreducible. Counting the dimensions over $\widehat{K}$, one finds that
$d=n$ and that the morphism is bijective. The right hand side is a
one-dimensional difference module over $\widehat{K}_n$ and hence is 
isomorphic to $E(cz^{t/n})$ for some $c\in {\mathbb C}^*$ with
$|q|^{1/n}<|c|\leq 1$ (compare [vdP-S1], p. 149-150).

 \smallskip

We have to show that an isomorphism between 
$E_1:=Res(E(c_1z^{\lambda _1}))$ and $E_2:=Res(E(c_2z^{\lambda _2}))$ implies
that $\lambda _1=\lambda _2$ and $c_1^n=c_2^n$. The first statement is obvious
since $\lambda _i$ is the slope of $E_i$. We write $\lambda _1=\lambda _2=t/n$
with $n\geq 1,\ (t,n)=1$. Let 
$F:E_1\rightarrow E_2$ be an isomorphism. Then $F$ is unique up to 
multiplication by a scalar in ${\mathbb C}^*$. Both modules have the structure
of a difference module over $\widehat{K}_n$. Consider the map $z^{-1/n}\circ
F\circ z^{1/n}$. This is also an isomorphism between the two difference 
modules over $\widehat{K}$. Hence $z^{-1/n}\circ F\circ z^{1/n}=cF$ for some
$c\in {\mathbb C}^*$. Clearly $c^n=1$. We change the $\widehat{K}_n$ structure
of the module $E(c_2z^{t/n})$ by applying a suitable automorphism of 
$\widehat{K}_n$ over $\widehat{K}$. Now $E(c_2z^{t/n})$ is changed into
$E(c_3z^{t/n})$ with $c_3=\zeta c_2$ for some $\zeta $ with $\zeta ^n=1$.
Moreover, we have now $z^{-1/n}\circ F\circ z^{1/n}=F$. Thus
$E(c_1z^{t/n})$ and $E(c_3z^{t/n})$ are isomorphic as difference modules
over $\widehat{K}_n$. This implies $c_1=c_3$ since 
$|q|^{1/n}<|c_1|,|c_3|\leq 1$.   \end{proof}

\begin{remarks}$\;$\\ {\rm 
(1) We note that Proposition 1.3 extends to the case where the field
${\mathbb C}$ 
is replaced by any field $C$ (of characteristic 0, with $q\in C^*$ not a root 
of unity and $C$ not necessarily equal to its algebraic closure 
$\overline{C}$). This can be formulated
as follows. One extends the action of $\phi$ on $K:=C((z))$ to the field
$\bigcup \overline{C}((z^{1/n!}))$ in the obvious way. This field
contains the algebraic closure $\overline{K}$ of $K$. Take a non zero 
element $\alpha \in \overline{K}$ of degree $m$ over $K$ and consider the
difference module $K(\alpha )e$ given by $\Phi (e)=\alpha e$. Then 
$K(\alpha )e$, viewed as a difference module over $K$, has dimension $m$
and is irreducible. It depends only on the Galois orbit of $\alpha$. 
Every irreducible difference module over $K$ is obtained in this way.

\smallskip

\noindent (2)
Put $K_n=K(z^{1/n})$. The difference module over $K$ obtained by viewing 
$K_ne$ with $\Phi e=cz^{t/n}e$ as a difference module over $K$, will also 
be denoted by $Res(E(cz^{t/n}))$.
}\end{remarks}

\begin{corollary} Proposition {\rm 1.3} 
remains valid if $\widehat{K}$ is replaced by $K$.
\end{corollary}
\begin{proof} From the slope filtration it follows that an irreducible
difference module over $K$ is pure of some slope $t/n$.
Let $K_n$ denote $K(z^{1/n})$.
The factorization of $L$ as element of $K_n[\Psi ,\Psi ^{-1}]$ is valid
over $K_n$, because $L$ is in this context a regular singular 
difference operator. 
\end{proof}

\begin{corollary}[Indecomposable modules] $\;$\\
{\rm (1)} Let $M$ be an indecomposable difference module
over $\widehat{K}$. Then there are unique integers $t,\ n,\ m$
and $c^n\in {\mathbb C}^*$ with $n,m\geq 1$, $g.c.d.(t,n)=1$, 
$|q|<|c^n|\leq 1$ such that $M$ is isomorphic with 
$Res(E(cz^{t/n}))\otimes _{\widehat{K}}U_m$. \\
{\rm (2)} Let $M$ be an indecomposable pure difference module over $K$,
then there are unique $t,\ n,\ m,\ c^n$ 
as above such that $M\cong Res(E(cz^{t/n}))\otimes _KU_m$.
\end{corollary}
\begin{proof} (1) If the difference module $M$ over $\widehat{K}$ is 
indecomposable, then $M$ is pure. The proof of (2) that we will produce can be
copied verbatim to complete the proof of (1).\\ 
(2) Let $M/K$ be pure with slope $t/n$. We will concentrate on the
non trivial case where $n>1$. We consider now $K_n\otimes _KM$. This
difference module over $K_n$ has an action of the generator $\sigma$ of
the Galois group of $K_n/K$ defined by $\sigma z^{1/n}=\zeta z^{1/n}$ with
$\zeta =e^{2\pi i/n}$. One writes $K_n\otimes _KM$ as 
$K_ne\otimes _{\mathbb C}V$, where the action of $\Phi$ is given by
$\Phi (e\otimes v)=z^{t/n}e\otimes Av$ and where $A:V\rightarrow V$ is a
$\mathbb C$-linear map such that all its eigenvalues $\alpha$ satisfy
$|q|^{1/n}<|\alpha |\leq 1$. We note that this presentation of 
$K_n\otimes _KM$ is unique. Moreover, the subset 
${\mathbb C}[z^{1/n},z^{-1/n}]e\otimes V$ consists of the elements $f$ in
$K_n\otimes _KM$ such that the $\mathbb C$-vector space generated by
$\{(z^{-t/n}\Phi )^mf|\ m\in {\mathbb Z}\}$ has finite dimension. The vector
space $e\otimes V$ consists of the elements $f\in K_n\otimes _KM$ such that
there is a monic polynomial $L\in {\mathbb C}[(z^{-t/n}\Phi )]$ with 
$L(f)=0$ and all the roots $\alpha$ of $L$ satisfy
$|q|^{1/n}<|\alpha |\leq 1$. Since $\sigma$ commutes with $\Phi$ on
$K_n\otimes _KM$ one has that $e\otimes V$ is invariant under $\sigma$.
Hence we can write $\sigma (e\otimes v)=e\otimes B(v)$, where 
$B:V\rightarrow V$ is a linear map satisfying $B^n=1$. The fact that $\sigma$
and $\Phi$ commute translates into $BAB^{-1}=\zeta ^tA$. This induces a
decomposition $V=V_0\oplus V_1\oplus \cdots \oplus V_{n-1}$ into $A$-invariant
subspaces with the property $B(V_i)=V_{i+1}$ (where we use the cyclic notation
$V_n=V_0$). 

\smallskip
The submodules of $M$ are in bijection with the submodules of $K_n\otimes M$
that are invariant under $\sigma$. The latter are in bijection with
the $A$-invariant subspaces $W_0$ of $V_0$. This bijection 
associates to $W_0$ the $\sigma$-invariant submodule 
$K_ne\otimes _{\mathbb C}(\oplus _{i=0}^{n-1}B^iW_0)$. In particular,
$M$ is indecomposable if and only if the action of $A$ on $V_0$ has only one
Jordan block. Suppose that $A$ has this form and let $c$ be the eigenvalue of
$A$ on $V_0$, then one has $N\cong Res(E(cz^{t/n}))\otimes _K U_m$ with
$m=\dim V_0$. \end{proof}

We note that there are indecomposable difference modules over $K$ not
described in part (2) of Corollary 1.6.

\begin{corollary} Let $M$ be a pure difference module over $K$
with slope $\frac{t}{n}$ where $g.c.d.(t,n)=1$ and $n>1$. There exists a
difference module $N$ over $K_n$ such that $Res(N)\cong M$. The module $N$ is 
not unique. A similar statement holds for pure $q$-difference modules over
$\widehat{K}$.  
\end{corollary}

\begin{definition} $M_{global}$ for a pure module $M$ over $K$.\\{\rm
Suppose that the slope $\lambda$ of the pure module $M$ over $K$ is an integer.
Then $Kf\otimes _KM$, where $Kf$ is the module defined by $\Phi f=z^{-\lambda}
f$, is pure of slope 0. It follows that $M$ has a unique finite dimensional
$\mathbb C$-linear subspace $W$, such that $W$ is invariant under the operator
$z^{-\lambda}\Phi$ and the restriction $A\in {\rm GL}(W)$ has the property
that every eigenvalue $c$ of $A$ satisfies $|q|<|c|\leq 1$. Moreover, the
canonical $K$-linear map $K\otimes W\rightarrow M$ is a bijection.

For any $\mathbb C$-linear operator $L$ on $M$ and any $c\in \mathbb C$, one 
writes $Eigen(L,c)$ for the generalized eigenspace of $L$ for the eigenvalue
$c$. In other words $Eigen(L,c)=\bigcup _{s\geq 1}{\rm ker}( (L-c)^s,M)$. With
this terminology one has that 
\[W=\oplus _{c,\ |q|<|c|\leq 1}Eigen(z^{-\lambda}\Phi ,c)\ . \]
One defines $M_{global}:={\mathbb C}[z,z^{-1}]\otimes W=\oplus _{c\in 
{\mathbb C}^*} Eigen(z^{-\lambda}\Phi ,c)$. This is a free 
${\mathbb C}[z,z^{-1}]$-submodule of
$M$, invariant under $\Phi$ and $\Phi ^{-1}$. Thus $M_{global}$ is a global
difference module. Further, the canonical map
$K\otimes _{{\mathbb C}[z,z^{-1}]}M_{global}\rightarrow M$ is a bijection.

\smallskip
Now we consider a pure difference module $M$ with slope $\lambda =t/n$,
where $n\geq 1,\ (t,n)=1$. By Corollary 1.7, there exists a module $N$ over 
$K_n$ such that $M=Res(N)$. As above, $N$ has a unique finite dimensional
$\mathbb C$-linear subspace $W$ invariant under $z^{-\lambda }\Phi$, such that
all eigenvalues $c$ of  
the restriction $A$ of $z^{-\lambda}\Phi$ to $W$ satisfy $|q|^{1/n}<|c|\leq
1$. One defines $M_{global}:=N_{global}=
{\mathbb C}[z^{1/n},z^{-1/n}]\otimes W$.   
Thus $M_{global}=\oplus _{c\in {\mathbb C}^*}Eigen(z^{-\lambda}\Phi,c)$.
As before, $M_{global}$ is a global difference module and
the canonical map $K\otimes _{{\mathbb C}[z,z^{-1}]}M_{global}
\rightarrow M$ is an isomorphism.

In order to see that the definition of $M_{global}$ does not depend on the
choice of $N$ one considers the operator 
$(z^{-\lambda}\Phi)^n=q^\alpha z^{-t}\Phi ^n$, where $\alpha$ is some rational
number. It follows that $M_{global}$ is also equal to 
$\oplus _{c\in {\mathbb C}^*}Eigen(z^{-t}\Phi ^n,c)$. This expression is
clearly independent of the choice of $N$. Thus we can formulate the definition
of $M_{global}\subset M$ for a pure module over $K$ of slope $\lambda =t/n$
with $n\geq 1,\ g.c.d.(t,n)=1$ by the statement:\\
{\it The following properties of $m\in M$ are equivalent.\\
{\rm (1)} $m\in M_{global}$.\\ 
{\rm (2)} The $\mathbb C$-vector space generated by
$\{(z^{-t}\Phi ^n)^sm \ | s\geq 0 \}$ has finite dimension.\\ 
{\rm (3)}  There exists a $L\in {\mathbb C}[T],\ L\neq 0$ such that
$L(z^{-t}\Phi ^n)(m)=0$.} }\hfill $\square$\end{definition}

The main technical difficulties in this paper arise from pure modules $M$ with
non integer slope $\lambda =t/n$. There are two methods to handle these. 
The first one (Corollary 1.7) is to write $M=Res(N)$ for some difference module
$N$ over $K_n$. The second one, used in the proof of Corollary 1.6, replaces
$M$ by $K_n\otimes M$ provided with the action of the Galois group of
$K_n/K$. Both methods have their good and weak points. Now we develope the
second method in more detail. The main idea is to replace a pure differential
module $N$ over $K$ by $M=K_\infty \otimes _KN$ with decent data $D$. Here 
$K_\infty$ denotes the algebraic closure of $K$. With this method one can more
easily describe tensor products of pure modules over $K$.

\subsubsection{Pure difference modules over $K_\infty$ with descent data}
$K_\infty$ denotes the algebraic closure of $K$ and $Gal$ denotes
the Galois group of $K_\infty/K$. Let $M$ be a difference module over
$K_\infty$. {\it Descent data} $D$ for $M$ means a map 
$\sigma \in Gal \mapsto D(\sigma )$ satisfying:\\  
$D(\sigma )$ is a $\sigma$-linear bijection on $M$,
$D(\sigma )$ commutes with $\Phi$,\\
$D(\sigma _1)D(\sigma _2)=D(\sigma _1\sigma _2)$ and
the stabilizer of any $m\in M$, i.e., the group $\{\sigma \in Gal\ |\
D(\sigma )m=m\}$, is an open subgroup of $Gal$.

\medskip

One associates to a difference module $N$ over $K$ the module 
$M:=K_\infty \otimes N$ with descent data given by 
$D(\sigma )(f\otimes n)=\sigma (f)\otimes n$ for all $f\in K_\infty$ and
$n\in N$. This induces a functor from the category of the difference modules
over $K$ to the category of the difference modules over $K_\infty$ provided
with descent data.

\begin{proposition} $N\mapsto (K_\infty\otimes N,D)$ is an equivalence of
Tannakian categories.
\end{proposition}
\begin{proof} 
The essential thing to prove is that any pair $(M,D)$ is isomorphic to
$(K_\infty\otimes N,D)$ for some difference module $N$ over $K$. 

Take a basis $e_1,\dots ,e_r$ of $M$ over $K_\infty$. Let the
open subgroup $H:=\{\sigma \in Gal \ |\ \sigma (e_j)=e_j\mbox{ for all }j\}$
have index $m$ in $Gal$. Then $K_\infty ^H=K_m$ and 
$M^H=K_me_1+\cdots +K_me_r$. The cyclic group $Gal/H=Gal(K_m/K)$ acts on
$M^H$. This action induces an element of $H^1(Gal(K_m/K),{\rm GL}(r,K_m))$.
By Hilbert 90, this cohomology set is trivial. It follows that $M^H$
contains a basis $f_1,\dots ,f_r$ over $K_m$, consisting of $Gal$-invariant
elements. Now $N:=Kf_1\oplus \cdots \oplus Kf_r$ is equal to $M^{Gal}$ and
the natural map $K_\infty \otimes _KN\rightarrow M$ is an isomorphism.
Since $\Phi$ commutes with the action of $Gal$, one has $\Phi (N)=N$. Thus
$N$ is a difference module over $K$ and clearly induces the pair
$(M,D)$.

The tensor product of two pairs $(M_1,D_1),\ (M_2,D_2)$ is defined as
$(M_1\otimes _{K_\infty}M_2,D_1\otimes D_2)$. We note that the tensor product
$D_1(\sigma )\otimes D_2(\sigma )$ of two $\sigma$-linear maps 
makes sense. It is easily seen that the above equivalence respects tensor  
products.  \end{proof}

For a pure difference module $N$ over $K$ of slope $\lambda$, the module
$M=K_\infty \otimes N$ is also pure with slope $\lambda$ and has the form
$K_\infty \otimes _{\mathbb C}V$, where $V$ is a finite dimensional
${\mathbb C}$-vector space provided with an element $A\in {\rm GL}(V)$. 
The action of $\Phi$ on $M$ is given by 
$\Phi (f\otimes v)=z^\lambda \phi (f)\otimes A(v)$.

The subspace $V$ is not unique. By changing $V$, the eigenvalues
of $A$ are multiplied by arbitrary, rational powers of $q$. We
normalize $A$ and $V$ as follows.

Choose a ${\mathbb Q}$-linear subspace $L\subset {\mathbb C}$ such that
$L\oplus {\mathbb Q}={\mathbb C}$. One requires that every eigenvalue $c$ of 
$A$ has the form $e^{2\pi i(a_0(c)+a_1(c)\tau )}$
with $a_0(c),a_1(c)\in {\mathbb R}$ and $a_1(c)\in L$.

After this normalization the subspace $V$ of $M$ is
unique. Indeed, $V$ is the direct sum of the kernels of 
$(z^{-\lambda}\Phi -c)^s$ with $s>>0$ and $c\in {\mathbb C}^*$ with 
$a_1(c)\in L$.

We note that the use of this subspace $L$ is somewhat
artificial. It can be avoided at the cost of verifying that 
formulas that we will produce are independent of certain choices.

One observes that $V$ is invariant under $D(\sigma )$ for all 
$\sigma \in Gal$. The group $Gal$ is identified with $\widehat{\mathbb Z}$ 
and the
action of $Gal$ is expressed by $\sigma (z^\lambda )=
e^{2\pi i\lambda \sigma} z^\lambda$. For the operators $A$ and $D(\sigma )$,
restricted to $V$, one finds the equality
$AD(\sigma )=e^{2\pi i\lambda \sigma}D(\sigma )A.$ 
Thus we have associated to a pure difference module $N$ over $K$ 
a tuple $data(N):=(\lambda ,V,A,\{D(\sigma )\})$ with
\begin{itemize}
\item $\lambda \in {\mathbb Q}$. 
\item $V$ a vector space over $\mathbb C$ of finite dimension.
\item $A\in {\rm GL}(V)$ with eigenvalues in the subgroup \\
$\{c=e^{2\pi i(a_0(c)+a_1(c)\tau )}|a_0(c)\in {\mathbb R},\ a_1(c)\in L\}$ of
${\mathbb C}^*$.
\item a homomorphism $\sigma \in Gal\cong \widehat{\mathbb Z}
\mapsto D(\sigma )\in {\rm  GL}(V)$ satisfying\\ 
$AD(\sigma )=e^{2\pi i\lambda \sigma}D(\sigma )A$.
\end{itemize}

On the other hand an object $(\lambda ,V,A,\{D(\sigma )\})$ as above defines
a pure module $N$ over $K$ of slope $\lambda$ in the following way.
Consider $M:=K_\infty\otimes V$ with $\Phi$ given by 
$\Phi (f\otimes v)=z^\lambda \phi (f)\otimes A(v)$ and
with descent data given by $D(\sigma )(f\otimes v)=\sigma (f)\otimes 
D(\sigma )v$. Then $N:=M^{Gal}$.

Consider a morphism  $f:N_1\rightarrow N_2,\ f\neq 0$ between pure modules. 
Then $N_1,N_2$ have the same slope $\lambda$ and $f$ induces a morphism from
$data(N_1)$ to $data(N_2)$, i.e., a linear map $F$ between the two 
$\mathbb C$-vector spaces equivariant for the maps of the data. On the other
hand, a $\mathbb C$-linear map $F$, equivariant for the maps of the data, comes
from a unique morphism $f:N_1\rightarrow N_2$. 

{\it Thus $N\mapsto data(N)$ is an equivalence between the category of the pure
modules over $K$ and the category of tuples $(\lambda ,V,A, \{D(\sigma )\})$
defined above}. One observes the following useful properties.

\smallskip
For pure difference modules $N_i$ with
$data(N_i)=(\lambda _i,V_i,A_i,\{D_i(\sigma )\})$
for $i=1,2$ one has the nice formula \\
$data(N_1\otimes N_2)=(\lambda _1+\lambda _2,V_1\otimes V_2,
A_1\otimes A_2,\{D_1(\sigma )\otimes D_2(\sigma )\})$.

\smallskip
Let the pure modules $N$ have $data(N)=(\lambda ,V,A,\{D(\sigma )\})$.
Then the dual module $N^*$ has data $(-\lambda ,V^*,B,\{E(\sigma )\})$, where
$V^*$ is the dual of $V$; $B=(A^{-1})^*$ and $E(\sigma )=(D(\sigma )^{-1})^*$.

The $\Phi$-equivariant pairing $N\times N^*\rightarrow K$, given by
$(n,\ell )\mapsto \ell (n)\in K$ translates for the data of $N$ and $N^*$ 
into the usual pairing $V\times V^*\rightarrow {\mathbb C}$, given by
$(v,\ell )\mapsto \ell (v)\in {\mathbb C}$. This pairing is equivariant
with respect to the prescribed actions on $V$ and $V^*$.

\section{Vector bundles and $q$-difference modules}

We recall that $O$ denotes the algebra of the holomorphic functions on 
${\mathbb C}^*$ and that a difference module $M$ over $O$ is a left module 
over the ring $O[\Phi ,\Phi ^{-1}]$, free of some rank $m<\infty $ over 
$O$. Further $pr:{\mathbb C}^*\rightarrow E_q:={\mathbb C}^*/q^{\mathbb Z}$
denotes the canonical map. One associates 
to $M$ the vector bundle $v(M)$ of rank $m$ on $E_q$ given by
$v(M)(U)=\{f\in O(pr ^{-1}U)\otimes _OM|\ \Phi (f)=f\}$,
where, for any open $V\subset {\mathbb C}^*$, $O(V)$ is the algebra of the
holomorphic functions on $V$.

On the other hand, let a vector bundle $\mathcal M$ of rank
$m$ on $E_q$ be given. Then ${\mathcal N}:=pr^*{\mathcal M}$ is a vector
bundle on ${\mathbb C}^*$ provided with a natural isomorphism
$\sigma _q^*{\mathcal N}\rightarrow {\mathcal N}$, where $\sigma _q$ is the
map $\sigma _q(z)=qz$. One knows that ${\mathcal N}$ is in fact a free
(or trivial) vector bundle of rank $m$ on ${\mathbb C}^*$ (see [For], p. 204).
 Therefore, $M$, the collection
of the global sections of $\mathcal N$, is a free $O$-module
of rank $m$ provided with an invertible action $\Phi$ satisfying
$\Phi (fm)=\phi (f)\Phi (m)$ for $f\in O$ and $m\in M$.
It is easily verified that the above describes an equivalence $v$ of tensor
categories. 

\smallskip

The equivalence $v$ extends to an equivalence between the 
left $O[\Phi, \Phi ^{-1}]$-modules which are finitely generated
as $O$-module and the coherent sheaves on $E_q$. This is an
equivalence of Tannakian categories.

By an {\it admissible difference module over $O$} we will mean a left
$O[\Phi, \Phi ^{-1}]$-module which is a
direct limit of left $O[\Phi, \Phi ^{-1}]$-modules of finite
type over $O$. The equivalence $v$ extends to a Tannakian equivalence
between category of the admissible difference modules over $O$ and the 
category of the quasi-coherent sheaves on $E_q$.

\begin{lemma} There are isomorphisms 
${\rm ker}(\Phi -1,M)\rightarrow H^0(E_q,v(M))$ and 
${\rm coker}(\Phi -1,M)\rightarrow H^1(E_q,v(M))$ between these functors
defined on the category of the admissible difference modules $M$ over $O$.
\end{lemma}
\begin{proof} The isomorphism ${\rm ker}(\Phi -1,M)\rightarrow H^0(E_q,v(M))$
follows from the definition of $v$. Let $v^{-1}$ denote the `inverse'
of the functor $v$. Then ${\mathcal M}\mapsto {\rm ker}(\Phi
-1,v^{-1}({\mathcal M}))$ is canonically isomorphic to 
${\mathcal M}\mapsto H^0(E_q,{\mathcal M})$.
One observes that an exact sequence
$0\rightarrow M_1\rightarrow M_2\rightarrow M_3\rightarrow 0$
 of admissible difference modules over $O$ induces (by the snake lemma) an 
exact sequence
\[0\rightarrow {\rm ker}(\Phi-1,M_1)\rightarrow {\rm ker}(\Phi-1,M_2)
\rightarrow {\rm ker}(\Phi-1,M_3)\rightarrow \]\[{\rm coker}(\Phi-1,M_1)
\rightarrow {\rm coker}(\Phi-1,M_2)\rightarrow {\rm coker}(\Phi-1,M_3)
\rightarrow 0 \ .\]
>From this it easily follows that the first right derived functor of the 
functor $M\mapsto {\rm ker}(\Phi -1,M)$, on the category of admissible modules
over $O$, is equal to ${\rm coker}(\Phi -1,M)$. Now the second isomorphism of
functors follows. \end{proof}

\begin{examples}$\;$\\ {\rm
(1) Consider $M=Oe$ with $\Phi (e)=e$. Then $v(M)$ is the structure sheaf
$O_{E_q}$ of $E_q$. Any element in $m\in M$ can be written uniquely as 
$m=\sum _{n\in {\mathbb Z}}a_nz^ne$. Then 
$(\Phi -1)m=\sum _{n\in {\mathbb Z}}(q^n-1)a_nz^ne$. One observes that
${\rm ker}(\Phi -1,M)={\mathbb C}e$ and that ${\rm coker}(\Phi -1,M)$
is represented by ${\mathbb C}e$. This illustrates Lemma 2.1.

\smallskip

\noindent
(2) Consider difference module $M=Oe$ with $\Phi e=ce$ and 
$c\in {\mathbb C}^*,\  |q|<|c|\leq 1$. If $c\neq 1$, then $\theta _ce$ is a 
meromorphic section of $v(M)$ with divisor  $-pr(c^{-1})+pr(1)$. One concludes
that $v(M)\cong O_{E_q}(pr(c^{-1})-pr(1))$. Thus one finds all line bundles of
degree 0 on $E_q$ in this way.

\smallskip

\noindent (3)
Consider the difference module $M:=Oe$ with $\Phi e=(-z)e$. There is a
$\Phi$-invariant element, namely $\Theta e$. This is a global section of
$v(M)$. The cokernel of the morphism $O_{E_q}\rightarrow v(M)$, given by
$1\mapsto \Theta e$, is a skyscraper sheaf with support $\{1\}$ and stalk
$\mathbb C$ at that point. Indeed, the function $\Theta$ has simple zeros at 
$q^{\mathbb Z}$. One concludes that $v(M)\cong O_{E_q}([1])$. Using tensor
products one obtains that the line bundle $v(M)$, with $M=Oe,\ \Phi e=cz^te$,
has degree $t$. Moreover every line bundle on $E_q$ is obtained in this way.

\smallskip

\noindent (4) Let $M=O/O\Theta$ with the $\Phi$-action induced by the 
usual one of $O$. Then $v(M)$ is the skyscraper sheaf on $E_q$ with support
$\{1\}$ and with stalk $\mathbb C$ at that point.\hfill $\square$
}\end{examples}

We recall that a {\it split difference module $M$ over $K$} is a direct sum
of pure modules $M_i$. The global module $M_{global}$ over ${\mathbb
  C}[z,z^{-1}]$ associated to $M$ is by definition the direct sum of the 
global modules $(M_i)_{global}$. A morphism $f:M\rightarrow N$ between split
modules is easily seen to be the direct sum of morphisms between the pure
components of $M$ and $N$. In particular $f$ maps $M_{global}$ to
$N_{global}$. 

One associates to a split difference module $M$ the difference module 
$O\otimes _{{\mathbb C}[z,z^{-1}]}M_{global}$ and, by Lemma 2.1,  a vector 
bundle on $E_q$. For notational convenience we write again $v(M)$ for this
vector bundle. In this way we obtain a functor $v$ from the category of the
split difference modules over $K$ to the category of vector bundles on $E_q$.
One observes that ${\rm Hom}(M_1,M_2)\rightarrow {\rm Hom}(v(M_1),
v(M_2))$ is $\mathbb C$-linear and injective. Moreover, one easily sees that 
$v$ preserves tensor products.

\begin{theorem} The functor $v$ from the category of the split
difference modules over $K$ to the category of the vector bundles on
$E_q$ is bijective on isomorphy classes of objects. This bijection respects
tensor products. 
\end{theorem}
\begin{proof} We have to show that $v$ induces a bijection between the
isomorphy classes of the indecomposable objects in the two categories.
We start by proving that for an indecomposable pure difference module
$M$ the corresponding vector bundle $v(M)$ is indecomposable.

\smallskip

>From Examples 2.2 one concludes that $v$ provides a bijection between the
isomorphy classes of the difference modules of dimension 1 over $K$ and the 
isomorphy classes of all line bundles on $E_q$.  

\smallskip

By Corollary 1.6, an indecomposable pure difference module has 
the form $M=Res(E(cz^{t/n}))\otimes U_m$, with unique 
$n\geq 1,\ (t,n)=1,\ m\geq 1$ and $c^n$ such that $|q|<|c^n|\leq 1$. The 
vector bundle $v(M)$ has clearly rank $nm$. The 
exterior product $\Lambda ^{nm}M$ is equal to $Kf$ with $\Phi f=sz^{tm}f$ for
some $s\in {\mathbb C}^*$. Thus $v(M)$ has degree $tm$. The case $nm=1$ has
been treated above and we suppose now $nm>1$.

One can present $M_{global}$ as ${\mathbb C}[z^{1/n},z^{-1/n}]\otimes W$,
with $W$ a $\mathbb C$-linear space of dimension $m$ and $\Phi$ given by
$\Phi (1\otimes w)=cz^{t/n}\otimes U(w)$ where $U$ is a unipotent map with
minimal polynomial $(U-1)^m=0$. Then 
$O\otimes _{{\mathbb C}[z,z^{-1}]}M_{global}$ can
be represented as $H_n\otimes W$, where $H_n$ consists of the convergent
Laurent series in $z^{1/n}$. Thus $H_n$ consists of the expressions 
$\sum _{k=-\infty}^{+\infty}a_kz^{k/n}$ with 
$\lim _{|k|\rightarrow \infty}|a_k|^{1/|k|}=0$. One provides $W$ with some
norm $\|\ \|$. The elements of $H_n\otimes W$ have the form
$\sum _{k=-\infty}^{+\infty}z^{k/n}\otimes w_k$ with 
$\lim _{|k|\rightarrow \infty}\|w_k\|^{1/|k|}=0$.
Then $\Phi$ acts on $H_n\otimes W$ by 
\[\Phi (\sum z^{k/n}\otimes w_k)=
\sum q^{k/n}z^{k/n}cz^{t/n}\otimes U(w_k) \ .\] 
Write $\Psi = z^{-t}\Phi ^n$. Then 
\[\Psi (\sum z^{k/n}\otimes w_k)=
\sum z^{k/n}\otimes dq^kU^n(w_k)\mbox{, with }d=c^nq^{t(n-1)/2}\ .\]
For each
$k$, the vector space $z^{k/n}\otimes W$ is invariant under $\Psi$ and this
operator has eigenvalues $q^kd$ on this vector space. One concludes from this
that the subset $M_{global}\subset H_n\otimes W$ consists of the elements
$fe$ such that there exists a non zero polynomial $L\in {\mathbb C}[T]$ with
$L(\Psi) (fe)=0$. This has as consequence that every 
$O$-linear endomorphism $A$ of 
$O\otimes _{{\mathbb C}[z,z^{-1}]}M_{global}$, commuting with $\Phi$, is the 
$O$-linear extension of a unique ${\mathbb C}[z,z^{-1}]$-linear
endomorphism $B$ of $M_{global}$ commuting with $\Phi$. 

A direct sum decomposition of $v(M)$ induces a $O$-linear
endomorphism $A$ of $O\otimes _{{\mathbb C}[z,z^{-1}]} M_{global}$ 
commuting with
$\Phi$ and such that $A^2=A$. The corresponding $B$ induces a direct sum
decomposition of $M_{global}$, contradicting that $M$ is indecomposable.
Thus $v(M)$ is indecomposable. A similar reasoning proves that for
indecomposable $M_1,M_2$ the relation $v(M_1)\cong v(M_2)$ implies that
$M_1\cong M_2$. In this way we have found a collection of indecomposable
vector bundles on $E_q$. That we have found all of them follows at once from
the classification given in [At], Theorem 10. Indeed, Atiyah 
constructs a certain indecomposable vector bundle of rank $r$ and degree
$d$, called $E_A(r,d)$. Let $h=(r,d)$. Then every indecomposable vector bundle
of rank $r$ and degree $d$ has the form $L\otimes E_A(r,d)$ with $L$ a line
bundle of degree 0. This $L$ is unique up to multiplication with a line bundle
$N$ such that $N^{\otimes r/h}$ is the trivial line bundle.   
\end{proof}

The final part of the proof of Theorem 2.3 depends on [At]. We present now
a proof which only uses a simple result of this paper, namely Lemma 11,
formulated as follows:\\
{\it Let $W$ be an indecomposable vector bundle of rank $m$ and degree 0
on an elliptic curve $E$, then $W=L\otimes W'$, with $L$ a line bundle of
degree 0 and such that the indecomposable $W'$ has a sequence of subbundles
$0=W'_0\subset W'_1\subset \cdots \subset W'_m=W'$ such that each quotient
$W'_{i+1}/W'_i$ is isomorphic to $O_E$}.

\begin{proof} 
$V$ is an indecomposable vector bundle on $E_q$, rank $nm$ and degree $tm$
with $n,m\geq 1,\ g.c.d.(t,n)=1$. As before we consider
$pr:{\mathbb C}^*_z\rightarrow E_q={\mathbb C}_z^*/q^{\mathbb Z}$. The index
$z$ means that we use $z$ as variable on this copy of ${\mathbb C}^*$.
Write $M:=H^0({\mathbb C}^*,pr^*(V))$. It suffices to produce a 
${\mathbb C}[z,z^{-1}]$-submodule $M_0\subset M$, invariant under $\Phi$ and
$\Phi ^{-1}$, such that the natural map
$O\otimes _{{\mathbb C}[z,z^{-1}]}M_0\rightarrow M$ is bijective and
$K\otimes _{{\mathbb C}[z,z^{-1}]}M_0$ is a pure module.

Let $\beta :{\mathbb C}^*_s\rightarrow {\mathbb C}_z^*$ be given by 
$s\mapsto s^n=z$ (or by $s=z^{1/n}$). Define the elliptic curve $E$ by
$\gamma :{\mathbb C}^*_s\rightarrow E:={\mathbb C}^*_s/(q^{1/n})^{\mathbb Z}$.
There is an induced morphism of $\alpha :E\rightarrow E_q$, of degree $n$, 
such that $\alpha \circ \gamma = pr\circ \beta$.\\
The map $\alpha$ is an unramified cyclic covering of degree $n$. Let $\sigma$
denote a generator of the automorphism group of this covering. One can take
for $\sigma$ the automorphism of ${\mathbb C}^*_s\rightarrow {\mathbb C}^*_z$,
given by $s\mapsto e^{2\pi i/n}s$. The last map will also be called 
$\sigma$. 

\smallskip

The vector bundle $\alpha ^*V$ on $E$ has rank $nm$ and degree $tnm$. Then
$\alpha ^*V$ is a direct sum of indecomposable vector bundles 
$W_1\oplus \cdots \oplus W_r$ on $E$. Let $W_1,W_2,\dots ,W_{r'}$ with
$r'\leq r$ be all the $W_i$ which have the same rank and degree as $W_1$.
The direct sum $W_1\oplus \cdots \oplus W_{r'}$ is invariant under the
action of $\sigma ^*$. Since $V$ is indecomposable, one has $r'=r$. Thus
all $W_i$ have the same rank and degree. It follows that 
$O_E(-t[1_E])\otimes \alpha ^*V$ is the direct sum of
indecomposable vector bundles $W_i'$ on $E$ of degree 0 and rank $nm/r$.

\smallskip

Using Lemma 11 of [At], we conclude that the
difference module $H:=H^0({\mathbb C}^*_s,\gamma ^*W'_i)$ over
$O({\mathbb C}^*_s)$ has a sequence of submodules $0=H_0\subset \cdots \subset
H_{nm/r}=H$ such that each quotient has the form $O({\mathbb C}^*_s)e$ with 
$\Phi e=e$. Thus $H$ has a basis $e_1,\dots ,e_{mn/r}$
over $O({\mathbb C}^*_s)$ such that the matrix of $\Phi$ w.r.t. this basis is 
upper triangular and all its diagonal entries are 1. One easily verifies that 
a base change turns this matrix into a matrix with constant coefficients.

\smallskip

The difference module over $O({\mathbb C}^*_s)$ associated to $O_E(t[1_E])$
has the form $O({\mathbb C}^*_s)e$ with $\Phi e=cs^{t}e$ for some constant
$c$. By taking the tensor product with this module and direct sums, we 
conclude that the difference 
module $N$ over $O({\mathbb C}^*_s)$, 
associated to $\gamma ^*\alpha ^* V=\beta ^*pr^*V$, 
has a basis for which the matrices of $\Phi$ and $\Phi ^{-1}$ have 
coefficients in ${\mathbb C}\cdot s^{\mathbb Z}$. In particular, there is 
a ${\mathbb C}[s,s^{-1}]$-module $N_0$ of $N$, invariant under
$\Phi$ and $\Phi ^{-1}$, such that
$O({\mathbb C}^*_s)\otimes _{{\mathbb C}[s,s^{-1}]}N_0=N$. This submodule is,
by construction, invariant under the action of $\sigma$. Therefore the 
${\mathbb C}[z,z^{-1}]$-module $M_0:=N_0^\sigma$ is finitely generated and
invariant under the actions of $\Phi$ and $\Phi ^{-1}$. Moreover 
$O({\mathbb C}^*_z)\otimes _{{\mathbb C}[z,z^{-1}]}M_0=N^\sigma =M$ and
$K\otimes _{{\mathbb C}[z,z ^{-1}]}M_0$ is pure.  
\end{proof}

\begin{remarks} Tensor products of pure difference modules over $K$.\\
{\rm One rediscovers Part III of [At] (for the base field $\mathbb C$) by 
using Theorem 2.3 and some calculations for difference modules. We will
give some results.

\smallskip
\noindent 
(1) The indecomposable vector bundle corresponding to $U_m$ is called 
$F_m$ in [At]. The tensor product $U_a\otimes U_b$ corresponds to the
tensor product of two unipotent operators $A,B$ on vector spaces
$V,W$ of dimensions $a,b$, having each only one Jordan block. One can find
a decomposition of the unipotent operator $A\otimes B\in {\rm GL}(V\otimes W)$
as a direct sum of Jordan blocks. This will produce the formulas in [At], 
Theorem 8 and our method is close to the remarks [At], p. 438-439. 

\smallskip
\noindent
(2) Consider the module $(K_ne,\Phi e=cz^{t_1/n_1}e)$ with $(t_1,n_1)=1$ and
$n=dn_1$ with $d>1$. We note that $c$ is determined up to an element of
$\mu _{n_1}\times q^{{\mathbb Z}/n}$, where $\mu _k$ denotes the group of the
$k$th roots of unity. This module decomposes as the direct sum of the
irreducible modules $(K_{n_1}e_j,\Phi e_j=q^{j/n}cz^{t_1/n_1}e_j)$ for
$j=0,\dots ,d-1$. A change of $c$ will permute these modules. 

\smallskip
\noindent
(3) $M_i:=(K_{n_i}e_i, \Phi e_i=c_iz^{t_i/n_i}e_i)$ for $i=1,2$ with
$(t_i,n_i)=1$. Let $d=(n_1,n_2)$ and $n=\frac{n_1n_2}{d}$. We note that
$c_i$ is determined up to an element in 
$\mu _{n_i}\times q^{{\mathbb Z}/n_i}$. Write $t_1/n_1+t_2/n_2=t_3/n_3$
with $(t_3,n_3)=1$ and $n=kn_3$.
Then $M_1\otimes M_2$ is the direct sum of $(K_nf_j,\Phi f_j
=\zeta _jc_1c_2z^{t_3/n_3}f_j)$ 
for $j=0,\dots ,d-1$ and suitable $n$th roots of
unity $\zeta _j$. Each term has a further decomposition, according to (2), if
$k>1$. 

\smallskip
\noindent
(4) $M=(K_ne,\Phi e=cz^{t/n}e)$ with $(t,n)=1$. The dual $M^*$ is equal to
$(K_ne^*,\Phi e^*=c^{-1}z^{-t/n}e^*)$.
Further $M\otimes M^*$ is the direct sum of the $n^2$ difference modules
$(Ke_{s,t},\Phi e_{s,t}=q^{s/n}e^{2\pi it/n}e_{s,t})$ with
$0\leq s,t<n$. This corresponds to the
direct sum of all line bundles $L$ of order dividing $n$. }\end{remarks}

\smallskip
\noindent
(5) The homomorphism
${\rm Hom}(M_1,M_2)\rightarrow {\rm Hom}(v(M_1),v(M_2))$, where $M_1,M_2$
are split modules over $K$, is in general not surjective. Indeed, the category
of the split difference modules over $K$ is Tannakian and the category of the
vector bundles on $E_q$ is not even an abelian category.

\begin{theorem} $\;$\\
Let $M$ be a pure difference module over $K$ with slope $\lambda <0$.
 The maps\\
{\rm (1)} ${\rm coker}(\Phi -1,M_{global})\rightarrow {\rm coker}
(\Phi -1,O\otimes _{{\mathbb C}[z,z^{-1}]} M_{global})$ \\
{\rm (2)} ${\rm coker}(\Phi -1,M_{global})\rightarrow 
{\rm coker}(\Phi -1,K\otimes _{{\mathbb C}[z,z^{-1}]}M_{global})$\\ 
are isomorphisms. Moreover, ${\rm coker}(\Phi -1,M_{global})$ is canonical
isomorphic to $H^1(E_q,v(M))$ and has dimension $-\lambda \cdot dim _KM$
over $\mathbb C$.
\end{theorem}

Before starting the technical proof we study the simplest example:\\ 
\[M_{global}={\mathbb C}[z,z^{-1}]e\mbox{ with }\Phi e=z^{-1}e \ .\] 
An element of $O\otimes M_{global}=Oe$ has the form $\sum a_nz^ne$. By
Examples 2.2, $v(O\otimes M_{global})=O_{E_q}(-[1])$ and by 
Lemma 2.1, one has 
${\rm ker}(\Phi -1,O\otimes M_{global})=H^0(E_q,O_{E_q}(-[1]))=0$.
Consider the equation 
\[(\Phi -1)\sum a_nz^ne=\sum (a_{n+1}q^{n+1}-a_n)z^ne=\sum b_nz^ne \ ,\]
for a given $\sum b_nz^ne\in Oe$. There is at most one solution. If there is a
solution then its coefficients satisfy the recurrence 
$a_{n+1}q^{n+1}-a_n=b_n$ which can also be written as
\[ a_{n+1}q^{(n+2)(n+1)/2}-a_nq^{(n+1)n/2}=b_nq^{(n+1)n/2} \ .\]
Summation of these recurrences
over all $n\in \mathbb Z$ implies $\sum b_nq^{(n+1)n/2}=0$.

On the other hand, the condition $\sum b_nq^{(n+1)n/2}=0$ leads to a formula 
$a_n:=\sum _{s=n}^\infty -b_sq^{(s+1)s/2-(n+1)n/2}$. One easily verifies that
this infinite sum converges, that the series $\sum a_nz^ne$ belongs to $Oe$
and solves the equation.
If, moreover, $\sum b_nz^ne\in {\mathbb C}[z,z^{-1}]e=M_{global}$, then also
$\sum a_nz^ne$ lies in $M_{global}$. This proves part (1) of Theorem 2.5 for 
this example.

\begin{proof}  (1) Let $-t/n <0$, with $t,n\geq 1,\ g.c.d.(t,n)=1$, be the 
slope of $M$. Then $M_{global}$ can be
written as ${\mathbb C}[z^{1/n},z^{-1/n}]\otimes _{\mathbb C}W$ with $\Phi$
action given by $\Phi (f\otimes w)=z^{-t/n}\phi (f)\otimes A(w)$ for some 
$A\in {\rm GL}(W)$. Then $O\otimes M_{global}$ can be written as
$O_n\otimes _{\mathbb C}W$, where $O_n$ consists of the convergent Laurent
series in $z^{1/n}$. The action of $\Phi$ is again given by 
$\Phi (f\otimes w)=z^{-t/n}\phi (f)\otimes A(w)$. Consider the equation 
\[(\Phi -1)(\sum _{k\in \mathbb Z}z^{k/n}\otimes x_k)= 
\sum _{k\in \mathbb Z}z^{k/n}\otimes (q^{(k+t)/n}A(x_{k+t})-x_k)=
\sum _{k\in \mathbb Z}z^{k/n}\otimes w_k ,\]
where the given series $\sum z^{k/n}\otimes w_k$ is either finite or 
convergent. This yields recurrence relations for the $x_k$. Write
$k=k_0+st$ with $k_0\in \{0,1,\dots ,t-1\}$ and $s\in \mathbb Z$.
The recurrence relations are
\[q^{(k_0+(s+1)t)/n}A(x_{k_0+(s+1)t})-x_{k_0+st}=w_{k_0+st} \mbox{ for } 
k_0\in \{0,1,\dots ,t-1\},\  s\in {\mathbb Z}\ .\] 
Put $m(k_0,s)=\frac{sk_0}{n}+\frac{s(s+1)t}{2n}$. Then the recurrences
can be rewritten as
\[q^{m(k_0,s+1)}A^{s+1}(x_{k_0+(s+1)t})-q^{m(k_0,s)}A^s(x_{k_0+st})=
q^{m(k_0,s)}A^s(w_{k_0+st}) \ .\]
Suppose that there exists a convergent solution $\sum z^{k/n}\otimes x_k$,
then, for each $k_0$, the sum over all $s\in \mathbb Z$ of the left hand side
is 0. Thus a necessary condition for the existence of a convergent 
solution is
\[ \sum _{s\in \mathbb Z}q^{m(k_0,s)}A^s(w_{k_0+st})=0\mbox{ for }
k_0=0,1,\dots ,t-1\ .\]

This condition is also sufficient for the existence of a convergent solution.
Indeed, the recurrences and the above condition on the  coefficients
$\{w_{k_0+st}\}$ lead to the formula 
\[x_{k_0+st}=\sum _{a=s}^\infty -q^{m(k_0,a)-m(k_0,s)}A^{a-s}(w_{k_0+st}) \ .\]
The possibly infinite sum in this formula converges and one can verify that
the resulting series $\sum z^{k/n}\otimes x_k$ lies in $O_n\otimes W$ and
solves the above equation.

If there are ony finitely many $w_k\neq 0$ and the above condition is 
satisfied, then the unique solution $\sum _kz^k\otimes x_k$ has only finitely
many $x_k\neq 0$. From this observation statement (1) follows.\\ 
The proof of (2) is similar. Lemma 2.1 and the above computation yield a proof
of the last statement of the theorem. \end{proof}

\begin{leeg} Canonical representatives for ${\rm coker}(\Phi
  -1,M_{global})$.\\
{\rm Let $M$ be a pure difference module over $K$ with slope $-t/n$
and $n,t\geq 1,\ g.c.d.(t,n)=1$. As in the proof of Theorem 2.5, one writes 
$M_{global}={\mathbb C}[z^{1/n},z^{-1/n}]\otimes W$ with $\Phi$ action
given by $\Phi (f\otimes w)=z^{-t/n}\phi (f)\otimes A(w)$. One requires
that $|q|^{1/n}<|c|\leq 1$ for every eigenvalue of $A$. This makes the
presentation unique.

Define 
$W^+:=({\mathbb C}1+{\mathbb C}z^{1/n}+\cdots +{\mathbb C}z^{(t-1)/n})\otimes
W$. We claim that
\[W^+\rightarrow {\rm coker}(\Phi -1,M_{global})\mbox{ is a bijection.}\]
Indeed, write $w^+\in W^+$ as $w^+=\sum _{i=0}^{t-1}z^{i/n}\otimes w_i$.
Then $w^+$ lies in the image of $\Phi -1$ if and only if for all
$k_0\in \{0,1,\dots ,t-1\}$ one has 
$\sum _{s\in \mathbb Z}q^{m(k_0,s)}w_{k_0+st}=0$. This implies that $w^+=0$.
Hence $W^+\rightarrow {\rm coker}(\Phi -1,M_{global})$ is injective. Since
both spaces have dimension $t\cdot \dim W$, the map is a bijection.

For any $\mathbb C$-linear operator $T$ on $M$, one writes $Eigen(T,c)$
for the generalized eigenspace of $T$ for the eigenvalue $c\in {\mathbb C}$.
Thus $Eigen(T,c)=\bigcup _{s\geq 1}{\rm ker}((T-c)^s,M)$. With this notation
one observes that 
\[W^+=\oplus _{|q|^{t/n}<|c|\leq 1}Eigen(z^{t/n}\Phi ,c)\mbox{ and }
 W^+=\oplus _{|q|^{t}<|c|\leq 1}Eigen((z^{t/n}\Phi )^n,c)\ .\]

We write $Repr(M_{global})$ or $Repr(M)$ for the vector space $W^+$ of
representative of ${\rm coker}(\Phi -1,M)$. We note that $Repr(M)$ does not
depend on the choice of the module $N$ over $K_n$ with $M=Res(N)$. 

\medskip

$M\mapsto Repr(M)$ has some functorial properties, namely:

\smallskip
\noindent
{\it A morphism $f:M_1\rightarrow M_2$ between pure modules of the same
negative slope maps $Repr(M_1)$ to $Repr(M_2)$}.\\
This follows easily from the second description of $W^+$.

\smallskip
\noindent
{\it Let $M_1,M_2$ denote pure modules of slopes $\lambda _1,\lambda _2<0$.
Then\\ $Repr(M_1)\otimes _{\mathbb C}Repr(M_2)$ is mapped to
$Repr(M_1\otimes M_2)\subset M_1\otimes M_2$}. 

\smallskip
\noindent 
The proof goes as follows. The pure module $M_3:=M_1\otimes _KM_2$ has slope 
$\lambda _3=\lambda _1+\lambda _2$. Choose an integer $n\geq 1$ such that 
$n\lambda _1,n\lambda _2\in {\mathbb Z}$. For each $i=1,2,3$, the vector space
$Repr(M_i)$ is equal to 
$\oplus _{|q|^{-n\lambda _i}<|c|\leq 1}Eigen((z^{-\lambda _i}\Phi)^n,c)$.
The statement that we want to prove follows from
$(z^{\lambda _i}\Phi )^n=q^{-\lambda  _in(n-1)/2}z^{-n\lambda _i}\Phi ^n$ and 
the observation that 
$q^{-\lambda _3n(n-1)/2}z^{-n\lambda _3}\Phi ^n(w_1\otimes w_2)$ equals
\[q^{-\lambda _1n(n-1)/2}z^{-n\lambda  _1}\Phi ^n(w_1)\otimes 
q^{-\lambda  _2n(n-1)/2}z^{-n\lambda _2}\Phi ^n(w_2)\]
for $w_1\in Repr(M_1),\ w_2\in Repr(M_2).$

\bigskip

There is a second method, based on subsection 1.4.1, to define a subspace
of representatives for ${\rm coker}(\Phi -1,N)$ with $N$ pure of slope
$\lambda <0$. Let $data(N)=(\lambda ,V,A,\{D(\sigma ) \})$ and 
$K_\infty \otimes N=K_\infty \otimes _{\mathbb C}V$. Then 
\[V^+:=\sum _{s\in {\mathbb Q},\ 0\leq s<-\lambda }z^s\otimes V  \]
can be seen to be a vector space of representative for 
${\rm coker}(\Phi -1,K_\infty \otimes V)$. This set of representatives
is invariant under the action of $Gal$ and therefore 
$Repr^*(N):=(V^+)^{Gal}$ (depending on the choice of $L\subset {\mathbb C}$) 
is a $\mathbb C$-subspace of $N_{global}$ representing 
${\rm coker}(\Phi -1,N)$. Again $N\mapsto
Repr^*(N)$ has the same functorial properties as $M\mapsto Repr(M)$.
However the behaviour of $Repr^*(N)$ with respect to tensor products is
more transparent. }\end{leeg}

\section{Moduli spaces for $q$-difference equations}
For a difference module $M$ over $K$ we consider again the slope filtration
$0=M_0\subset M_1\subset \cdots \subset M_r$. One associates to $M$ 
the graded module $gr(M):=\oplus _{i=1}^r M_i/M_{i-1}$. The aim is to
produce a moduli space for the collection of all $q$-difference modules
over $K$ with a fixed $gr(M)$. This problem is the theme of
[R-S-Z].  As Ch.~Zhang has remarked, the problem is already present in a paper 
of G.~Birkhoff and P.E.~Guenther, see [Bir]. 
Here we treat the general case (i.e., arbitrary, rational slopes). 

\smallskip

Fix a split difference module $S=P_1\oplus P_2\oplus \cdots \oplus P_r$
(with $r\geq 2$) over $K$ such that each $P_i$ is pure of slope
$\lambda _i$ and $\lambda _1<\lambda _2<\cdots <\lambda _r$. One 
considers the pairs $(M,f)$ consisting of a difference module over
$K$ and an isomorphism $f:gr(M)\rightarrow S$. Two pairs 
$(M(i),f(i)),\ i=1,2$ are called equivalent if there exists an isomorphism
$f:M(1)\rightarrow M(2)$ such that 
$gr(M(1))\stackrel{gr(f)}{\rightarrow} gr(M(2))
\stackrel{f(2)}{\rightarrow}S$ coincides with $f(1)$. One wants to
give the collection of equivalence classes a structure of algebraic
variety over $\mathbb C$. A naive, but useful way is to produce 
in every equivalence class a unique representative $(M,f)$. Our aim is to
define a moduli functor ${\mathcal F}$ and to find a fine moduli space for
${\mathcal F}$.

\smallskip
It is useful to give another formulation for the pairs $(M,f)$.
Let $U$ be the group of the $K$-linear maps $u:S\rightarrow S$ such that
$u$ respects the filtration $P_1\subset P_1\oplus P_2\subset
\cdots \subset S$ of $S$ and moreover $gr(u)=id$. 
Let $\Phi =\Phi _S$ denote the given action on $S$.
For any $u\in U$ one defines the $q$-difference module $M$ by
$M:=S$ provided with a new $\Phi$-action, namely $u\Phi$. Let 
${\mathcal F}^+({\mathbb C})$ denote the set of all actions $\{u\Phi\}$. Two 
actions $u_1\Phi$ and $u_2\Phi$ are equivalent if there exists an $a\in U$
such that $u_1\Phi a=au_2\Phi$. Let ${\mathcal F}({\mathbb C})$ denote the set
of equivalence classes. We note that this coincides with the set of
equivalence classes for the pairs $\{(M,f)\}$.

For any ${\mathbb C}$-algebra $R$ (i.e., commutative and with a unit element)
one defines $U(R)$ as the set of the $R\otimes _{\mathbb C}K$-linear maps $u$
from $R\otimes _{\mathbb C}S$ to itself such $u$
respects the filtration $R\otimes _{\mathbb C} P_1\subset R\otimes _{\mathbb
  C} (P_1\oplus P_2)\cdots
\subset R\otimes _{\mathbb C} S$ and moreover $gr(u)=id$. Further
${\mathcal F}^+(R)$ denotes the set of all actions $\{u\Phi \}$ with 
$u\in U(R)$ on 
$R\otimes _{\mathbb C}S$. As before, two actions $u_1\Phi$ and $u_2\Phi$ are
equivalent if there exists an $a\in U(R)$ such that 
$a^{-1}u_1\Phi a\Phi ^{-1}=u_2$.
The set of equivalence classes is denoted by ${\mathcal F}(R)$. Thus we obtain
a covariant functor ${\mathcal F}$ defined on the category of the
$\mathbb C$-algebras. One can view ${\mathcal F}$ as a contravariant functor
on affine schemes over $\mathbb C$ and extend ${\mathcal F}$ to a contravariant
functor from $\mathbb C$-schemes to the category of sets. Our aim is to show
that $\mathcal F$ is representable by a certain $\mathbb C$-algebra, or
in other terms by an affine scheme over $\mathbb C$. It will turn out that
the scheme representing $\mathcal F$ is in fact ${\mathbb A}_{\mathbb C}^N$
for some $N$.

\begin{example} $S=P_1\oplus P_2$, with pure modules $P_1,P_2$ of slopes 
$\lambda _1<\lambda _2$.\\ {\rm 
Any element in $U$ has the form $u=1+u_{1,2}$ with
$u_{1,2}:P_2\rightarrow P_1$ a $K$-linear map. Further 
$v:=a^{-1}u\Phi a\Phi ^{-1}$ satisfies $v_{1,2}=u_{1,2}-a_{1,2}+\Phi
(a_{1,2})$. Here $\Phi (a_{1,2})$ denotes the action of $\Phi$ on the element
$a_{1,2}$ of the pure module $B:={\rm Hom}_K(P_2,P_1)=P_2^*\otimes P_1$. 
The map $B\rightarrow {\mathcal F}({\mathbb C})$ yields an isomorphism 
${\rm coker}(\Phi -1,B)\rightarrow {\mathcal F}({\mathbb C})$.  

For any $\mathbb C$-algebra $R$, one obtains in the same way one isomorphism
${\rm coker}(\Phi -1,R\otimes _{\mathbb C}B)\rightarrow {\mathcal F}(R)$. 
In 2.6, one has defined the ${\mathbb C}$-vector space 
$Repr(B)\subset B_{global}\subset B$ of representatives for 
${\rm coker}(\Phi -1,B)$. One introduces the functor 
${\mathcal F}^o$ by ${\mathcal F}^o(R)= R\otimes _{\mathbb C}Repr(B)$.
One obtains an isomorphism of functors 
${\mathcal F}^o(R)\rightarrow {\mathcal F}(R)$. 

The finite dimensional $\mathbb C$-vector space $Repr(B)$ is seen as  
complex algebraic variety. Its algebra of regular functions $O(Repr(B))$ 
is the symmetric algebra of $Repr(B)^*$. Then
\begin{small} 
\[{\rm Hom}_{{\mathbb C}-{\rm algebra}}(O(Repr(B)),R)=
{\rm Hom}_{{\mathbb C}-{\rm vector space}}(Repr(B)^*,R)=
Repr(B)\otimes R \ .\]
\end{small}
Thus ${\mathcal F}^o$ and ${\mathcal F}$ are represented by $O(Repr(B))$. In 
terms of schemes, $\mathcal F$ is represented by $Repr(B)$, seen as affine
space over $\mathbb C$. We note that $Repr(B)$ is in fact the $\mathbb
C$-vector space $Ext^1(P_2,P_1)$, where $P_1,P_2$ are seen as left 
$K[\Phi ,\Phi ^{-1}]$-modules. The universal family is made explicit in the
following examples.

\medskip 

\noindent (1) $P_i=Ke_i$ for $i=1,2$ and $\Phi e_1=e_1,\ \Phi e_2=z^te_2$.
Then the universal family is $K[x_0,\dots ,x_{t-1}]e_1+
K[x_0,\dots ,x_{t-1}]e_2$ with $\Phi$ given by $\Phi e_1=e_1$ and
$\Phi e_2=z^te_2+(x_0+x_1z+\cdots +x_{t-1}z^{t-1})e_1$.

\smallskip

\noindent (2) $P_1=Res(K_ne_1)$ with $\Phi e_1=z^{-t/n}e_1,\ t,n\geq 1,\
g.c.d.(t,n)=1$ and $P_2=Ke_2$ with $\Phi e_2=e_2$. The universal family is:\\
$K_n[x_0,\dots ,x_{t-1}]e_1+K[x_0,\dots ,x_{t-1}]e_2$ with $\Phi$ given by
$\Phi e_1=z^{-t/n}e_1$ and $\Phi e_2=e_2+
(x_0+x_1z^{1/n}+\cdots +x_{t-1}z^{(t-1)/n})e_1$. \hfill $\square$
}\end{example}

\begin{theorem} Let $S=P_1\oplus \cdots \oplus P_r$ be a direct sum of
pure modules with slopes $\lambda _1<\cdots <\lambda _r$. The functor 
$\mathcal F$ associated to $S$ is represented by a free polynomial ring
over $\mathbb C$ in
$N:=\sum _{i<j}(\lambda _j-\lambda _i)\cdot \dim _KP_i\cdot \dim _K P_j$
variables. Equivalently, $\mathcal F$, seen as a contravariant functor
from $\mathbb C$-schemes to sets, is represented by the affine space
${\mathbb A}^N_{\mathbb C}$. 
\end{theorem}
\begin{proof} We will use induction with respect to $r$. 
The case $r=2$ is dealt with in Example 3.1. Take $r=3$. Write 
${\mathcal F}_{1,2,3}$ for
the functor associated to $P_1\oplus P_2\oplus P_3$. Further 
${\mathcal F}_{1,2},\ {\mathcal F}_{2,3}, {\mathcal F}_{1,3}$ are the functors
associated to $P_1\oplus P_2,\ P_2\oplus P_3,\ P_1\oplus P_3$. 
There is a morphism of functors $T:
{\mathcal F}_{1,2,3}\rightarrow {\mathcal F}_{1,2}\times {\mathcal F}_{2,3}$.
An element of ${\mathcal F}_{1,2,3}(R)$, represented by a filtration
$M_1\subset M_2\subset M_3$, is mapped to the pair of equivalence classes
$([M_1\subset M_2],[M_2/M_1\subset
M_3/M_1])$ in ${\mathcal F}_{1,2}(R)\times {\mathcal F}_{2,3}(R)$. We 
{\it claim} that ${\mathcal F}_{1,2,3}$ is a {\it trivial torsor} over 
${\mathcal G}:={\mathcal F}_{1,2}\times {\mathcal F}_{2,3}$ for the algebraic
group ${\mathcal F}_{1,3}$. In other words, we will produce an
action $m$, functorial in $R$, 
\[m:{\mathcal F}_{1,3}(R)\times {\mathcal F}_{1,2,3}(R)\rightarrow
{\mathcal F}_{1,2,3}(R)\ ,\]
such that the map $(g,h)\mapsto (m(g,h),h)$ from 
${\mathcal F}_{1,3}(R)\times {\mathcal F}_{1,2,3}(R)$ to the fibre product
${\mathcal F}_{1,2,3}(R)\times _{{\mathcal G}(R)}{\mathcal F}_{1,2,3}(R)$
is a bijection. We note that this definition becomes the usual one, after
introducing the algebraic group $({\mathcal F}_{1,3})_{\mathcal G}$ over 
${\mathcal G}$, by the formula $({\mathcal F}_{1,3})_{\mathcal G}(R)=
{\mathcal F}_{1,3}(R)\times {\mathcal G}(R)\; $.

The triviality of the torsor means that there is an isomorphism of functors
 ${\mathcal F}_{1,3}\times ({\mathcal F}_{1,2}\times {\mathcal F}_{2,3})
\rightarrow {\mathcal F}_{1,2,3}$ compatible with $T$. The last statement
and the case $r=2$ imply the theorem for the case $r=3$.

\medskip

We represent, as before, an element of ${\mathcal F}_{1,2,3}(R)$ as an
equivalence class of actions $u\Phi$ on $S$ with $u\in U(R)$. Let $\sim$
denote the equivalence relation on $U(R)$, given by $u_1\sim u_2$ if
there exists an $a\in U(R)$ with $a^{-1}u_1\Phi a \Phi ^{-1}=u_2$.
Then $U(R)/\sim$ identifies with ${\mathcal F}_{1,2,3}(R)$. 
Let $U'(R)$ denote the subgroup of $U(R)$ consisting of the elements
$u$ such that $u-1$ is $0$ on $R\otimes _{\mathbb C}(P_1\oplus P_2)$ and
$R\otimes _{\mathbb C}P_3$ is mapped to $R\otimes _{\mathbb C}P_1$. On 
$U'(R)$ we
introduce the equivalence relation $\sim$ by $u_1\sim u_2$ is there exists an 
$a\in U'(R)$ with $a^{-1}u_1\Phi a\Phi ^{-1}=u_2$. 
Now $U'(R)/\sim$ identifies in an obvious way with ${\mathcal F}_{1,3}(R)$.

\smallskip
A small calculation, based on the observation that $U'(R)$ 
lies in the center of $U(R)$, shows that the map 
$U'(R)\times U(R)\rightarrow U(R)$, given by 
$(u',u)\mapsto u'u$, has the property that $u_1'\sim u_2'$ and $u_1\sim u_2$
imply $u_1'u_1\sim u_2'u_2$. Thus there is an induced action
$m:(U'(R)/\sim )\times (U(R)/\sim )\rightarrow (U(R)/\sim )$. One can verify 
that $m$ defines a group action, that $m$ has the torsor property over
${\mathcal G}(R)$, and that the construction is
functorially in $R$. Finally, we want to show that the torsor is trivial.
Thus we need to define a functorial isomorphism
${\mathcal F}_{1,3}(R)\times {\mathcal G}(R)\rightarrow {\mathcal
  F}_{1,2,3}(R)$. 
This is done as follows. One considers the map
\[{\mathcal F}_{1,3}^o(R)\times {\mathcal F}^o_{1,2}(R)\times {\mathcal
 F}^o_{2,3}(R) \rightarrow {\mathcal F}^+_{1,2,3}(R)\mbox{ (see Example 3.1) },
\]  given by $(u_{1,3},u_{1,2},u_{2,3})\mapsto 
\left(\begin{array}{ccc}1&u_{1,2}&u_{1,3}\\ 0&1& u_{2,3}\\ 0 &
    0&1\end{array}\right) \cdot \Phi$.
One easly verifies that the resulting map
\[{\mathcal F}_{1,3}^o(R)\times {\mathcal F}^o_{1,2}(R)\times {\mathcal
 F}^o_{2,3}(R) \rightarrow {\mathcal F}_{1,2,3}(R),\] 
is bijective and depends functorially on $R$. 
Since the maps ${\mathcal F}^o_{i,j}(R)\rightarrow {\mathcal F}_{i,j}(R)$
are isomorphisms functorial in $R$, we have found a trivialization of the
torsor. This ends the proof for the case $r=3$. 

\medskip
 
For $r=4$ and $S=P_1\oplus \cdots \oplus P_4$, one defines in a 
similar way functors ${\mathcal F}_{1,2,3,4},\ {\mathcal F}_{2,3,4}$ etc.
With the same methods one proves that the morphism of functors
${\mathcal F}_{1,2,3,4}\rightarrow {\mathcal G}$, where 
${\mathcal G}:={\mathcal F}_{1,2,3}\times _{{\mathcal F}_{2,3}}
{\mathcal F}_{2,3,4}$, is a trivial torsor for the algebraic group
${\mathcal F}_{1,4}$. It is clear how to extend this to any $r>4$.
\end{proof}

\begin{remarks} {\rm
(1) Let $S=P_1\oplus \cdots \oplus P_r$ 
with corresponding functor $\mathcal F$.
The method of the proof of Theorem 3.2 yields maps 
$\prod _{i<j}{\mathcal F}_{i,j}^o(R)\rightarrow {\mathcal F}^+(R)$, functorial
in $R$. Let ${\mathcal F}^o(R)$ denote the image of the this map.
Then ${\mathcal F}^o$ is a functor and the obvious morphism
${\mathcal F}^o\rightarrow {\mathcal F}$ is an isomorphism. In particular one
finds an isomorphism of functors 
$\prod _{i<j}{\mathcal F}_{i,j}\rightarrow {\mathcal F}$ which is obtained
by trivializing the sequence of torsors involved in $\mathcal F$. \\
(2) The case of Theorem 3.2 where the slopes are integers is present
in a paper of G.~Birkhoff and P.E.~Guenther, see [Bir], where they present
normal forms for these equations. This case is also treated in 
[R-S-Z, Sau2, Sau3].\\
(3) Each ${\mathcal F}_{i,j}$ has the natural structure of algebraic group
over $\mathbb C$ and is in fact equal to $Ext ^1(P_i,P_j)$. However,
$\mathcal F$ has no evident structure of (unipotent) algebraic group. In 
contrast with this, the obvious injective map 
$\prod _{i<j}{\mathcal F}^o_{i,j}(R)\rightarrow U(R)$ has as image a
{\it subgroup} of $U(R)$. 

In proving this, one has to consider indices $a<b<c$ and one has to show that
the the obvious map ${\mathcal F}^o_{a,b}(R)\times {\mathcal F}^o_{b,c}(R)$
to $R\otimes {\rm Hom}(P_c,P_a)$ has image in 
${\mathcal F}^o_{a,c}(R)$. This follows from the statement:
$Repr(M_1)\otimes Repr(M_2)$ is mapped to 
$Repr(M_1\otimes M_2)\subset M_1\otimes M_2$, proved in 2.6.

The functorial isomorphism $\prod _{i<j}{\mathcal F}^o_{i,j}\rightarrow
{\mathcal F}$ provides the latter with a structure of unipotent linear
algebraic group.\hfill $\square$ }\end{remarks}

\section{Global difference modules}

The skew ring ${\mathbb D}:={\mathbb C}[z,z^{-1}][\Phi ,\Phi
 ^{-1}]$ is defined by the relation $\Phi z= qz\Phi$. We recall that a {\it
global difference module} is a left ${\mathbb D}$-module $N$, which is as 
${\mathbb C}[z,z^{-1}]$-module, finitely generated and free.

It suffices in fact to assume that $N$ is finitely generated as
${\mathbb C}[z,z^{-1}]$-module. Indeed, let $N_0$ denote the torsion
submodule of $N$ and let $I\subset {\mathbb C}[z,z^{-1}]$ denote the 
anihilator ideal of $N_0$. Then $N_0$ and $I$ are invariant under the action
of $\Phi$ and $\phi$. This implies $I=(1)$, $N_0=\{0\}$ and $N$ is free.

\begin{lemma} Any global difference module has a cyclic vector.
\end{lemma}
\begin{proof} We imitate the proof of [vdP-S2], Proposition 2.9. Let the 
global difference module
$N$ have rank $n$ over ${\mathbb C}[z,z^{-1}]$.  It suffices to find an
element $m\in N$ such that the ${\mathbb C}[z,z^{-1}]$-module $N'$ generated by
$\{\Phi ^sm|\ s\geq 0\}$ has rank $n$. Indeed, this implies that the
global difference module $N/N'$ has rank 0 and is therefore equal to 0.
This condition on $m$ translates into $m\wedge \Phi m\wedge \Phi ^2m\wedge
\cdots \wedge \Phi ^{n-1}m\neq 0$ as element of the global difference
module $\Lambda ^nN$.

We recall that for any integer $d\geq 1$, the global difference 
module $\Lambda ^dN$ is defined as the $d$th exterior product of
$N$ as ${\mathbb C}[z,z^{-1}]$-module with $\Phi$-action given by
$\Phi (m_1\wedge \cdots \wedge m_d)=(\Phi m_1)\wedge \cdots \wedge (\Phi
m_d)$. 

Suppose that we have found an element $e\in N$ such that
$\{\Phi ^se|\ s\geq 0\}$ generates a ${\mathbb C}[z,z^{-1}]$-module $N'$
of rank $m<n$. Then $e\wedge \cdots \wedge \Phi ^{m-1}e\neq 0$ and
$e\wedge \cdots \wedge \Phi ^me=0$. Take an element $f\in N\setminus N'$ and
consider $\tilde{e}=e+\lambda z^sf$ with $\lambda \in {\mathbb Q}$ and
$s\in {\mathbb Z}$. We claim that for suitable $\lambda$ and $s$ one has
that $E(\lambda ,s):=\tilde{e}\wedge \cdots \wedge \Phi ^m\tilde{e}\neq
0$. From this claim the lemma follows. One can write $E(\lambda ,s)$ as a sum
of terms, which are wedge products of some $\Phi ^ue$ and some $\Phi ^vf$
(like $e\wedge \Phi e\wedge \Phi ^2f\wedge \cdots \wedge \Phi ^mf$)
with coefficients $(\lambda z^s)^aq^b$ for $a=0,\dots ,m+1$ and suitable
$b$. If $E(\lambda ,s)=0$ for all $\lambda$ and $s$, then each of these
wedge products is 0. This contradicts 
$e\wedge \Phi e\wedge \cdots \wedge \Phi ^{m-1}e\wedge  \Phi ^mf\neq 0$.  
\end{proof}

A global module $N$ will be called {\it pure} of slope $\lambda$ 
if ${\mathbb C}(\{z\})\otimes _{{\mathbb C}[z,z^{-1}]}N$ is pure with slope
$\lambda$ and 
${\mathbb C}(\{z^{-1}\})\otimes _{{\mathbb C}[z,z^{-1}]}N$ is pure with
slope $-\lambda$. An example of a pure global module with slope $t/n$ (and 
$n\geq 1,\ g.c.d.(t,n)=1$) and rank $n$ is ${\mathbb C}[z^{1/n},z^{-1/n}]e$ 
with $\Phi e=cz^{t/n}e$ and $c\in {\mathbb C}^*$.  

\begin{lemma} {\rm Pure global modules.}\\
 {\rm (1)} A global module $N$ is pure with slope $t/n$ 
{\rm (}with $n\geq 1$ and $g.c.d.(t,n)=1${\rm )} if and only if 
$N\cong {\mathbb D}/{\mathbb D}L(z^{-t}\Phi ^n)$ for a monic
$L\in {\mathbb C}[T]$ with $L(0)\neq 0$.\\
{\rm (2)} Let $L\in {\mathbb C}[T]$ have the form 
$\prod _{j=1}^s(T-c_j)^{m_j}$ with distinct $c_1,\dots c_s\in {\mathbb C}^*$.
Then ${\mathbb D}/{\mathbb D}L(z^{-t}\Phi ^n)$ is the direct sum of
the indecomposable global modules 
${\mathbb D}/{\mathbb D}(z^{-t}\Phi ^n-c_j)^{m_j}$. \\
{\rm (3)}
Let $N$ be a pure global module of slope $t/n$ and $m\in N$. The 
operator $S$ in ${\mathbb C}[z,z^{-1}][\Phi ]$ of minimal degree in
$\Phi$, satisfying $S(m)=0$, has the form $S=P(z^{-t}\Phi ^n)$ with
monic $P\in {\mathbb C}[T ]$ and $P(0)\neq 0$.
\end{lemma}
\begin{proof} (1) Suppose that $N$ is pure with slope $t/n$ and has rank 
$m\cdot n$. Consider a cyclic vector $e$ of $N$ and let $L\in {\mathbb
  C}[z,z^{-1}][\Phi ]$ denote a non zero element of minimal degree (namely
$m\cdot n$) in $\Phi$  such that $Le=0$. Clearly the constant term $L(0)\in
{\mathbb C}[z,z^{-1}]$ is different from 0. After multiplication by an
invertible element of ${\mathbb C}[z,z^{-1}]$, one may suppose that 
$L(0)\in {\mathbb C}[z]$ and that the constant term of $L(0)$ is 1. 
Now $\Phi ^{-1}$ acts bijectively on $N\cong {\mathbb D}/{\mathbb D}L$.
This implies that $L(0)=1$. The Newton polygon of $L$, considered over
${\mathbb C}(\{z\})$, contains only terms $z^a\Phi ^b$ with 
$a+b\frac{t}{n}\geq 0$. The Newton polygon of $L$, considered over the field
${\mathbb C}(\{z^{-1}\})$ contains only terms $z^a\Phi ^b$ with
$a+b\frac{t}{n}\leq 0$. Thus $L\in {\mathbb C}[z^{-t}\Phi ^n ]$ and $L(0)=1$. 
After multiplication by a constant, we may suppose that $L$ is monic.

On the other hand, a global module, given as ${\mathbb D}/{\mathbb D}L$
with $L\in {\mathbb C}[\Psi ]$ and $L(0)\neq 0$, is clearly a
pure global module with slope $t/n$. 

\medskip
 The proof of (2) is straightforward and (3) follows from the observation that
 ${\mathbb D}m$ is again a pure global module of slope $t/n$. 
\end{proof}

A global difference module has, in general, no slope filtration! We consider 
the category of the global modules $N$ that have {\it an ascending slope 
filtration}, i.e., a sequence of submodules 
$0=N_0\subset N_1\subset \cdots \subset N_r=N$, such
that each quotient $N_i/N_{i-1}$ is pure with slope $\lambda _i$ and 
$\lambda _1<\cdots <\lambda _r$.

\begin{theorem} The functor $T$ from the category of the global difference 
modules with ascending slope filtration to the category of the difference
modules over $K$, given by $N\mapsto K\otimes _{{\mathbb C}[z,z^{-1}]}N$,
is an equivalence of Tannakian categories. \end{theorem}

\begin{remarks}{\rm
For a given module $M$ over $K$, the global module with ascending slope
filtration $N$ with $T(N)=K\otimes N\cong M$ has the property that
${\mathbb C}(\{z^{-1}\})\otimes N$ is a difference module over the field
${\mathbb C}(\{z^{-1}\})$ which is a direct sum of pure modules. Indeed,
${\mathbb C}(\{z^{-1}\})\otimes N$ has both an ascending and a descending slope
filtration. In other words, $N$ is an algebraic vector bundle above 
${\mathbb P}^1-\{0,\infty \}$ with a $\Phi$-action and a filtration according
to the slopes at $z=0$.  

One can see Theorem 4.3 as an analogue of a theorem of G.~Birkhoff which
states that every differential module over $K$ comes from a connection
on ${\mathbb P}^1$ which has only singularities at $0$ and $\infty$.
Moreover the singularity at $\infty$ is regular singular. 

{\it The proof of Theorem 4.3 is given in the following series of observations 
}}\end{remarks}

\begin{observations}$\;$ \\ {\rm 
(1) Let $N$ be a global difference module with ascending slope filtration
$0=N_0\subset N_1\subset \cdots \subset N_r=N$ such that $N_i/N_{i-1}$ is
pure of slope $\lambda _i=t_i/n_i$ with $n_i\geq 1,\ g.c.d. (t_i,n_i)=1$
and $\lambda _1<\cdots <\lambda _r$.

\smallskip
\noindent
{\it Let $m\in N$. The difference operator 
$L\in {\mathbb C}[z,z^{-1}][\Phi ]\subset {\mathbb D}$, monic and of minimal 
degree in $\Phi$, satisfying $L(m)=0$, 
has the form 
\[L=q^{a}z^{b}L_1(z^{-t_1}\Phi ^{n_1})\cdots L_r(z^{-t_r}\Phi ^{n_r})\] 
with $L_1,\dots ,L_r$ monic elements of ${\mathbb C}[T]$ and suitable
$a,b\in {\mathbb Z}$.}  
\begin{proof} $L_r$ is the polynomial of Lemma 4.2, applied to the image of
$m$ in $N_r/N_{r-1}$. If $m$ happens to be an element of $N_{r-1}$, then
$L_r=1$. Lemma 4.2 applied to the image of $L_r(z^{-t_r}\Phi ^{n_r})(m)$
in $N_{r-1}/N_{r-2}$ produces $L_{r-1}$. Induction finishes the proof. One
multiplies the above operator with a suitable
term $q^az^b$ to obtain a monic operator.\end{proof}

\noindent (2) Using (1) one deduces that any morphism between global
difference modules with ascending slope filtrations respects the filtrations.  
Further the full subcategory of the category of all global modules, whose
objects are the global modules with ascending filtration is closed under
direct sums, tensor products, duals, submodules and quotients. In particular,
this subcategory is a Tannakian category.

\medskip

\noindent (3) The indecomposable difference module 
$M:=Res(K_ne)\otimes _KU_m$ with $\Phi e=cz^{t/n}$ has a cyclic vector $f$ such
that $L=(z^{-t}\Phi ^n-c^nq^{t(n-1)/2})^m$ is a polynomial of minimal degree in
$\Phi$ with $L(f)=0$. Clearly, $M=T(N)$ with
$N={\mathbb D}/{\mathbb D}(z^{-t}\Phi ^n-c^nq^{t(n-1)/2})^m$ and 
$M_{global}$ coincides, according to Definition 1.8, with $N$. 

One concludes that for a split difference module $M$ over $K$ there exists
a unique (up to isomorphism) global module $N$, direct sum of pure global
ones, such that $T(N)\cong M$.

\medskip

\noindent (4) Let $N$ be as in (1). The associated graded global module 
$gr(N)$ is defined as $\sum _{j=1}^rN_j/N_{j-1}$.

Let a direct sum of pure global modules $U=R_1\oplus \cdots \oplus R_r$ with
slopes $\lambda _1<\cdots <\lambda _r$ be given. As in section 3, one
considers the set of equivalence classes of the pairs $(N,f)$ consisting
of a global module with an ascending slope filtration and an isomorphism
$f:gr(N)\rightarrow U$. Two pairs $(N_i,f_i)$ are equivalent if there exists
an isomorphism $\alpha :N_1\rightarrow N_2$ such that 
$gr(N_1)\stackrel{gr(\alpha )}{\rightarrow}gr(N_2)\stackrel{f_2}{\rightarrow}
S$ coincides with $f_1$. As in section 3 one proves that the set of
equivalence classes is in a natural way isomorphic to
$\prod _{i<j}{\rm coker}(\Phi -1,{\rm Hom}(R_j,R_i))$.

\medskip

\noindent (5) We use the notation of (4) and put $S=T(U)$ and $P_j=T(Q_j)$.
The functor $T$ maps the set of equivalence classes for
$U=R_1\oplus \cdots \oplus R_r$ to the set of equivalence classes
for the graded module $S=P_1\oplus \cdots \oplus P_r$ over $K$, studied in
section 3. From Theorem 2.5 one easily concludes that $T$ defines an
isomorphism between the two classes of objects.

This implies that the functor $T$ induces a bijection between the isomorphy
classes of global difference modules with an ascending filtration and the
isomorphy classes of difference modules over $K$.

\medskip
\noindent (6) The final part of the proof of Theorem 4.3 consists of verifying
that the map ${\rm Hom}(N_1,N_2)\rightarrow {\rm Hom}(T(N_1),T(N_2))$ is
bijective. Since $T$ clearly commutes with tensor products we may suppose
that $N_1={\mathbb C}[z,z^{-1}]e$ with $\Phi e=e$. Thus we have to show, for
any global module $N$ with ascending filtration, that
${\rm ker}(\Phi -1,N)\rightarrow {\rm ker}(\Phi -1,T(N))$ is bijective.
The injectivity is obvious.

Let $0=N_0\subset N_1\subset \cdots \subset N_r=N$ be the slope filtration.
Then $0=M_0\subset M_1\subset \cdots \subset M_r=M$ with $M_j=T(N_j)$ is the
slope filtration of $M=T(N)$.

Let $N'$ be the ${\mathbb C}[z,z^{-1}]$-submodule generated by
$\{m\in N|\ \Phi (m)=m\}$. Then $N'$ is a pure global submodule of $N$
with slope 0. Define $s\leq r$ to be the smallest integer such that
$N'\subset N_s$. The map $\alpha :N_s\rightarrow N_s/N_{s-1}$ is not zero on
$N'$. Hence $\lambda _s=0$. The restriction of $\alpha$ to $N'$ is injective.
Indeed, according to (1), a non zero element of $N_{s-1}$ cannot have minimal
polynomial $\Phi -1$. 

Let $M'$ denote the $K$-subspace of $M$ generated by $\{m\in M|\Phi (m)=m\}$.
Then $M'$ is a pure submodule of $M$ with slope 0. Thus $M'\subset M_s$ and
the restriction of $T(\alpha ): M_s\rightarrow M_s/M_{s-1}$ to $M'$ 
is injective. The image $T(\alpha )(M')$ contains $T(\alpha )(N')$. 
Let $M^-\subset M_s$ denote the preimage $T(\alpha )^{-1}T(\alpha )(M')$.
Then $M'\rightarrow M^-/M_{s-1}$ is a bijection and thus 
$M^-$ is a direct sum $M_{s-1}\oplus M'$. The global module $N^-$ with 
ascending filtration and $T$-image $M_{s-1}\oplus M'$ is clearly 
$N_{s-1}\oplus {\mathbb C}[z,z^{-1}]\otimes _{\mathbb C}V$ where $V$ is a
$\mathbb C$-vector space where $\Phi$ acts as the identity. Then 
${\rm ker}(\Phi -1,N)={\rm ker}(\Phi -1,N^-)=V={\rm ker}(\Phi -1,M^-)=
{\rm ker}(\Phi -1,M)$.\hfill $\square$
}\end{observations}

\begin{remarks}$\;$\\ {\rm
(1) Let $N$ be as in part (1) of Observation 4.5 and put 
$M=T(N)=K\otimes _{{\mathbb C}[z,z^{-1}]}N$. The elements $m$ of $N=M_{global}$
are characterized by the 
condition that the monic polynomial $L$ of minimal degree
in $\Phi$, satisfying $L(m)=0$, has the form
\[L=q^{a}z^{b}L_1(z^{-t_1}\Phi ^{n_1})\cdots L_r(z^{-t_r}\Phi ^{n_r})\] 
with $L_1,\dots ,L_r$ monic elements of ${\mathbb C}[T]$ and suitable
$a,b\in {\mathbb Z}$.

\medskip
\noindent
(2) Let $N$ be a pure global module of slope $\lambda >0$, then
${\rm coker}(\Phi -1,K\otimes N)$ is zero. Moreover, ${\rm coker}(\Phi -1,N)$
is a $\mathbb C$-vector space of dimension $\lambda \cdot {\rm rank} N$ and  
the natural  map ${\rm coker}(\Phi -1,N)\rightarrow
{\rm coker}(\Phi -1,{\mathbb C}(\{z^{-1}\})\otimes N)$ is a bijection.
}\end{remarks}

\section{The Galois group of a $q$-difference module}

\subsection{Universal Picard-Vessiot extensions}

We recall that a difference module $M$ over $K$ is {\it split} if $M$ is a 
direct sum of pure modules. A split difference module $M$ over $K$ has
a Picard-Vessiot extension $PV(M)$ and the (difference) Galois group
is the group of the automorphisms of $PV(M)/K$ which commute with the
action of $\phi$ on $PV(M)$. The {\it universal Picard-Vessiot extension
$Univ$} for the family of all split modules over $K$ is the direct limit of
all $PV(M)$. The universal difference Galois group ${\mathbb G}_{univ}$
is the group of all $K$-automorphisms of $Univ$ which commute with the action
of $\phi$ on $Univ$. The restriction of ${\mathbb G}_{univ}$ to the subring
$PV(M)\subset Univ$ is the Galois group of $M$.

 We note that the Picard-Vessiot ring of the difference module 
$\widehat{K}\otimes M$ over $\widehat{K}$ is simply 
$\widehat{K}\otimes _KPV(M)$. In particular, $M$ and $\widehat{K}\otimes _KM$ 
have the same Galois group. Moreover, $\widehat{K}\otimes _KUniv$ is the
universal Picard-Vessiot extension for the collection of {\it all} difference
modules over $\widehat{K}$ since every difference module over $\widehat{K}$
is split.

We give the explicit description of $Univ$, given in [vdP-S1], p.150, in a
slightly changed form. $Univ:=K[\{e(cz^\lambda )\},\ell ]$, with $c\in
{\mathbb C}^*$ and $\lambda \in {\mathbb Q}$. The only relations are:
\[e(c_1z^{\lambda _1})e(c_2z^{\lambda _2})=e(c_1c_2z^{\lambda _1+\lambda _2}),
\ \ e(1)=1,\ \ e(q)=z^{-1}.\]
The algebraic closure of $K$ embeds in $Univ$, by identifying $z^\lambda$ with
$e(e^{-2\pi i \tau \lambda})$ for all $\lambda \in {\mathbb Q}$. 
The $\phi$-action on $Univ$ is given by
\[\phi (e(cz^\lambda ))=(cz^\lambda )^{-1}\cdot e(cz^\lambda ),\  \ \phi (\ell
)=\ell +1.\] 
The group ${\mathbb G}_{univ}$ consists of elements $\sigma =(h,s,a)$ with 
$h:{\mathbb C}^*\rightarrow {\mathbb C}^*$
a homomorphism with $h(q)=1$, $s:{\mathbb Q}\rightarrow {\mathbb C}^*$ a
homomorphism and $a\in {\mathbb C}$. The action of $\sigma$ is given by
\[\sigma (e(cz^\lambda ))=s(\lambda )\cdot h(c)\cdot 
e(h(e^{2\pi i \tau \lambda}))\cdot e(cz^\lambda),\ \  \sigma (\ell )=\ell +a.\]
We will produce topological generators (for the Zariski topology) for this 
rather complicated group. Define the elements 
$\Gamma ,\Delta \in {\mathbb G}_{univ}$ by 
\[\Gamma (e(cz^\lambda ))=
e^{2\pi i a_1} \cdot e(e^{2\pi i \lambda})\cdot e(cz^\lambda )
\mbox{ and } \Gamma (\ell )=
\ell +\frac{1}{\tau }\]  
\[\Delta (e(cz^\lambda ))=e^{-2\pi i a_0}\cdot e(cz^\lambda ) \mbox{ and }
\Delta (\ell )=\ell +1 ,\]
\[ \mbox{ where }c=e^{2\pi i(a_0+a_1\tau )}\mbox{ with } a_0,a_1\in {\mathbb
  R}.\] 
For the commutator $Com:=\Gamma \Delta \Gamma ^{-1}\Delta ^{-1}$ one
calculates the formulas
\[Com (e(cz^\lambda ))=e^{2\pi i\lambda}\cdot e(cz^\lambda )\mbox{ and }
Com(\ell )=\ell .\]
To any homomorphism $s:{\mathbb Q}\rightarrow {\mathbb C}^*$ one
associates the element $\tilde{s}\in {\mathbb G}_{univ}$ by the formulas
\[\tilde{s}(e(cz^\lambda ))=s(\lambda )\cdot e(cz^\lambda )\mbox{ and }
\tilde{s}(\ell )=\ell \ .\] 
We note that (any) $\tilde{s}$ commutes with $\Gamma$ and $\Delta$.
Write $D=\tilde{s}$ with $s$ given by $s(\lambda )=e^{2\pi i\tau \lambda}$.

\begin{proposition} $\Gamma ,\Delta ,D$ are topological generators
for ${\mathbb G}_{univ}$.
\end{proposition}
\begin{proof} Let $H$ be the group generated by $\Gamma ,\Delta , D$ and
let $Qt(Univ)$ denote the total quotient ring of $Univ$. We recall that 
$Univ$ is the direct limit of Picard-Vessiot rings $PV(M)$ of split modules 
$M$ over $K$. The group ${\mathbb G}_{univ}$ is the projective limit of the 
automorphism groups $G_M$ of these $PV(M)$. We have to show that the image
of $H$ in $G_M$ is Zariski dense. This is equivalent to the statement
that $K$ is the set of the $H$-invariant elements of the total quotient
ring of $PV(M)$. Thus the statement of the proposition is equivalent with
the set of the $H$-invariant elements of $Qt(Univ)$ is $K$. 

\smallskip
\noindent
(1) {\it If $\xi \in K[\{e(c)\},\ell ]$ is invariant under $\Gamma$ and
$\Delta$, then $\xi \in K$.}

\smallskip
\noindent {\it Proof of} (1).
Write $\xi =\sum _{i=0}^m \xi _i\ell ^i$ and all $\xi _i\in K[\{e(c)\}]$
and we may suppose that $\xi _m\neq 0$.
Now
$R:=\{e^{2\pi i (a_0+a_1\tau )}|\ a_0,a_1\in {\mathbb R},\ 0\leq a_0,a_1 <1\}$
is a set of representatives of ${\mathbb C}^*/q^{\mathbb Z}$.
Each $\xi _i$ has uniquely the form $\sum _{c\in R}a_{i,c}e(c)$. One sees that
$\xi _m$ is invariant under $\Gamma$ and $\Delta$. From 
$\Gamma e(c)=e^{2\pi ia_1}\cdot e(c)$ and $\Delta e(c)=e^{-2\pi i a_0}\cdot
e(c)$ it follows that $\xi _m\in K$. We may suppose that $\xi _m=1$ and we
have to prove that $m=0$. Suppose that $m>0$ and apply $\Gamma$ to $\xi$.
Comparing the coefficient of $\ell ^{m-1}$, one finds that 
$\Gamma (\xi _{m-1})+\frac{m}{\tau}=\xi _{m-1}$. This is clearly impossible
for $m>0$ and an element in $K[\{e(c)\}]$. 
We conclude that $m=0$ and $\xi \in K$.\hfill $\square$ 

\smallskip

\noindent (2) {\it Let $\xi \in Univ^H$, then $\xi \in K$.}

\smallskip
\noindent {\it Proof of }(2). One write 
$\xi =\sum _{\lambda \in {\mathbb Q}}a_\lambda  e(z^\lambda )$ with all
$a_\lambda \in K[\{e(c)\},\ell ]$. Then $D (\xi )=\sum a_\lambda e^{2\pi
  i\lambda \tau}e(z^\lambda )$. It follows that $a_\lambda =0$ for
$\lambda \neq 0$. By (1), $\xi \in K$.\hfill $\square$

\smallskip
\noindent (3) {\it Let $\xi \in Qt(Univ)^H$. Then 
$I:=\{a\in Univ\ |\ a\xi \in Univ\}$ is the unit ideal}. 

\smallskip 
\noindent (3) implies $\xi \in Univ$ and, by (2), $\xi \in K$. This
finishes the proof of 5.1.

\smallskip 

\noindent {\it Proof of }(3). The ideal $I$ is invariant under $H$. Write any
element 
$a\in I$ as $\sum _{\lambda \in {\mathbb Q}}a_\lambda e(z^\lambda )$ with
all $a_\lambda \in K[\{e(c)\},\ell ]$. Then 
$D^m(a)=\sum a_\lambda e^{2\pi i\tau m\lambda}\cdot e(z^\lambda )\in I$ for
every $m\geq 0$. 
It follows that all $a_\lambda \in I$. Thus $I$ is generated by an ideal
$J\subset K[\{e(c)\},\ell ]$ that is also $H$-invariant. Let $m\geq 0$ be
minimal such that $J$ contains a non zero element of degree $m$ in $\ell$.
Let $J_m\subset K[\{e(c)\}]$ denote the ideal of the coefficients of $\ell ^m$
of the elements in $J$ with degree $\leq m$. Then $J_m$ is again $H$-invariant.
Consider a non zero element $A=\sum _{c\in R}a_ce(c)\in J_m$ with all $a_c\in
K$. From $\Gamma ^nA\in J_m$ and $\Delta ^nA\in J_m$ for all $n\geq 0$ one
concludes that each $e(c)$ with $a_c\neq 0$ lies in $J_m$. Hence $J_m$ is the
unit ideal. If $m=0$, then the proof of (2) is finished. If $m>0$, then,
by construction, $J$ contains an element $B$ of degree $m$ in $\ell$
and with leading coefficient 1. Let $B'$ be another element having these
properties. Then $B-B'=0$ since its degree in $\ell$ is $<m$. Thus $B$ is
unique and therefore invariant under
$\Gamma$ and $\Delta$. Using (1), one obtains the contradiction $B\in K$.
\end{proof}

The proof of Corollary 5.2 is similar to that of Proposition 5.1. 

\begin{corollary} The universal ring $Univ_{rs}$ for the category of
regular singular $q$-difference modules over $K$ is
$K[\{e(c)\},\ell ]$. The restrictions $\Gamma _{rs},\Delta _{rs}$ of 
$\Gamma ,\Delta$ to the universal ring $Univ_{rs}$ commute and are topological
generators for the group ${\mathbb G}_{rs}$
of all automorphisms of the difference ring $Univ_{rs}$.
\end{corollary}

We introduce some Tannakian categories. The first one $RegSing_K$ consists
of the regular singular difference modules over $K$. 

The second one, called $Tuples_1$, has as objects the tuples  $(V,\Gamma
_V,\Delta _V)$ where $V$ is a finite dimensional vector space over ${\mathbb
  C}$ and $\Gamma _V,\Delta _V$ are commuting elements in ${\rm GL}(V)$ and
such that all their eigenvalues have absolute value 1. A morphism
$(V,\Gamma _V,\Delta _V)\rightarrow (W,\Gamma _W,\Delta _W)$ is 
a linear map $f:V\rightarrow W$ satisfying $f\circ \Gamma _V=\Gamma _W\circ f$
and $f\circ \Delta _V=\Delta _W\circ f$.

The third one, called $Unitcirc$, consists of the finite dimensional complex
representations $\rho$ of the group ${\mathbb Z}^2$ such that every element 
in the image of $\rho$ has all its {\it eigenvalues on the unit circle}
$S^1:=\{z\in {\mathbb C}|\ |z|=1\}$.

\begin{lemma} There are equivalences of Tannakian categories
\[ RegSing_K\rightarrow Tuples_1\rightarrow Unitcirc \ . \] 
\end{lemma}
\begin{proof} To a regular singular difference module $M$ over $K$ one
associates its `solution space' $V:={\ker}(\Phi -1,Univ_{rs}\otimes _KM)$. 
Then $V$ is a complex vector space and the canonical map
$Univ_{rs}\otimes _{\mathbb C}V\rightarrow Univ_{rs}\otimes _KM$ is bijective.
The group ${\mathbb G}_{rs}$ acts on $V$. Let $\Gamma _V,\ \Delta _V$ denote
the actions of $\Gamma ,\ \Delta $ on $V$. The eigenvalues of
$\Gamma _V,\ \Delta _V$ have absolute value 1, according to the definition
of $\Gamma$ and $\Delta$. Then we associate to $M$ the tuple $(V,\Gamma
_V,\Delta _V)$.

To a tuple $(V,\Gamma _V,\Delta _V)$ we associate the representation
$\rho :{\mathbb Z}^2\rightarrow {\rm GL}(V)$ with $\rho (1,0)=\Gamma _V$
and $\rho (0,1)=\Delta _V$. 

It is easily verified that the functors, defined above, are equivalences 
of tannakian categories.
\end{proof}

We recal that $E_q={\mathbb C}^*/q^{\mathbb Z}$ and consider the
maps 
\[can:{\mathbb C}_u\stackrel{u\mapsto e^{2\pi iu}}{\rightarrow} {\mathbb C}_z^*
\stackrel{pr}{\rightarrow} E_q \ .\]
The indices $u$ and $z$ denote the global variables for these spaces. 
The map $pr$ is the obvious map. The map $can$ is the universal covering. Its
kernel ${\mathbb Z}+{\mathbb Z}\tau$ can be identified with $\pi _1(E_q)$.
Let $a,b$ denote the generators of $\pi _1(E_q)$ corresponding to $1,\tau$ 
(or to the two circles ${\mathbb R}/{\mathbb Z}1$ and 
${\mathbb R}/{\mathbb Z}\tau$ on $E_q$).  

We identify the group  ${\mathbb Z}^2$ (or equivalently the group generated
by $\Gamma _{rs},\ \Delta _{rs}$) with
$\pi _1(E_q)$, by $(1,0)$ (or $\Gamma _{rs}$) is mapped to $a$ and $(0,1)$ 
(or $\Delta _{rs}$) is mapped to $b$. In this way one finds a Tannakian
equivalence between $RegSing_K$ and the category of the `unit circle'
representations of $\pi _1(E_q)$.

To a connection $({\mathcal M},\nabla )$ on $E_q$ one associates its monodromy
representation, i.e., a representation of $\pi _1(E_q)$. In this way
the category of the connections on $E_q$ is equivalent to the category of the
representations of $\pi _1(E_q)$. Combining this with Lemma 5.3 one obtains
 
\begin{corollary} The category of the regular singular difference modules
over $K$ is Tannakian equivalent to a full subcategory of connections on $E_q$.
\end{corollary}

In section 6, this equivalence will be made explicit. This result is, in a
different form, present in the work of J.~Sauloy, see [R-S-Z].

\subsection{Galois groups for split modules over $K$}

Now we extend the above to larger categories.
We introduce a category $Tuples_2$ as follows. The objects are tuples 
\[(V,\{V_\lambda \}_{\lambda \in {\mathbb Q}},\Gamma _V,\Delta _V),
\mbox{ satisfying :}\]
(a) $V$ is a finite dimensional vector space over $\mathbb C$.\\
(b) $V$ is the direct sum of the subspaces $V_\lambda$ (and thus there
are only finitely many $\lambda$ with $V_\lambda \neq 0$).\\
(c) $\Gamma _V,\Delta _V$ are invertible operators on $V$ respecting
the direct sum decomposition and such that all their eigenvalues are on the
unit circle $S^1=\{z\in {\mathbb C}|\ |z|=1\}$. \\
(d) $\Gamma _V\Delta _V\Gamma _V^{-1}\Delta ^{-1}$ acts on each
$V_\lambda$ as multiplication by $e^{2\pi i\lambda}$.

\smallskip
A morphism $(V,\{V_\lambda \},\Gamma _V, \Delta _V)\rightarrow
(W,\{W_\lambda \},\Gamma _W,\Delta _W)$ is a linear map $f:V\rightarrow W$
which respects the direct sum decompositions and satisfies
$f\circ \Gamma _V=\Gamma _W\circ f,\ f\circ \Delta _V=\Delta _W\circ f$.
The tensor product of two objects 
$(V,\{V_\lambda \},\Gamma _V, \Delta _V)$ and
$(W,\{W_\lambda \},\Gamma _W,\Delta _W)$ is $(U,\{U_\lambda\},\Gamma _U,
\Delta _U)$, given by $U=V\otimes W,\ U_\lambda =
\sum _{\mu _1+\mu _2=\lambda}V_{\mu _1}\otimes W_{\mu _2}$ and
$\Gamma _U=\Gamma _V\otimes \Gamma _W,\ \Delta _U=\Delta _V\otimes \Delta _W$.
It is easily seen that $Tuples_2$ is a Tannakian category. Further, the Galois 
group of an object $(V,\{V_\lambda \},\Gamma _V, \Delta _V)$ is the algebraic
subgroup of ${\rm GL}(V)$ generated by the maps $\Gamma _V,\Delta _V, D_V$,
where the last map is defined by: $D_V$ is multiplication by 
$e^{2\pi i\tau \lambda} $ on the direct summands $V_\lambda$ of $V$. 

\smallskip

One defines a functor $\mathcal F$ from the category of the split
$q$-difference modules over $K$ to the category $Tuples_2$ as follows.
For a module $M$, ${\mathcal F}(M)=(V,\{V_\lambda \},\Gamma _V,\Delta _V)$,
where\\ 
(i) $V:={\rm ker}(\Phi -1,\ Univ\otimes _KM)$.\\
(ii) $V_\lambda :={\rm ker}(\Phi -1,\ Univ_\lambda \otimes _KM)$, where
$Univ_\lambda :=K[\{e(c)\},\ell ]e(z^\lambda )$.\\
(iii) $\Gamma _V, \Delta _V$ are induced by the action 
of $\Gamma, \Delta$ on $Univ$ and $Univ\otimes _KM$. Since $\Phi -1$ commutes
with $\Gamma, \Delta$, the latter maps leave $V$ invariant.\\
The definition of ${\mathcal F}(f)$, for a morphism $f$, is obvious. The
verification that $\mathcal F$ is a functor between Tannakian categories
is straightforward.

\begin{proposition} $\mathcal F$ is an equivalence between the Tannakian
categories of the split difference modules over $K$ and $Tuples_2$.
\end{proposition}
\begin{proof} One defines a functor ${\mathcal G}$ from the category 
$Tuples_2$ to the category of the split difference modules over $K$, in the
following way. The image $M$ of an object 
$(V,\{V_\lambda \},\Gamma _V, \Delta _V)$ is the set of the 
${\mathbb G}_{univ}$-invariant elements (or the elements invariant under 
$\Gamma ,\Delta, D$) of $Univ\otimes _{\mathbb C}V$. The action on the last
object is defined by $\Gamma (u\otimes v)=\Gamma (u)\otimes \Gamma _V(v)$ for
$u\in Univ,\ v\in V$ and similarly for $\Delta$ and $D$. The proof that 
$\mathcal G$ is the `inverse' of $\mathcal F$ is similar to the proof of
Proposition 5.1. \end{proof}

The (difference) Galois group of a split module $M$ over $K$ coincides
with that of the object ${\mathcal F}(M)$.
We note that $Tuples_2$ is equivalent to the category of all difference
modules over $\widehat{K}$, too.

\begin{example}{\rm
 Let $M=Res(E(cz^{t/n}))$ with $n>1$, $g.c.d.(t,n)=1$
and $c=e^{2\pi i(a_0+a_1\tau)}$ with $a_0,a_1\in {\mathbb R}$ and
$0\leq a_0,a_1<1$. Write $M$ as $K_ne$ with $\Phi e=cz^{t/n}e$. Then 
$\{z^{j/n}e|\ j=0,\dots ,n-1\}$ is a basis of $M$ over $K$.
To obtain a basis for 
$V=V_{t/n}={\rm ker}(\Phi -1, Univ_{t/n}\otimes M)$ one observes that $V$ is 
in fact
${\rm ker}(\Phi -1,\ K[e(e^{-2\pi i\tau /n}),e(e^{2\pi i
 /n})]e(cz^{t/n})\otimes _KK_ne)$. 
After identifying $z^{1/n}$ with $e(e^{-2\pi i\tau  /n})$, 
this space
takes the form $K_n[e(e^{2\pi i/n})]e(cz^{t/n})\otimes K_ne$ and the 
$\Phi$-action is given by
$\Phi (ae(cz^{t/n})\otimes be)=\phi (a)z^{-t/n}e(cz^{t/n})\otimes 
\phi (b)z^{t/n}e$. From this expression one finds a basis 
$\{v_j|\ j=0,\dots ,n-1\}$ of $V$ of the form
\[v_j=\sum _{s=0}^{n-1} \Phi ^s\{e(e^{2\pi ij/n})e(cz^{t/n})\otimes e\}\ . \] 
We use the `cyclic' notation $v_{j+n}=v_j$. 
The actions of $\Gamma$ and $\Delta$ are given by 
\[\Gamma v_j=e^{2\pi ia_1}v_{j+t}\mbox{ and }
\Delta v_j=e^{-2\pi i(a_0+j/n)}v_j \ .\]
The difference Galois group $G$ is topologically generated by $\Gamma ,\Delta$
and $D$ (the latter is here the multiplication 
by $e^{2\pi i\tau t/n}$). One finds
the exact sequence $1\rightarrow G^o\rightarrow G\rightarrow 
({\mathbb Z}/n{\mathbb Z})^2\rightarrow 0$ with $G^o={\mathbb G}_m={\mathbb
  C}^*id $. Further $\Gamma ^n=e^{2\pi ina_1}id$, 
$\Delta ^n=e^{2\pi i na_0}id$, $\Gamma \Delta \Gamma ^{-1}\Delta ^{-1}=
e^{2\pi it/n}id$ and $G/G^o$ is generated by the images of $\Gamma$ and
$\Delta$. }\end{example}

\subsection{Galois groups for modules over $K$ with two slopes}
Let the $q$-difference module $M$ over $K$ be given by an exact sequence
$0\rightarrow P_1\rightarrow M\rightarrow P_2\rightarrow 0$, where $P_1,P_2$
are pure with slopes $\lambda _1 < \lambda _2$. The Picard-Vessiot ring 
$PV(M)$ of $M$ contains $PV(P_1\oplus P_2)$, the Picard-Vessiot ring of 
$P_1\oplus P_2$. Thus there is an exact sequence
$1\rightarrow H\rightarrow G\rightarrow G'\rightarrow 1$, where $G$ and
$G'$ are the difference Galois groups of $M$ and $P_1\oplus P_2$.
The group $G'$ has been calculated in subsections 5.1 and 5.2.
Further, $H$ is the difference Galois group of the extension
$PV(P_1\oplus P_2)\subset PV(M)$. 

Let $V_1,V_2,V$ denote the solution spaces
for $P_1,P_2,M$. There is an obvious exact sequence $0\rightarrow
V_1\rightarrow V\rightarrow V_2\rightarrow 0$.  The space $V_1$ is invariant
under the action of $G$ on $V$ and the induced action of $G$ on
$gr(V):=V_1\oplus V_2$ coincides with $G'$. Hence $H$ can be identified with
a ${\mathbb C}$-vector space of linear maps from $V_2$ to $V_1$. Indeed, 
for $h\in H$, the kernel of $h-1$ contains $V_1$ and the image of $h-1$
lies in $V_1$. In other words, $H$ can be identified with a linear subspace
$W$ of $V_2^*\otimes _{\mathbb C} V_1$, which is the solution space of
$P_2^*\otimes _KP_1$.

\smallskip

For any difference module $N$ over $K$ we write $\widehat{N}$ for
the difference module $\widehat{K}\otimes _KN$ over $\widehat{K}$.
Consider the split exact sequence  
\[0\rightarrow \widehat{P}_1\rightarrow \widehat{M}\rightarrow \widehat{P}_2
\rightarrow 0 .\] 
The splitting $\widehat{P}_2\rightarrow \widehat{M}$ is unique since
the only morphism $\widehat{P}_2\rightarrow \widehat{P}_1$ of 
$q$-difference modules over $\widehat{K}$ is the zero map. Moreover,
the canonical morphism $\widehat{K}\otimes PV(P_1\oplus P_2)\rightarrow
PV(\widehat{M})$ is an isomorphism. The resulting embedding 
$PV(P_1\oplus P_2)\subset PV(\widehat{M})$ extends in a unique way
to an embedding $PV(M)\subset PV(\widehat{M})$. The solution space 
$V$ of $M$ can be  defined as ${\rm ker}(\Phi -1,\widehat{K}\otimes _K
Univ \otimes _KM)$. Similar expressions are valid for $V_1,V_2$, the solution
spaces of $P_1,P_2$. Write $G''$ for the difference Galois group of 
$PV(\widehat{M})$. The restriction map induces an isomorphism 
$G''\rightarrow G'$. We conclude that the exact sequences
\[1\rightarrow H\rightarrow G\rightarrow G'\rightarrow 1\mbox{ and }
0\rightarrow V_1\rightarrow V\rightarrow V_2\rightarrow 0\] 
have canonical splittings. Write $V=V_1\oplus V_2$ for this splitting, then
$G'$ (or $G''$) can be identified with the subgroup of $G$ leaving both
$V_1$ and $V_2$ invariant. Thus $G$ is the semi-direct product $H\rtimes G'$.
The exact sequence, defining $M$, is given by an 
element $\xi \in {\rm coker}(\Phi -1,{\rm Hom}(P_2,P_1))$. The main issue
is to derive $H$ (or $W$) from the given $\xi$ and to formulate $\xi$ in terms
of the solution spaces $V_1,V_2$.

\medskip

Let $N$ be a pure difference module with slope $\lambda <0$ and solution
space $V(N)$. Our first aim is to formulate ${\rm coker}(\Phi -1,N)$ in terms 
of $V(N)$.

\begin{lemma} Let $N$ be a pure module with slope $\lambda <0$. The canonical 
map $a\mapsto 1\otimes a$ from $N$ to $Univ\otimes _KN$ induces an injective
map  
\[{\rm coker}(\Phi -1,N)\rightarrow {\rm coker}(\Phi -1, Univ\otimes _KN)\]
whose image consists of the elements which are invariant under the
universal group ${\mathbb G}_{univ}$. 
\end{lemma}
\begin{proof} 
 It suffices to consider a pure, indecomposable 
$N=Res(E(dz^{t/n}))\otimes U_m$. By writing $N$ as an extension of
$Res(E(dz^{t/n}))\otimes U_{m-1}$ by $Res(E(dz^{t/n}))$, one is reduced to
verify the lemma for $Res(E(dz^{t/n}))$. 
For notational reasons we write 
${\mathbb R}=L\oplus {\mathbb Q}$ for some ${\mathbb Q}$-vector space $L$.
Let $C\subset {\mathbb C}^*$ denote the elements $c$ of the form 
$e^{2\pi i(a_0+a_1\tau )}$ with $a_1\in L$. Then we can write
$Univ=\oplus _{c\in C,\lambda \in {\mathbb Q}}K_{\infty}[\ell ]e(cz^\lambda )$.
Then $Univ\otimes _KN =\oplus K_{infty}[\ell ]e(cz^\lambda )\otimes _KN$. Each
direct summand is invariant under the action of $\Phi$ and 
${\mathbb G}_{univ}$. The cokernel of $\Phi -1$ of each direct summand has
a certain ${\mathbb G}_{univ}$-action. Using the topological generators
$\Gamma ,\Delta ,D$ of this group, one finds that only the direct summand 
$K_{\infty}[\ell ]\otimes _K N$ can produce ${\mathbb G}_{univ}$-invariant
elements in the cokernel of $\Phi -1$. A further inspection shows that
this contribution comes from the subspace $K\otimes _KN$. Thus we find the
required bijection. \end{proof}

We observe that $Univ\otimes _KN=Univ \otimes _{\mathbb C}V(N)$, where
the solution space $V(N)$ of $N$ is defined as 
${\rm ker}(\Phi -1,Univ\otimes _KN)$.
The ${\mathbb G}_{univ}$-invariant part of ${\rm coker}(\Phi -1, Univ 
\otimes _{\mathbb C}V(N))$, comes from 
$K[\{e(c)\},\ell ]e(z^{\lambda})\otimes _{\mathbb C}V(N)$, where $D$ acts as 
the identity. We conclude the following.

\begin{lemma} Let $N$ be a pure module with slope $\lambda <0$ and solution 
space $V(N)$. Then ${\rm coker}(\Phi -1,N)$ can be identified with the 
$\mathbb C$-linear subspace of  
\[({\rm coker}(\phi -1,K[\{e(c)\},\ell ]e(z^{\lambda}))\otimes _{\mathbb C}
V(N),\]
consisting of the elements invariant under $\Gamma$ and $\Delta$.
\end{lemma}

\begin{lemma} Let $\lambda <0$. The subspace 
$H(\lambda ):=
\{\oplus _{|q|^{-\lambda}<|c|\leq 1}{\mathbb C}e(cz^\lambda )\}[\ell ]$
of $K[\{e(c)\},\ell ]e(z^\lambda )$ has the following properties:\\
$H(\lambda )\rightarrow {\rm coker}(\phi -1,K[\{e(c)\},\ell ]e(z^\lambda ))$
is bijective and\\
$H(\lambda )$ is invariant under the actions of $\Gamma, \Delta$ and $D$.
\end{lemma}
\begin{proof} From Theorem 2.5 it follows that we may, for the calculation
of this cokernel, replace $K$ by ${\mathbb C}[z,z^{-1}]$. After replacing
$z$ by $e(q^{-1})$, we have to compute the coker of $\phi -1$ on the space
$\{\oplus _{c\neq 0}{\mathbb C}e(cz^\lambda )\}[\ell ]$. Using 
2.6 one finds that $H(\lambda )$ has the first property.
An inspection of the actions of $\Gamma ,\Delta , D$ yields the second
property. \end{proof}

\begin{corollary} The canonical map
${\rm coker}(\Phi -1,N)\rightarrow 
(H(\lambda )\otimes V(N))^{<\Gamma ,\Delta >}$ {\rm (}i.e., the elements
invariant under the group generated by $\Gamma$ and $\Delta${\rm )} is an 
isomorphism. 
\end{corollary}

\begin{corollary} The exact sequences $0\rightarrow P_1\rightarrow M
\rightarrow P_2\rightarrow 0$ {\rm (}with $P_1,P_2$ pure modules with
slopes $\lambda _1<\lambda _2$ and solution spaces $V_1,V_2${\rm )} are in
bijection with the elements of $(H(\lambda _1-\lambda _2)\otimes 
{\rm Hom}(V_2,V_1))^{<\Gamma ,\Delta >}$.
\end{corollary}

\begin{theorem} Let the difference module $M$ be given by an
exact sequence $0\rightarrow P_1\rightarrow M\rightarrow P_2\rightarrow 0$
with $P_1,P_2$ pure modules with slopes $\lambda _1< \lambda _2$.
Let $V_1,V_2$ denote the solution spaces of $P_1,P_2$ and write $V=V_1\oplus
V_2$ for the solution space of $M$.

Let $\xi \in {\rm coker}(\Phi -1,{\rm Hom}(P_2,P_1))$ and its image\\
 $\xi ' \in (H(\lambda _1-\lambda _2)\otimes 
{\rm Hom}(V_2,V_1))^{<\Gamma ,\Delta >}$
represent the exact sequence.

Define $W\subset {\rm Hom}(V_2,V_1)$ to be the smallest subspace such that\\
$\xi ' \in H(\lambda _1-\lambda _2)\otimes _{\mathbb C}W$. Then:

\smallskip

The difference Galois group $G\subset {\rm GL}(V)$ of $M$ is the semi-direct 
product $G=H\rtimes G'$ with: $G'=\{g\in G| g(V_i)=V_i\mbox{ for }i=1,2\}$ is
the difference Galois group of $P_1\oplus P_2$ and 
$H$ consists of the linear maps of the form $id+\tilde{w}$ where 
$\tilde{w}$ is equal to
$V\stackrel{pr_2}{\rightarrow}V_2\stackrel{w}{\rightarrow}V_1\subset V$ with 
$w\in W$.
\end{theorem}
\begin{proof} It suffices to prove the last statement. For any submodule
$N$ of ${\rm Hom}(P_2,P_1)$ one considers the exact sequence 
$0\rightarrow N\rightarrow {\rm Hom}(P_2,P_1)\rightarrow N'\rightarrow 0$.
This gives rise to an exact sequence
\[0\rightarrow {\rm coker}(\Phi -1,N)\rightarrow 
{\rm coker}(\Phi -1,{\rm Hom}(P_2,P_1))\rightarrow 
{\rm coker}(\Phi -1,N')\rightarrow 0\ .\]
It follows that there exists a smallest submodule $N_0$ such that $\xi$ lies
in ${\rm coker}(\Phi -1,N_0)$. The solution space $V(N_0)$ lies in 
${\rm Hom}(V_2,V_1)$ and  $\xi '$ belongs to
$H(\lambda _1-\lambda _2)\otimes V(N_0)$ (and is invariant under $<\Gamma
,\Delta >$). 

It is clear that $H$ identifies with $\{id +\tilde{w}|\ w\in W_1\}$ for some 
$W_1\subset {\rm Hom}(V_2,V_1)$ which is invariant under the action of the
difference Galois group of ${\rm Hom}(P_2,P_1)$, acting upon its solution
space ${\rm Hom}(V_2,V_1)$. In other words, $W_1$ is invariant under
$<\Gamma ,\Delta >$. By Tannakian correspondence, $W_1$ corresponds
to a submodule $N_1$ of ${\rm Hom}(P_2,P_1)$. Moreover, $\xi$ lies in
${\rm coker}(\Phi -1,N_1)$.   

One concludes that $N_0=N_1$ and that $W=W_1$.
\end{proof}

The elements $\{e(cz^\lambda )\ell ^m\}_{|q|^{-\lambda }<|c|\leq 1;\ m\geq 0}$,
with $\lambda :=\lambda _1-\lambda _2$, form a $\mathbb C$-basis of
$H(\lambda _1-\lambda _2)$. Thus $\xi '$ can be written uniquely as
$\sum e(cz^\lambda )\ell ^m\otimes w(c,m)$ with all $w(c,m)\in 
{\rm Hom}(V_2,V_1)$. Then $W$ is the subspace of ${\rm Hom}(V_2,V_1)$ generated
by $\{w(c,m)\}$. This observation makes it possible to compute the difference
Galois group for explicit examples.

\subsection{Modules over $K$ with more slopes}
We describe here the difference Galois group for the general case. The proof
follows straightforward from the case of two slopes.
Consider a $q$-difference module $M$ over $K$ with slope filtration
$0=M_0\subset M_1\subset \cdots \subset M_r=M$. Write 
$gr(M)=P_1\oplus \cdots \oplus P_r$ where $P_i$ is pure of slope 
$\lambda _i$ and $\lambda _1<\cdots <\lambda _r$. According to Remarks 2.5,
part (1), $M$ is given by an element $\xi =\{\xi _{i,j}\}$ in 
$\prod _{i<j}{\mathcal F}_{i,j}^o({\mathbb C})\cong \prod _{i<j}{\rm coker}
(\Phi -1,{\rm Hom}(P_j,P_i))$. Using the unique direct
sum decomposition of $\widehat{K}\otimes M=\oplus \widehat{K}\otimes P_i$
one finds a canonical decomposition $V=V_1\oplus \cdots \oplus V_r$ of the
solution space $V$ of $M$, where $V_i$ stands for the solution space of $P_i$.
Moreover the difference Galois group $G\subset {\rm GL}(V)$ of $M$ is a
semi-direct product $G=H\rtimes G'$. The group $G'$ is the difference Galois
group of $P_1\oplus \cdots \oplus P_r$ and consists of the $g\in G$ that
leave each subspace $V_i$ invariant. The normal subgroup $H$ is generated by
subgroups $H_{i,j}$ with $i<j$. This subgroup $H_{i,j}$ consists of
the maps $id +\tilde{w}_{i,j}$ with $\tilde{w}_{i,j}=
V\stackrel{pr _j}{\rightarrow}V_j\stackrel{w_{i,j}}{\rightarrow} V_i\subset V$ 
with $w_{i,j}\in W_{i,j}\subset {\rm Hom}(V_j,V_i)$. Further $W_{i,j}$
is the smallest $\mathbb C$-linear subspace of ${\rm Hom}(V_j,V_i)$ such that
the element $\xi _{i,j}\in {\rm coker}(\Phi -1,{\rm Hom}(P_j,P_i))$ is
represented by $\xi _{i,j}'\in H(\lambda _i-\lambda _j)\otimes W_{i,j}$. 
According to Remarks 3.3, part (3), $H$ is in fact equal to set of elements
$id +\sum _{i<j}\tilde{w}_{i,j}$ with all $w_{i,j}\in W_{i,j}$. Indeed,
$w_{i,j}\circ w_{j,k}\in W_{i,k}$ for $i<j<k$ and any $w_{i,j}\in W_{i,j}$, 
$w_{j,k}\in W_{j,k}$.

\medskip

The Tannakian category of difference modules over $K$ has the following 
description in terms of solution spaces. One can attach to a difference
module $M$ the tuple
$(V,\{V_\lambda \},\Gamma _V ,\Delta _V,\{\xi _{\lambda ,\mu}\})$,
where $(V,\{V_\lambda \},\Gamma _V,\Delta _V)$ is the tuple associated to
the split module $gr(M)$ and $\xi _{\lambda ,\mu}$ is, for each
$\lambda <\mu $, an element of $H(\lambda -\mu )\otimes {\rm Hom}(V_\mu
,V_\lambda )$, invariant under the action of $<\Gamma _V,\Delta _V>$. The
category of the above tuples is, in an obvious way, a Tannakian category. The
data $\{\xi _{\lambda ,\mu}\}$ come from divergent solutions (i.e., with
coefficients in $\widehat{K}$). They can be seen as the equivalent of the 
Stokes matrices in the theory of irregular differential equations over 
$\mathbb C$.

\section{Realizing $q$-difference modules over $K$ as 
connections on the elliptic curve $E={\mathbb C}^*/q^{\mathbb Z}$}

As in section 2 we associate to any difference module $M$ over $K$ a vector
bundle $v(M)$ on the elliptic curve $E_q$. The aim is to provide $v(M)$
with a suitable connection. The following theorem gives an explicit version
of Corollary 5.4.

\begin{theorem}[Regular singular $q$-difference modules] $\;$\\
Let $a,b$ denote the generators of $\pi _1(E_q)$, corresponding to the shifts
$u\mapsto u+1,\ u\mapsto u+\tau$ on the universal covering ${\mathbb C}_u$ of
$E_q$. Let $i:<\Gamma _{rs},\Delta _{rs}>\rightarrow \pi _1(E_q)$ denote the
isomorphism given by $\Gamma _{rs}\mapsto a,\ \Delta _{rs} \mapsto b$.

For every regular singular difference module $M$ over $K$, corresponding
to a representation $\rho$ of $<\Gamma _{rs},\Delta _{rs}>$,
there exists a unique `unit circle' connection $\nabla _M$ on $v(M)$,
corresponding to a representation $\rho '$ of $\pi _1(E_q)$, 
such that $\rho =i\circ \rho '$. 

This induces a Tannakian equivalence between the category of the regular
singular difference modules over $K$ and a full subcategory of the
category of all `unit circle' connections on $E_q$. 
\end{theorem}

\begin{proof} 
We recall that a connection $\nabla :{\mathcal M}\rightarrow \Omega _{E_q}
\otimes {\mathcal M}$ is called {\it unit circle} if the corresponding 
representation $\rho ':\pi _1(E_q)\rightarrow {\rm GL}(V)$ has the property
that every eigenvalue of every $\rho '(\alpha )$ has absolute value 1.\\
(1) We start by giving an explicit construction for rank one difference
modules $M$ over $K$. Write $M=Ke$ with $\Phi e=ce$ and $c\in {\mathbb C}^*$.
A connection $\nabla _M$ on $v(M)$ is equivalent with a connection
$\nabla : Oe\rightarrow \frac{dz}{z}\otimes Oe$ commuting
with $\Phi$ (we recall that $O$ is the algebra of holomorphic functions on
${\mathbb C}^*$).
One concludes that $\nabla (e)=a(c)\frac{dz}{z}\otimes e$, where
$a:{\mathbb C}^*\rightarrow {\mathbb C}$ is a homomorphism. 
The differential equation 
$\nabla e=a(c)\frac{dz}{z}\otimes e$, considered on ${\mathbb C}$, reads
$\nabla e=2\pi i a(c)du\otimes e$ where $u$ is the variable of ${\mathbb C}$
and $z=e^{2\pi i u}$. The basic solution of this equation is
$e^{-2\pi i a(c)u}e$ and $e$ can be interpreted as 
$\theta _{c^{-1}}(e^{2\pi i u})$. The shifts 
$u\mapsto u+1$ and $u\mapsto u+\tau$
multiply this equation with $e^{-2\pi i a(c)}$ and 
$e^{-2\pi i a(c)\tau }c$. Write $c=e^{2\pi i (a_0+a_1\tau )}$ with
$a_0=a_0(c),a_1=a_1(c)\in {\mathbb R}$. The `unit circle' condition implies 
that 
$a(c)=a_1(c)$ and the basic solution is multiplied by $e^{-2\pi i a_1}$ and
$e^{2\pi i a_0}$. Thus we found a unique unit circle connection on $v(Ke)$
and moreover the actions of $\Gamma _{rs},\Delta _{rs}$ coincide with the
action of the fundamental group $\pi _1(E_q)$.     

\smallskip

\noindent (2) 
The next case to consider is $M=K\otimes W$ with $\Phi (f\otimes w)=
\phi(f)\otimes A(w)$ and $A\in {\rm GL}(W)$ is unipotent. Define
$\nabla :O\otimes W\rightarrow O\frac{dz}{z}\otimes W$ by $\nabla (1\otimes w)=
\frac{1}{2\pi i \tau}\frac{dz}{z}\otimes (1\otimes \log (A)(w))$. Clearly
$\nabla$ commutes with $\Phi$. One can verify that this definition
of $\nabla$ induces a regular connection on $v(M)$. The pullback of this 
connection to $\mathbb C$ has the fundamental matrix $A^{u/\tau}$. The shifts
$u\mapsto u+1,\ u\mapsto u+\tau$ multiply this matrix with $A^{1/\tau}$ and
$A$, in accordance with the actions of $\Gamma _{rs}$ and $\Delta _{rs}$. 
 
The general case is obtained from these two special cases by taking 
tensor products and direct sums.

\smallskip

\noindent (3)
The map ${\rm Hom}(M,N)\rightarrow {\rm Hom}((v(M),\nabla _M),(v(N),\nabla
_N))$ is a bijection since both groups translate into homomorphisms between
representations.\end{proof}

\begin{remarks} $\;$\\ {\rm
(1) We note that not every unit circle connection on $E$ is isomorphic to
some $(v(M),\nabla _M)$. E.g., the free vector bundle 
$O_Ee_1\oplus O_Ee_2$, provided with the unit circle 
connection $\nabla$ given by 
$\nabla e_1=0,\ \nabla e_2=\frac{dz}{z}\otimes e_1$, is not isomorphic to
some $(v(M),\nabla _M)$. \\   
(2) Our aim is to extend Theorem 6.1 to the category of all difference
modules over $K$. The first step is to extend Theorem 6.1 to the category of
the split difference modules. Explicitly, one wants to construct for every
split difference module $M$ a connection $\nabla _M$ on $v(M)$ such that
$\nabla _M :v(M)\rightarrow \Omega ([1])\otimes v(M)$ (i.e., $\nabla _M$
has a regular singularity at $1\in E_q$). Moreover, one requires that
$M\mapsto \nabla _M$ is functorial, commutes with direct sums and tensor
products, and extends the construction in Theorem 6.1 of $\nabla _M$ for
regular singular difference modules.\\
(3) The pure module $M=Ke$ with $\Phi e=(-z)e$ is the first candidate for a
construction of $\nabla _M:v(M)\rightarrow \Omega ([1])\otimes v(M)$.
We note that $v(M)=O_{E_q}([1])e$. Any
connection $\nabla$ on this line bundle with at most a regular singularity
at $1$ is given by $\nabla e=a\frac{dz}{z}\otimes e$ for some constant
$a$. The assumption $a=0$ is the natural choice for $\nabla _M$. Thus
$\nabla _M$ is the unique extension of the trivial connection on
$O_{E_q}e$, given by $\nabla e=0$ to a connection 
$\nabla _M:v(M)\rightarrow \Omega ([1])\otimes v(M)$. One can describe this 
connection on $Oe$ by $\nabla e=-\frac{d\Theta}{\Theta}\otimes e$.  This leads
to a formula for $\nabla _M$ for the difference modules 
$M=Ke$ with $\Phi e=c(-z)^te$, namely 
\[\nabla e=(-t\frac{d\Theta}{\Theta}+a_1\frac{dz}{z})\otimes e
\mbox{ where }c=e^{2\pi i(a_0+a_1\tau )}\mbox{ and } a_0,a_1\in {\mathbb R} 
\ .\]
One observes that for difference modules of rank one , the map $M\mapsto
(v(M),\nabla _M)$ respects tensor products. \\
(4) Instead of continuing the method of (3), we will use subsubsection 1.4.1
to give a general construction of $\nabla _M$ for split difference modules
over $K$.}\end{remarks}

\begin{theorem} There exists a functor $N\mapsto (v(N),\nabla _N)$ from the
category of the split difference modules over $K$ to the category of the
connections on $E_q$ with at most a regular singularity at $1\in E_q$. This
functor extends the one defined in Theorem 6.1. Moreover, the functor 
$N\mapsto (v(N),\nabla _N)$ commutes with tensor products and is faithful.
\end{theorem}
\begin{proof}  We use the notation and the results of subsubsection 1.4.1.

\smallskip

\noindent (1)
{\it Construction of $\nabla _M$ for pure difference modules $M$ 
over $K_\infty$.}\\
Consider a pure difference module $M$ over $K_\infty$ with slope $\lambda$.
Then $M=K_\infty \otimes _{\mathbb C}V$ with $\Phi$ given by
$\Phi (f\otimes v)=z^\lambda \phi (f)\otimes A(v)$ where $
A\in {\rm GL}(V)$ has the
property that every eigenvalue $c$ of $A$ has the form
$c=e^{2\pi i (a_0(c)+a_1(c)\tau )}$ with
$a_0(c)\in \mathbb R$ and $a_1(c)\in L\subset {\mathbb R}$. Let 
$a_1(A_{ss})\in {\rm End}(V)$ be obtained from $A_{ss}$ by replacing every
eigenvalue $c$ of $A_{ss}$ by $a_1(c)$.  We introduce the notation:
$O_n$ is the algebra of the convergent Laurent series in the variable 
$z^{1/n}$ and $O_\infty =\bigcup O_n$.
With these notation one defines the connection 
\[\nabla _M:O_\infty \otimes _{\mathbb C}V\rightarrow 
\frac{1}{\Theta }O_\infty \frac{dz}{z}\otimes  V \mbox{ by }\]
\[\nabla _M (v)=-\lambda \frac{d\Theta}{\Theta}\otimes v+
\frac{dz}{z}\otimes (\; a_1(A_{ss})+\frac{1}{2\pi i\tau}\log (A_u)\; )(v)\ .\]
The last formula is extended by 
$\nabla _M(f\otimes v)=df\otimes v+f\nabla _M(v)$ and by additivity to
a $\nabla _M$ defined on $O_\infty \otimes _{\mathbb C}V$. By construction,
$\nabla _M$ commutes with the action of $\Phi$.

For two pure difference modules $M_i=K_\infty \otimes _{\mathbb C}V_i,\
i=1,2$ over $K_\infty$ with $\Phi$-actions given by the slopes $\lambda _1,
\lambda _2$ and $A_i\in {\rm GL}(V_i),\ i=1,2$, the pure module 
$M_3=M_1\otimes _{K_\infty}M_2$ has the form $K_\infty \otimes _{\mathbb
  C}(V_1\otimes V_2)$ with $\Phi$-action given by the slope 
$\lambda _1+\lambda _2$ and  $A_3:=A_1\otimes A_2 \in {\rm GL}(V_1\otimes
V_2)$. One observes that the eigenvalues $c$ of $A_3$ satisfy again
$a_1(c)\in L\subset {\mathbb R}$. Further $(A_3)_{ss}=(A_1)_{ss}\otimes 
(A_2)_{ss}$ and $a_1((A_1)_{ss}\otimes (A_2)_{ss})$
is equal to $(a_1((A_1)_{ss})\otimes id)+(id\otimes a_1((A_2)_{ss}))$.
There is a similar formula for $\frac{1}{2\pi i \tau}\log (A_3)_u$.

One concludes that $\nabla _{M_3}$ is the tensor product 
$\nabla =\nabla _{M_1}\otimes \nabla _{M_2}$. The latter is defined by
$\nabla (m_1\otimes m_2)=(\nabla _{M_1}(m_1)\otimes m_2)+
(m_1\otimes \nabla _{M_2}m_2)$ (for $m_i\in O_\infty \otimes V_i,\ i=1,2$).

A morphism $f:M_1\rightarrow M_2$ between two difference modules 
$M_i=K_\infty \otimes _{\mathbb C}V_i,\ i=1,2$ 
with the same slope $\lambda$ and $\Phi$ given by
$A_i\in {\rm GL}(V_i),\ i=1,2$ is induced by a $\mathbb C$-linear map
$F:V_1\rightarrow V_2$ satisfying $F\circ A_1=A_2\circ F$. Therefore $f$ 
induces a morphism between the two connections $\nabla _{M_i}$.

\medskip
\noindent
(2) {\it Construction of $(v(N),\nabla _N)$ for a pure difference module $N$
over  $K$}.\\
The connection $\nabla _N :v(N)\rightarrow \Omega _{E_q}([1])\otimes v(N)$, 
that we want to construct, translates into a connection with the same name
\[\nabla _N:O\otimes _{{\mathbb C}[z,z^{-1}]}N_{global}\rightarrow 
\frac{1}{\Theta}O\frac{dz}{z}\otimes_{{\mathbb C}[z,z^{-1}]}N_{global}\ ,\]
which commutes with the action of $\Phi$.

Put $M=K_\infty \otimes _KN$. This difference module is equipped with the data
of $N$, i.e., $data(N)=(\lambda ,V,A,\{D(\sigma )\})$. We note that
$O_\infty \otimes _{\mathbb C}V=O_{\infty}\otimes _{{\mathbb
    C}[z,z^{-1}]}N_{global}$ and that $(O_\infty \otimes _{\mathbb C}V)^{Gal}$
is equal to $O\otimes _{{\mathbb C}[z,z^{-1}]}N_{global}$. 
The $\nabla _M$, constructed above, descends to a $\nabla _N$ for $N$ if and
only if $\nabla _M$ commutes with the $\{D(\sigma )\}$.

The formula $D(\sigma )^{-1}AD(\sigma )=e^{2\pi i \lambda \sigma }A$ implies
that $D(\sigma )^{-1}A_uD(\sigma )=A_u$ and 
$D(\sigma )^{-1}A_{ss}D(\sigma )=e^{2\pi i\lambda \sigma }A_{ss}$.   
The eigenspace of $A_{ss}$ for the eigenvalue $c$ is mapped by
$D(\sigma )$ to the eigenspace for the eigenvalue $ce^{2\pi i\lambda \sigma}$.
>From $a_1(ce^{2\pi i\lambda \sigma})=a_1(c)$ it follows that $D(\sigma )$
leaves every eigenspace invariant of $a_1(A_{ss})$ for the eigenvalues of this
map. Thus $D(\sigma )$ commutes with $a_1(A_{ss})$ and with $A_u$, too. This
implies that $\nabla _M$ commutes the $\{D(\sigma )\}$.

\medskip

\noindent 
(3) {\it $(v(N),\nabla _N)$ for a split difference module $N$ over $K$}.\\
For $N=N_1\oplus \cdots \oplus N_r$ with all $N_i$ pure of slopes
$\lambda _1<\cdots <\lambda _r$ one defines 
$(v(N),\nabla _N):=\oplus _{i=1}^r(v(N_i),\nabla _{N_i})$. 

A morphism between split modules over $K$ is the sum of morphisms between 
pure modules with the same slope. The latter induces a morphism between the 
corresponding connections, according to (1). Thus $N\mapsto (v(N),\nabla _N)$
is a functor. According to (1) and 1.4.1, this functor preserves tensor
products.

For proving `faithful' it suffices to show that, for a pure module $N$, 
the map ${\rm ker}(\Phi -1,N)\rightarrow \{\xi \in H^0(E_q,v(N))|\ \nabla _N\xi
=0\}$ is injective. If the slope of $N$ is not 0, then the left hand side is 0.
If the slope is 0, then by Theorem 6.1, the above map is bijective.
\end{proof}

\noindent {\it Remark}. The following example shows that the functor of 
Theorem 6.3 is not fully faithful. Put $N=Ke$ with $\Phi (e)=(-z)^te$ and 
$t>0$. Then $v(N)=O_{E_q}(t\cdot [1])$ and $\nabla _N$ is induced by $d:
O_{E_q}\rightarrow \Omega _{E_q}$, using the inclusion 
$O_{E_q}\subset O_{E_q}(t\cdot [1])$. In this case, ${\rm ker}(\Phi -1,N)=0$
and $ \{\xi \in H^0(E_q,v(N))|\ \nabla _N\xi =0\}$ has dimension 1.

\bigskip

We want to extend the functor of Theorem 6.3 to the category of all difference
modules over $K$. We start with an example.
\begin{example} $N=Ke_1+Ke_2$ with $\Phi e_1=(-z)^te_1,\ \Phi e_2
  =e_2+pe_1$, $t\in {\mathbb Z},\ t>0$ and $p\in K$.\\
{\rm We may suppose $p\in {\mathbb C}[z,z^{-1}]$. The aim is to produce a
connection 
\[\nabla :Oe_1+Oe_2\rightarrow \frac{1}{\Theta ^t}O\frac{dz}{z}\otimes _O
(Oe_1+Oe_2)\ ,\]
such that $\nabla$ commutes with $\Phi$ and $\nabla$ induces for the pure
module $Ke_1$ and $Ke_2$ the connections of Theorem 6.3.

Consider the inclusion $Oe_1+Oe_2\subset Of_1+Of_2$ with 
$f_1=\Theta ^{-t}e_1,\ f_2=e_2$. Then $\Phi f_1=f_1$ and 
$\Phi f_2=f_2+p\Theta ^tf_1$. We propose $\nabla f_1=0$ and 
$\nabla f_2=\omega \otimes f_1$ with $\omega \in O\frac{dz}{z}$. The condition
$\Phi \nabla =\nabla \Phi$ is equivalent to 
$(\phi -1)(\omega )=d(p\Theta ^t)$.

The equation is solved as follows. There exists $f\in O$ with 
$(\phi -1)(f)=p\Theta ^t-c$ where $c$ is the constant term of $p\Theta ^t$.
This $f$ is unique up to its constant term. Now $\omega :=d(f)$ satisfies
$(\phi -1)(df)=d((\phi -1)(f))=d(p\Theta ^t)$.

Thus $\nabla e_1=t\frac{d\Theta }{\Theta}\otimes e_1$ and 
$\nabla e_2=\Theta ^{-t}\omega \otimes e_1$. Since $\nabla$ commutes with
$\Phi$ this induces a connection 
$\nabla :v(N)\rightarrow \Omega _{E_q}(t[1])\otimes v(N)$ with a pole at
$1\in E_q$ of order at most $t$. Moreover $\nabla$ induces on $v(Ke_1)$
and $v(Ke_2)$ the connections of Theorem 6.3. \hfill $\square$ }\end{example}

In the next example we give a construction of $\nabla _N$ for difference
modules $N$ over $K$ with two integer slopes.
\begin{example} $N=P_1\oplus P_2$ with $P_1,P_2$ global modules with
slopes $t_1<t_2, \ t=t_2-t_1\in \mathbb Z$ and 
$\Phi$ given by $\Phi (p_1+p_2)=\Phi _1 p_1+\ell (\Phi _2p_2)+\Phi _2p_2$ for
some ${\mathbb C}[z,z^{-1}]$-linear map $\ell :P_2\rightarrow P_1$.   

\smallskip

{\rm  Put $Q_1={\mathbb C}[z,z^{-1}]e\otimes P_1$ with $\Phi (e)=(-z)^te$ and 
$t=t_2-t_1$ and put $Q_2=P_2$. Then $O\otimes P_1$ is embedded into
$O\otimes Q_1$ by $p\mapsto \Theta ^te\otimes p$. Now $\Phi$ on
$O\otimes Q_1+O\otimes Q_2$ is given by 
\[\Phi (e\otimes p_1+p_2)=(-z)^te\otimes \Phi _1(p_1)+
\Theta ^te\otimes \ell (\Phi _2p_2)+\Phi _2(p_2) \ . \]
$\nabla $ on $O\otimes Q_1+O\otimes Q_2$ will be given by a formula of the
type: $\nabla (e\otimes p_1)=-t\frac{d\Theta}{\Theta}e\otimes p_1+e\otimes 
\nabla _1p_1$, where $\nabla _i$ are the connections for $P_i$ imposed by
Theorem 6.3. Further $\nabla p_2=\nabla _2p_2+m(p_2)$ where 
$m:O\otimes P_2\rightarrow O\frac{dz}{z}\otimes Q_1$ is a yet unknown 
$O$-linear map. Thus $m$ is an element of 
$O\frac{dz}{z}\otimes _OT$ with 
$T=O\otimes _{{\mathbb C}[z,z^{-1}]}(P_2^*\otimes Q_1)$. We note that
$P_2^*\otimes Q_1$ is a pure global module of slope 0. 

Write $\tilde{\ell }:O\otimes P_2\stackrel{\ell}{\rightarrow}O\otimes P_1
\subset O\otimes Q_1$. Thus $\tilde{l}$ is an element of $T$. The condition
$\Phi \nabla =\nabla \Phi$ is equivalent to 
$\Phi _T (m)-m=\nabla _T(\tilde{\ell })$, where $\Phi _T$ and $\nabla _T$
denote the $\Phi$-action and the connection for $T$.

According to Lemma 6.6, there is a canonical solution $m$ for this equation.
The corresponding $\nabla$ induces a connection 
$\nabla _N:v(N)\rightarrow \Omega _{E_q}(t[1])\otimes v(N)$ with a pole
at $1\in E_q$ of order at most $t$. Further $\nabla _N$ induces the
connections on $v(P_1)$ and $v(P_2)$ prescribed by Theorem 6.3.

The action of $\Phi$ on $O\otimes Q_1+O\otimes Q_2$  can be changed 
into an equivalent one  by adding to $\tilde{\ell }$ an expression
$(\Phi -1)(\xi )$ with $\xi \in T$. This is compatible with the construction
of $\nabla _N$, according to part (1) of Lemma 6.6. After such a change,
one may suppose that $\tilde{\ell}$ maps $Q_2$ into $Q_1$. In this situation 
the connection on $Q_1\oplus Q_2$ is the one prescribed by Theorem 6.3,
according to part (2) of Lemma 6.6 \hfill $\square$ }\end{example}

\begin{lemma}$\;$\\
{\rm (1)} Let $P$ be a pure global module of slope 0. There is a canonical
${\mathbb C}$-linear map $f\in O\otimes P\rightarrow \omega (f)\in
O\frac{dz}{z}\otimes P$ satisfying $(\Phi -1)\omega (f)=\nabla f$.
Further $\omega (\Phi (f))=\Phi (\omega (f))$.\\
{\rm (2)} Let $(Q_1,\Phi _1), (Q_2,\Phi _2)$ denote two pure global modules 
of the same slope and let $\ell :Q_2\rightarrow Q_1$ be a ${\mathbb
  C}[z,z^{-1}]$-linear map. Define the pure module $N=Q_1\oplus Q_2$ by
$\Phi (q_1+q_2)=\Phi _1(q_1)+\Phi _2(q_2)+\ell (\Phi _2(q_2))$. 

The connection on $O\otimes N$, defined by Theorem 6.3 coincides with the
connection $\nabla$ given by the formula $\nabla (q_1+q_2)=
\nabla _1(q_1)+\nabla _2(q_2)+\omega (\ell )(q_2)$.  
\end{lemma}
\begin{proof}  Write $P={\mathbb C}[z,z^{-1}]\otimes _{\mathbb C}W$ with
$\Phi$ defined by $\Phi (w)=A(w)$ and $\nabla (w)=
\frac{dz}{z}\otimes (a_1(A_{ss})+\frac{1}{2\pi i\tau}\log (A_u))(w)$.
For convenience we suppose that the eigenvalues $c$ of $A$ satisfy
$|q|<|c|\leq 1$. Write $f=\sum _{n\in \mathbb Z}z^n\otimes f_n$ and 
$\omega (f)=\sum _{n\in \mathbb Z}z^n\frac{dz}{z}\otimes \omega _n$
with all $f_n,\omega _n\in W$. \[\mbox{ Then } 
f=\sum _nz^n\frac{dz}{z}\otimes \{(n+a_1(A_{ss})+
\frac{1}{2\pi i \tau}\log A_u)(f_n)\}\] \[
\mbox{ and } (\Phi -1)\omega (f) =\sum _nz^n\frac{dz}{z}\otimes (q^nA-1)(\omega
_n)\ .\] 
This produces the equations
$(q^nA-1)\omega _n=(n+a_1(A_{ss})+\frac{1}{2\pi i \tau}\log A_u)f_n$. For 
$n\neq 0$,  the map 
$q^nA-1$ is invertible and there is a unique solution $\omega _n$.
For the equation 
$(A-1)\omega _0=(a_1(A_{ss})+\frac{1}{2\pi i \tau}\log A_u)f_0$ 
we write $W$ as a direct sum $\oplus W_c$, where $W_c$ is the generalized
eigenspace for $A$ and the eigenvalue $c$. For $c\neq 1$, the restriction of
the equation to $W_c$ has a unique solution since $A-1$ is invertible on $W_c$.
On the space $W_1$, the equation reads 
$(A-1)\omega _0^1=(\frac{1}{2\pi i\tau}\log A_u)f_0^1$, where $\omega _0^1,
f_0^1$ denote the components in $W_1$ of $\omega _0$ and $f_0$. On $W_1$
we have $A=A_u$. The canonical solution, that we propose, is given by
$\omega _0^1=\sum _{j=0}^\infty (-1)^j\frac{(A_u-1)^j}{2\pi i 
\tau (j+1)}(f_0^1)$.  

It is clear that $f\mapsto \omega (f)$ is $\mathbb C$-linear. The formula
$\Phi (\omega (f))=\omega (\Phi f)$ follows from the explicit definition of
$\omega (f)$.  The expression
`canonical' means the following. Let $\alpha :P_1\rightarrow P_2$ be a
morphism between global modules of slope 0 and let $f\in O\otimes P_1$.
Then $\alpha$ applied to $\omega (f)$ is equal to $\omega (\alpha (f))$. 

A straightforward calculation shows (2).
\end{proof}

\begin{theorem} There exists a $\mathbb C$-linear functor 
$N\mapsto (v(N),\nabla _N)$ from the 
category of the difference modules $N$ over $K$ with integer slopes to the
category of the connections on $E_q$ with at most a pole at the point
$1\in E_q$. This functor extends the one of Theorem 6.3, is faithful and 
commutes with tensor products.
\end{theorem}
\begin{proof} Consider pure global modules $P_j,\ j=1,\dots ,r$ with
integer slopes $\lambda _1<\cdots <\lambda _r$. Let $\Phi _j$ denote the
action of $\Phi$ on $P_j$. Let ${\mathbb C}[z,z^{-1}]$-linear maps
$\ell _{i,j}:P_j\rightarrow P_i$ for $i<j$, be given. 
Define the global module $N$
with ascending slope filtration by $N=P_1\oplus \cdots \oplus P_r$ and
$\Phi (p_1+\cdots +p_r)=\sum _{j=1}^r(\Phi _j(p_j)+\sum _{i<j}\ell _{i,j}\Phi
_j(p_j))$. We want to construct a canonical connection on $O\otimes N$
commuting with $\Phi$.

Define $Q_j:={\mathbb C}[z,z^{-1}]e_j\otimes P_j$ with 
$\Phi e_j=(-z)^{\lambda _r-\lambda _j}e_j$ for $j=1,\dots ,r$. Let $\Phi ^*_j$
be the action of $\Phi$ on $Q_j$. One embeds $O\otimes P_j$ into 
$O\otimes Q_j$ by $p_j\mapsto \Theta ^{\lambda _r-\lambda _j}e_j\otimes p_j$.
Then $O\otimes N=\oplus O\otimes P_j$ embeds into $O\otimes Q$ with 
$Q=\oplus _{j=1}^rQ_j$. Let 
$\tilde{\ell}_{i,j}:O\otimes Q_j\rightarrow O\otimes Q_i$ be the map
derived from $\ell _{i,j}$. Then  $\Phi$  on $O\otimes Q$ is given by 
\[\Phi (q_1+\cdots +q_r)=\sum _{j=1}^r(\Phi ^*_j(q_j)+
\sum _{i<j}\tilde{\ell} _{i,j}\Phi ^*_j(q_j))\ . \]
With this formula the embedding is $\Phi$-equivariant. On $O\otimes Q$
one wants to define a connection $\nabla$ of the form 
\[\nabla (q_1+\cdots +q_r)=\sum _{j=1}^r\nabla _j(q_j)+
\sum _{j=1}^r\sum _{i<j}m_{i,j}(q_j)\ ,\]
with $O$-linear maps $m_{i,j}:O\otimes Q_j\rightarrow O\frac{dz}{z}\otimes _O
(O\otimes Q_i)$. The condition $\Phi \nabla =\nabla \Phi$ translates into
$(\Phi -1)m_{i,j}=\nabla (\tilde{\ell} _{i,j})$
for all $i<j$. These equations are solved in the canonical way of Lemma 6.6.
The restriction of this $\nabla$ to $O\otimes N$ induce a connection 
$\nabla _N:v(N)\rightarrow \Omega _{E_q}((\lambda _r-\lambda _1)[1])\otimes
v(N)$.  
 
The maps $\ell _{i,j}$ in the definition of $N$ are not unique. They can be
changed by adding maps $(\Phi -1)r_{i,j}$ with $r_{i,j}:P_j\rightarrow P_i$.
It is easily seen that $\nabla _N$ only depends on the equivalence classes
of the $\ell _{i,j}$. Using Lemma 6.6, one shows that the above defines
a $\mathbb C$-linear functor. 

Let $N_1,N_2$ denote two global difference modules with ascending slope
filtration and with integer slopes. The above constuction embeds 
$O\otimes N_i$ into $O\otimes M_i$ (for $i=1,2$) where $O\otimes M_i$ has
only one slope. The connections on $N_1,N_2$ are the restrictions of the
connections on $O\otimes M_i$ prescribed by Theorem 6.3. The above
construction applied to $N_3=N_1\otimes N_2$ embeds $O\otimes N_3$ into
the tensor product of the pure modules $O\otimes M_i$. According to 
Theorem 6.3, the connection on this tensor product is the tensor product
of the connections on the $O\otimes M_i$. We conclude that the functor, 
construction above, respects tensor products. \end{proof}

\begin{remarks}$\;$\\ {\rm
(1) For general difference modules $M$ over $K$ it is also possible to
define a connection $\nabla _N$ on $v(N)$ with at most a pole at $1\in E_q$.
However for non integer slopes there seems not be a canonical choice for
$\nabla _N$. \\ 
(2) In the situation of Example 6.5, the map $\ell :P_2\rightarrow P_1$ is 
responsible for divergence, `Stokes matrices' and the unipotent part of the
difference Galois group of $N$. In general, the connection 
$\nabla _N:v(N)\rightarrow \Omega _{E_q}(t[1])\otimes v(N)$ has a pole of
order $t$. The irregularity of the connection $\nabla _N$ locally at
$1\in E_q$ will produce Stokes matrices (in the classical sense) and 
unipotent elements 
of the local analytic differential Galois group  which depend
on $\ell$. The precise relation remains to be investigated.\\
(3) It is interesting to apply another method to Example 6.4.
For any integer $t>0$ one defines 
$G_t=\sum _{n\in \mathbb Z}(q^t)^{n(n-1)/2} (-z^t)^n$. One observes that
$(-z^t)G_t(qz)=G_t(z)$ and that the set of the zeros of $G_t$ is
$\mu _t\times q^{\mathbb Z}$, where $\mu _t$ is the group of the $t$-th roots
of unity.

Let again $N=Ke_1+Ke_2$ with $\Phi (e_1)=(-z)^te_1,\ \Phi (e_2)=e_2+pe_1$
with $p\in {\mathbb C}[z,z^{-1}]$. Define now $f_1=G_t^{-1}e_1$ and $f_2=e_2$.
Then $Oe_1+Oe_2$ embeds into $Of_1+Of_2$ and $\Phi (f_1)=(-1)^{t-1}f_1$
and $\Phi (f_2)=f_2+pG_tf_1$. The connection $\nabla$ is defined by
$\nabla f_1=0$, $\nabla f_2=\omega \otimes f_2$ where $\omega$ is the
canonical solution of $(\Phi -1)\omega =(-1)^{t-1}d(pG_t)$. The resulting
connection $\nabla _N$ has at most simple poles at the image points of
$\mu _t$ in $E_q$.\\
(4) The variation in (3) on Example 6.4 extends to a functor on the category
of the difference modules over $K$ with integer slopes to the category of the
connections on $E_q$ having at most simple poles in the images of
$\bigcup _{t\geq 1}\mu _t$ in $E_q$. This functor is constructed as in the
proof of Theorem 6.7 and it has again the properties: $\mathbb C$-linear,
faithful and commuting with tensor products.

For this variation on Theorem 6.7 and in connection with Example 6.5,
one observes that $\ell :P_2\rightarrow P_1$ contributes to the poles of
$\nabla _N$ on the image points of $\mu _t$ in $E_q$. Thus $\ell$ contributes
to the monodromy group for the connection $\nabla _N$.   }\end{remarks}

\section{Positive characteristic}

Atiyah's paper makes some excursions to positive characteristic. Here, we 
do the same for $q$-difference equations. We replace the field $\mathbb C$
by a field $\mathbf C$ which is algebraically closed and complete for a non
trivial non archimedean valuation. The case where $\mathbf C$ has 
characteristic 0 (i.e., ${\mathbf C}\supset {\mathbb Q}_p$ 
for some prime $p$)  
is not very interesting since most of the preceeding results can be copied
from the complex case with the help of some rigid analysis.

In this section  
we consider an algebraically closed field $\mathbf C$ of characteristic $p>0$, 
complete with respect to a non trivial valuation. 
Further we choose a $q\in \mathbf C$ with $0<|q|<1$.

Over $K={\mathbf C}(\{z\})$, $\widehat{K}={\mathbf C}((z))$ and ${\mathbf
  C}(z)$  one can define $q$-difference modules and study their
properties. The elliptic curve associated to this is the Tate curve
$E_q:={\mathbf C}^*/q^{\mathbb Z}$, constructed with the help of rigid
analysis. We note that this curve is special in the sense that 
its $j$-invariant is transcendental over ${\mathbb F}_p$ and in
particular $E_q$ is ordinary. We make now a quick investigation of the main
results of this paper in this new context.

\begin{lemma} There exists an explicit pair $(F,\phi )$ of an 
algebraically closed field $F$ and an automorphism $\phi$, such that
$F\supset \widehat{K}$ and $\phi$ extends the given automorphism of 
$\widehat{K}$. The pair $(F,\phi )$ induces automorphisms of the algebraic
closures of $\widehat{K}$ and $K$ extending the given $\phi$.
\end{lemma}
\begin{proof} 
 The algebraic closure of $\widehat{K}$ has no explicit
description. However, there is an explicit algebraically closed field 
$F:={\mathbf C}((z^{\mathbb Q}))$ containing $\widehat{K}$. The elements of this field
are expressions $\sum _{\lambda \in {\mathbb Q}}a_\lambda z^{\lambda}$ with
all $a_\lambda \in {\mathbf C}$ and such that $\{\lambda |\ a_\lambda \neq 0\}$ is
a well ordered subset of ${\mathbb Q}$. It is well known that $F$ is a
maximally complete field with residue field $\mathbf C$ and value group 
${\mathbb Q}$.
In particular, $F$ is algebraically closed. Choose a homomorphism
$\lambda \mapsto q^\lambda$ from ${\mathbb Q}$ to ${\mathbf C}^*$ with $q^1=q$.
One defines an automorphism $\phi$ of $F$ by the formula
$\phi ( \sum _{\lambda \in {\mathbb Q}}a_\lambda z^{\lambda})=
\sum _{\lambda \in {\mathbb Q}}a_\lambda q^\lambda z^{\lambda}$.
This extends the action of $\phi$ on $\widehat{K}$. The algebraic closures
of $\widehat{K}$ and $K$ can be seen as subfields of $F$. They are obviously
invariant under $\phi$. This proves the assertion.\end{proof}

\begin{lemma} Let $\widehat{K}\subset L$ be an extension of degree $m<\infty $
such that $\phi$ extends to an automorphism of $L$. Then 
$L=\widehat{K}(z^{1/m})$. A similar statement holds for $K$ replacing
$\widehat{K}$.
\end{lemma}
\begin{proof} Write $L={\mathbf C}((t))$. Then 
$z=a_mt^m+a_{m+1}t^{m+1}+\cdots$ with $a_m\neq 0$.
The action of $\phi$ on $L$ has therefore the form $\phi (t)=q_1t+\cdots$
with $q_1^m=q$. Since $0<|q_1|<1$, one can produce an element 
$s\in {\mathbf C}[[t]]$ such that ${\mathbf C}[[s]]={\mathbf C}[[t]]$ and
$\phi (s)=q_1s$. Thus we may assume that $\phi (t)=q_1t$. Then 
\[q(a_mt^m+a_{m+1}t^{m+1}+\cdots )=qz=\phi (z)=a_mq_1^mt^m
+a_{m+1}q_1^{m+1}t^{m+1}+\cdots \]
implies that $z=a_mt^m$. This proves the statement for the case $\widehat{K}$.

For the case $K$, one has to show that $L={\mathbf C}\{t\}$ contains an element
$s$ with ${\mathbf C}\{s\}={\mathbf C}\{t\}$ and $\phi (s)=s$. Write 
$\phi (t)=q_1t+a_2t^2+\cdots$ and $s=t+b_2t^2+b_3t^3+\cdots$. Then 
$\phi (s)=q_1s$ leads to a sequence of equations for the $b_i$. An inspection
shows that the convergence of the series $q_1t+\sum _{n\geq 2}a_nt^n$ implies
the convergence of the series $t+\sum _{n\geq 2}b_nt^n$.   \end{proof}

Most of the preceeding sections remain valid, after a small adaptation, in the
present context. We give now some details.

In section 1, the fields $K_\infty$ and $\widehat{K}_\infty$ should be read
here, not as the algebraic closures but as the fields $\cup _{n\geq 1}K_n$ with
$K_n:={\mathbf C}(\{z^{1/n}\})$ and $\cup _{n\geq 1}\widehat{K}_n$ with 
$\widehat{K}_n:={\mathbf C}((z^{1/n}))$. The complex function theory for $E_q$ 
is replaced by the rigid analytic theory (see for instance [Fr-vdP]) and the 
formulas in subsection 1.1 remain valid. The only part of section 1 that has
no (obvious) translation is subsubsection 1.4.1.

All of section 2 remains valid with the exception of Remarks 2.4. Indeed,
the formulas for the decomposition of tensor products of indecomposable
modules (or for indecomposable vector bundles on $E_q$) are different in
positive characteristic. Especially, the decomposition of $U_a\otimes _KU_b$ 
poses a non trivial combinatorial problem, solved in [At] in the context
of vector bundles on $E_q$.

No changes are needed for the results of sections 3 and 4. However we will
rewrite Section 5 completely by developing a suitable Picard-Vessiot theory,
calculating difference Galois groups and a universal Picard-Vessiot ring.
We take [vdP-S 1,2] as guide for this.

It is {\it not} possible to attach, as in Theorem 6.1, to regular singular 
difference modules $M$ over $K$, connections on $v(M)$. Indeed, for the module
$Ke$ with $\Phi e=ce$ and $c\in {\mathbf C}^*$, the connection must have the
form $\nabla e=a(c)\frac{dz}{z}\otimes e$, where 
$a:{\mathbf C}^*\rightarrow {\mathbf C}$ is a homomorphism and satisfies
$a(q)=1$. However, $a(q)=p\cdot a(q^{1/p})=0$.

\subsection{Picard-Vessiot theory and examples}
Consider a field $F$ provided with an automorphism $\phi$ of infinite order.
The field of constants $C:=\{f\in F|\ \phi (f)=f\}$ is supposed to be
algebraically closed. A difference module $M$ is a finite dimensional vector
space over $F$, provided with a bijective additive map $\Phi :M\rightarrow M$
satisfying $\Phi (fm)=\phi (f)\Phi (m)$. After choosing a basis of $M$ over
$F$, the equation $\Phi (m)=m$ translates into a matrix difference equation
$y=A\phi (y)$ with $A\in {\rm GL}_n(F)$. 

A Picard-Vessiot ring $PV$ for $M$ (or $y=A\phi (y)$) is a commutative
$F$-algebra with unit element satisfying
\begin{enumerate} 
\item An extension of $\phi$ as automorphism of $PV$ is given.
\item $PV$ has no $\phi$-invariant ideals except $\{0\}$ and $PV$.
\item There is a $U\in {\rm GL}_n(PV)$ with $U=A\phi (U)$.
\item $PV$ is generated over $F$ by the entries of $U$ and $U^{-1}$.  
\end{enumerate}
With the methods of [vdP-S 1,2] one shows
the existence and unicity of $PV$ (up to $\phi$-isomorphisms). One observes
that $PV$ is reduced and has in general zero divisors. The field of constants 
of the total ring of fractions of $PV$ is again $C$.

The naive definition of the difference Galois group $G$ of $M$ is: $G$ is
the `abstract' group of the $F$-automorphisms of $PV$ commuting with $\phi$. 
This definition is insufficient in positive characteristic.
 A more precise definition is the following.
One defines a covariant functor $\mathcal G$ from the category of the finitely
generated commutative $C$-algebras $R$ to the category of groups by ${\mathcal
 G}(R)$ is the group of the automorphisms of $R\otimes _CPV$ which are
$R\otimes _CF$-linear and commute with the action of $\phi$. It can be shown
that $\mathcal G$ is represented by an affine group scheme of finite type over
$C$. The difference Galois group $G$ is by definition this group scheme. 

If the field $F$ has characteristic zero then $G$ is a reduced linear algebraic
group. In our case, where $F$ has characteristic $p>0$, the difference Galois 
group need not be reduced. In the following examples we calculate the 
Picard-Vessiot rings and their difference Galois groups for some typical
equations. For all examples we take $F=\widehat{K}={\mathbf C}((z))$ with 
$\mathbf C$ as before. Since the examples are pure modules, one may replace
$\widehat{K}$ by $K={\mathbf C}(\{z\})$ everywhere.  For convenience,
we calculate for modules $M$ the `contravariant solutions', i.e., 
${\rm ker}(\Phi -1,{\rm Hom}_(M,PV))\;$, in stead of the `covariant solutions',
i.e., ${\rm ker}(\Phi-1, PV\otimes M)\;$.

\medskip

\noindent
{\it Example} 1. 
The extension $\widehat{K}\subset \widehat{K}(z^{1/p})$ is the Picard-Vessiot
extension for the difference equation $\phi (y)=q^{1/p}y$. Its difference
Galois group is the group $\mu _{p,{\mathbf C}}$. More generally,
$\widehat{K}(z^{1/n})$ is the Picard-Vessiot extension of an equation
$\phi (y)=q^{t/n}y$ with $(t,n)=1$ and its difference Galois group is 
$\mu _{n,{\mathbf C}}$. We recall that 
$\mu _{n,{\mathbf C}} ={\rm Spec}({\mathbf C}[t]/(t^n-1))$, with
co-multiplication given by $t\mapsto t\otimes t$. 

\medskip

\noindent {\it Example} 2. Equation $\phi (y)=cy$
with $c\in {\mathbf C}^*$. Suppose that for all $n\geq 1$ the only solution
of $\phi (y)=c^ny$ in $\widehat{K}$ is $y=0$, then the Picard-Vessiot 
extension is $\widehat{K}[Y,Y^{-1}]$ with $\phi (Y)=cY$ and the difference
Galois group is ${\mathbb G}_{m,{\mathbf C}}$.\\
Suppose that there exists a non zero $y\in \widehat{K}$ such that
$\phi (y)=c^ay$ for some $a\geq 1$. Then $c$ has the form $\zeta q^{t/n}$
with $\zeta$ a primitive $d$th root of unity and $n\geq 1,\ (t,n)=1$. 
We consider the two equations $\phi (y)=q^{t/n}y$ and $\phi (y)=\zeta y$
separately. The first equation is considered in example 1. The second equation
has Picard-Vessiot ring $\widehat{K}[y]$ with equation $y^d=1$ and 
$\phi (y)=\zeta y$. This ring has obviously zero divisors (if $d>1$).
Its difference Galois group is $\mu _{d,{\mathbf C}}
\cong {\mathbb Z}/{\mathbb Z}d$ over $\mathbf C$ since 
$d$ is not divisible by $p$. The Picard-Vessiot
extension for $\phi (y)=cy$ is a subring of $\widehat{K}[z^{1/n}][y]$. The
difference Galois group is therefore a quotient of the group
$\mu _{n,{\mathbf C}}\times {\mathbb Z}/{\mathbb Z}d$.

\medskip

\noindent {\it Example} 3. $U_m=\widehat{K}\otimes V$ with
$\dim V=m$, $\Phi (1\otimes v)=1\otimes U(v)$ and $U\in {\rm GL}(V)$ the
indecomposable unipotent operator. There exists an element $e\in V$ such that
$e, (U-1)e,\dots ,(U-1)^{m-1}e$ is a basis of $V$. Hence $U_m\cong
\widehat{K}[\Phi ,\Phi ^{-1}]/\widehat{K}[\Phi ,\Phi ^{-1}]((\Phi -1)^m)$.
Thus we have to find the Picard-Vessiot ring for the equation 
$(\phi -1)^m(y)=0$.

\smallskip

(a) Suppose that $1<m\leq p$. The difference ring 
$A_1:=\widehat{K}[\ell ]$, defined by $\ell ^p-\ell =0$ and 
$\phi (\ell )=\ell +1$, has only trivial $\phi$-invariant ideals. The elements
${\ell \choose i}$ for $i=0,\dots ,m-1$ are $\mathbf C$-linear independent 
solutions of $(\phi -1)^m(y)=0$. 
Since $A_1$ is generated over $\widehat{K}$ by $\ell$ one
finds that $A_1$ is the Picard-Vessiot extension for the equation 
$(\phi -1)^m(y)=0$. Let $R$ be a $\mathbf C$-algebra and $\sigma$
a $R\otimes _{\mathbf C}$ automorphism of $R\otimes _{\mathbf C}K[\ell ]$,
which commute with $\phi$. Then $\sigma$ is determined by $\sigma (\ell )$
and $\sigma (\ell )=\ell +a$ where $a$ is any element in $R$ with $a^p=a$.
The difference Galois group is therefore the group 
${\mathbb Z}/{\mathbb Z}p$ over ${\mathbf C}$. 
In view of further equations we write $\ell =\ell _1$.

\smallskip

(b) Suppose that $p<m\leq p^2$. The Picard-Vessiot ring for 
$(\phi -1)^m(y)=0$ is $A_2:=\widehat{K}[\ell _1,\ell _2]$, defined by 
$\ell _1^p-\ell_1=0,\ \ell _2^p-\ell _2=0$, $\phi (\ell _1)=\ell_1+1,\ 
\phi (\ell _2)=\ell_2+{\ell _1\choose p-1}$. The set of maximal ideals of $A_2$
is $\{(\ell _1-a,\ell _2-b)| a,b\in {\mathbb F}_p\}\cong {\mathbb F}_p^2$.
One calculates that $\phi$ acts transitively on this set and one concludes
that $A_2$ has only trivial $\phi$-invariant ideals. One observes that
\[(\phi -1)^p(\ell_2)=(\phi -1)^{p-1}({\ell _1\choose p-1})=1 
\mbox{ and } 
(\phi -1)^{p^2-1}{\ell _2\choose p-1}{\ell _1\choose p-1}=1\ .\]
The conclusion is that ${\mathbf C}[\ell _1,\ell _2]$ is the kernel of 
$(\phi -1)^{p^2}$, acting on $A_2$. 
This shows that $A_2$ is generated by the solutions of $(\phi -1)^my=0$ and
thus that $A_2$ is indeed the Picard-Vessiot ring.

We have to represent the functor ${\mathcal G}$, given by
${\mathcal G}(R)$ is the group of the difference automorphism of 
$R\otimes _{\mathbf C}A_2$ over $R\otimes _{\mathbf C}K$. For the calculation
of ${\mathcal G}(R)$ we suppose for convenience that ${\rm Spec}(R)$ is
connected. Then $a\in R,\ a^p=a$ implies $a\in {\mathbb F}_p$. Further  
${\mathbb F}_p[\ell _1,\ell _2]=
\{\xi \in R\otimes _{\mathbf C}A_2|\ \xi ^p=\xi \mbox{ and } 
(\phi -1)^{p^2}\xi =0\}$. Thus any $\sigma \in {\mathcal G}(R)$ induces an
automorphism of ${\mathbb F}_p[\ell _1,\ell _2]$ commuting with $\phi$.
On the other hand, any automorphism of  ${\mathbb F}_p[\ell _1,\ell _2]$ 
commuting with $\phi$, extends uniquely to an element of ${\mathcal G}(R)$.

The algebra ${\mathbb F}_p[\ell _1,\ell _2]$ is seen as 
${\mathbb F}_p[\phi ]$-module. The element 
$\xi ={\ell _2\choose p-1 }\cdot {\ell _1\choose p-1 }$ is a generator and
the module can be written as ${\mathbb F}_p[t]/(t^{p^2})$, where $t=\phi -1$.
We recall that $t^{p^2-1}\xi =1$. Thus an
automorphism $\sigma$ satisfies $\sigma \xi =a\xi$ with 
$a\in {\mathbb F}_p[t]/(t^{p^2})$ and $a\equiv 1 \mod (t)$. On the other hand
any $a$, as above, produces a unique automorphism. The above group
is the cyclic group of order $p^2$, generated by $a=1+t\mod (t^{p^2})$. Thus
the difference Galois group is equal to the group ${\mathbb Z}/{\mathbb Z}p^2$
over $\mathbf C$.

\smallskip 
(c) In a similar way one obtains that the difference ring 
$A_k:=\widehat{K}[\ell _1,\cdots ,\ell _k]$, given by the equations
$\ell _i^p-\ell _i=0$ for $i=1,\dots ,k$ and 
\[(\phi -1)(\ell _i)={\ell _{i-1}\choose p-1}{\ell _{i-2}\choose p-1}\cdots
{\ell _1\choose p-1} \]
for $i=2,\dots ,k$ and $(\phi -1)\ell _1=1$,
is the Picard-Vessiot ring for $(\phi -1)^my=0$ for $m$ such that
$p^{k-1}<m\leq p^k$. Its difference Galois group
is the group ${\mathbb Z}/p^k{\mathbb Z}$ over ${\mathbf C}$. As in the case
$k=2$, the difference Galois group is identified with the group of the
automorphisms of $Z:={\mathbb F}_p[\ell _1,\dots ,\ell _k]$ which commute with
$\phi$. Now $Z$, as a ${\mathbb F}_p[\phi ]={\mathbb F}_p[t]$-module
(with $t=\phi -1$), has 
$\xi ={\ell _k\choose p-1}{\ell _{k-1}\choose p-1}\cdots {\ell _1\choose p-1}$
as cyclic element, and  is isomorphic to 
${\mathbb F}_p[t]/(t^{p^k})$. The automorphisms are given by the
elements $a\in {\mathbb F}_p[t]/(t^{p^k})$ with $a\equiv 1 \mod (t)$.
This group is cyclic of order $p^k$ and has generator $a=1+t\mod (t^{p^k})$.

\medskip

\noindent 
{\it Example} 4. $M=\widehat{K}(z^{1/n})e$, $\Phi e=z^{t/n}e$
with $n\geq 1$, $(t,n)=(p,n)=1$. The corresponding scalar equation is
$\phi ^n(y)=q^{t(n-1)/2}z^ty$. By definition, $PV$ contains an invertible
element $\alpha$ satisfying the equation. Any other solution $y$ has the form
$\tilde{y}\alpha$ with $\phi ^n(\tilde{y})=\tilde{y}$. Hence $PV$ contains
$\widehat{K}[y_1]$ with $y_1^n=1$ and $\phi (y_1)=\zeta _ny_1$ where
$\zeta _n$ is a primitive $n$th root of unity. The invertible element 
$u:=\phi (\alpha )\alpha ^{-1}$ satisfies the equation 
$\phi ^n(u)=q^tu$. All solutions of the latter equation have the form 
$z^{t/n}y$ with $\phi ^ny=y$. Thus $y$ is an invertible element of
${\mathbf C}[y_1]$. It follows that $z^{t/n}$ and $z^{1/n}$ are in $PV$. 
After changing $\alpha$ we may suppose that $u=z^{t/n}$. This leads to
the assertion that  the Picard-Vessiot ring is
$\widehat{K}(z^{1/n})[y_1,\alpha, \alpha ^{-1}]$ with the rules: 
$y_1^n=1,\ \phi (y_1)=\zeta _ny_1$ (with $\zeta _n$ a primitive $n$th root
of unity); $\alpha$ transcendental over 
$\widehat{K}$ and $\phi (\alpha )=z^{t/n}\alpha$.
The elements $\{y_1^i\alpha |\ i=0,\dots ,n-1\}$ form a ${\mathbf C}$-basis of
solutions. The inclusion $\widehat{K}[z^{1/n}][y_1]\subset 
\widehat{K}[z^{1/n}][y_1,\alpha ,\alpha ^{-1}]$ induces an exact sequence for
the difference Galois group, namely
\[1\rightarrow {\mathbb G}_m\rightarrow G\rightarrow \mu _n\times \mu
_n\rightarrow 1 \ .\]
The term $\mu _n\times \mu _n$ is the difference Galois group of the 
Picard-Vessiot extension $\widehat{K}(z^{1/n})[y_1]$ of the equation 
$\phi ^n(y)=qy$.

\medskip

\noindent 
{\it Example} 5. One considers the difference module
$M=\widehat{K}(z^{1/p})e$ with $\Phi e=z^{t/p}e$ (and $1\leq t<p$), seen as
difference module of dimension $p$ over $\widehat{K}$. A corresponding scalar
equation is $\phi ^p(y)=q^{t(n-1)/2}z^ty$. The method of example 4 yields
that the Picard-Vessiot
extension for $M$ is $PV=\widehat{K}(z^{1/p})[\ell _1,\alpha ,\alpha ^{-1}]$,
satisfying the following rules: $\ell _1 ^p-\ell _1 =0$ and $\phi (\ell _1
)=\ell _1 +1$, $\alpha$ is transcendental over $\widehat{K}$ and 
$\phi (\alpha )=z^{t/p}\alpha$. The inclusion 
$\widehat{K}(z^{1/p})[\ell _1]\subset \widehat{K}(z^{1/p})[\ell _1,\alpha
,\alpha ^{-1}]$, induces an exact for the difference Galois group $G$
\[1\rightarrow {\mathbb G}_m\rightarrow G\rightarrow {\mathbb Z}/p{\mathbb Z}
\times \mu _p\rightarrow 1 \ .\]
The group ${\mathbb Z}/p{\mathbb Z}\times \mu _p$ is 
the difference Galois group of the Picard-Vessiot extension
$\widehat{K}(z^{1/p})[\ell _1]$ of the equation $\phi ^p(y)=qy$.

\medskip

\noindent 
{\it Example} 6. $M=\widehat{K}(z^{1/p^k})e$ with $\Phi e=z^{t/p^k}e$ and
$(t,p)=1$. A corresponding scalar equation is 
$\phi ^{p^k}(y)=q^{t(p^k-1)/2}z^ty$. The Picard-Vessiot ring is
$\widehat{K}(z^{1/p^k})[\ell _1,\dots ,\ell _k,\alpha ,\alpha ^{-1}]$ with
$\phi (\alpha )=z^{t/p^k}\alpha$. The 
difference Galois group admits an exact sequence
\[1\rightarrow {\mathbb G}_m\rightarrow G\rightarrow 
{\mathbb Z}/p^k{\mathbb Z} \times \mu _{p^k}\rightarrow 1 \ .\]

\subsection{The universal difference ring over $\widehat{K}$} 
The above examples and the classification of the indecomposable difference 
modules over $\widehat{K}$ lead to the following description of the universal
Picard-Vessiot ring $Univ$ of $\widehat{K}$. Let $\widehat{K}^+$ denote the 
union of the fields $\widehat{K}(z^{1/n})$ (for all $n\geq 1$). One introduces
symbols $e(cz^\lambda )$ for $c\in {\mathbf C}^*,\ \lambda \in {\mathbb Q}$
and $\ell _k$ for $k\in {\mathbb Z},\ k\geq 1$. Then 
$Univ=\widehat{K}^+[\{\ell _k\}_{k\geq 1}\{e(cz^\lambda )\}]$. The only 
relations are $e(c_1z^{\lambda _1})\cdot e(c_2z^{\lambda _2})=
e(c_1c_2z^{\lambda _1+\lambda _2})$; $e(q^\lambda )=z^{-\lambda}$ for
$\lambda \in {\mathbb Q}$ (including $e(1)=1$); $\ell _k^p-\ell _k=0$ for all
$k\geq 1$. The action of $\phi$ on $Univ$ is given by:\\
$\phi$ acts on $\widehat{K}^+$ by $\phi (\sum a_\lambda z^\lambda )=
\sum a_\lambda q^\lambda z^\lambda , \ \ \ $
$\ \ \ cz^\lambda \phi (e(cz^\lambda ))= e(cz^\lambda ),\ \ \ $
$\phi \ell _1=\ell _1+1,\ \ \ $ 
$(\phi -1)\ell _k={\ell _{k-1}\choose p-1}{\ell _{k-2}\choose p-1}\cdots
{\ell _1\choose p-1}$ for $k\geq 2$.   

\smallskip

\noindent 
$Univ$ can also be described as the algebra 
$\widehat{K}[\{\ell _k\}_{k\geq 1},\{e(c)\}_{c\in {\mathbf C}^*},
\{e(z^\lambda )\}_{\lambda \in {\mathbb Q}}]$ with the relations:
$\ell _k^p=\ell _k$ for all $k\geq 1$; $e(c_1)e(c_2)=e(c_1c_2)$;
$e(1)=1,\ e(q)=z^{-1}$; $e(z^{\lambda _1})e(z^{\lambda _2})
=e(z^{\lambda _1+\lambda _2})$. The action of $\phi$ is given by
the above formulas for $\phi (\ell _k)$ and $\phi (e(c))$ and
$e(q^{-\lambda })\phi (e(z^{\lambda}))=e(z^\lambda )$. 

The universal
Picard-Vessiot ring $Univ_{rs}$ for the regular singular difference equations
over $\widehat{K}$ is
the subring $\widehat{K}[\{\ell _k \}_{k\geq 1},\{e(c)\}_{c\in {\mathbf C}}]$.

The difference
Galois group ${\mathbb G}_{rs}$ of $Univ_{rs}$ can be identified with 
${\mathbb Z}_p\times {\rm Hom}({\mathbf C}^*/q^{\mathbb Z},{\mathbf C}^*)$.
The factor ${\mathbb Z}_p$ is the difference Galois group of 
$\widehat{K}[\{ \ell _k  \}_{k \geq 1} ]$. The notation
${\rm Hom}({\mathbf C}^*/q^{\mathbb Z},{\mathbf C}^*)$ for the second factor
is an `abus de langage'. It  hides the non reduced part of the difference
Galois group. One writes $E_q:={\mathbf C}^*/q^{\mathbb Z}$ (with brute force)
as a product $D\times (E_q)_{tors}$, where $D$ is a divisible, torsion free
group. Since $D$ is a vector space over ${\mathbb Q}$, the term 
${\rm Hom}(D,{\mathbf C}^*)$ defines a reduced affine group scheme. More
precisely, $D$ is a direct limit of its free finitely generated subgroups
and the affine group ${\rm Hom}(D,{\mathbf C}^*)$ is the projective limit
of algebraic tori over ${\mathbf C}$.

 The group 
$(E_q)_{tors}$ is a product 
$\{a\in {\mathbf C}^*| a\mbox{ root of unity}\}\times 
q^{\mathbb Q}/q^{\mathbb Z}$. This torsion group is isomorphic to
$({\mathbb Q}/{\mathbb Z}[1/p])^2\times {\mathbb Q}_p/{\mathbb Z}_p$. The
first factor yields the reduced affine group scheme which is the projective
limit of $\mu _n\times \mu _n$, taking over all $n\geq 1$ not divisible by
$p$. The second factor yields the non reduced affine group scheme which is the
projective limit of the groups $\mu _{p^k}$.

The difference Galois group ${\mathbb G}_{univ}$ of $Univ$ admits an exact sequence
\[1\rightarrow {\rm Hom}({\mathbb Q},{\mathbf C}^*)  \rightarrow
{\mathbb G}_{univ}\rightarrow {\mathbb G}_{rs}\rightarrow 1 \]
The term ${\rm Hom}({\mathbb Q},{\mathbf C}^*)$ is the affine group scheme
which represents the $Univ_{rs}$-linear automorphisms of $Univ$, commuting
with $\phi$. We note that the affine group scheme ${\mathbb G}_{univ}$ is not
commutative. Unlike the complex case we do not find explicit topological
generators (like $\Gamma ,\Delta ,D$) for ${\mathbb G}_{univ}$. This is
due to the complicated structure of the group ${\mathbf C}^*$.

In the above we established a reasonable Picard-Vessiot theory and explicit
calculations for the difference Galois groups of difference modules over
$\widehat{K}$ (or equivalently for split difference modules over $K$). 
The explicit determination of the difference Galois group of (non split)
difference modules over $K$ can be copied from Section 5.

\vspace{1cm}

\noindent {\bf References}\\

\noindent [At] M.F.~Atiyah -{\it Vector bundles on elliptic curves } -
Proc.London.Math. Soc. {\bf 7}, 1957, 414-452.\\
\noindent [Bir] Collected mathematical papers of George David Birkhoff,
Volume 1. Dover publications, 1968.\\
\noindent [For] O.~Forster - {\bf Riemannsche Fl\"ache} - Heidelberger 
Taschenb\"ucher, Springer Verlag 1977.\\
\noindent [Fr-vdP] J.~Fresnel and M.~van der Put - {\bf Rigid Analytic
  Geometry and its Applications } - Progress in Math. {\bf 218}, 2004.\\
\noindent [vdP] M.~van der Put - {\it Skew differential fields, differential
  and difference equations} - Ast\'erisque {\bf 296}, 2004, p. 191-207\\
\noindent [vdP-R] M.~van der Put and M.~Reversat - {\it Krichever modules for
  difference and differential equations } - Ast\'erisque {\bf 296}, 2004,
  207-227.\\  
\noindent [vdP-S1] M. van der Put and M.F.~Singer -{\bf Galois theory of
  difference equations} - Lect. Notes in Math., vol 1666, Springer Verlag,
  1997. \\
\noindent [vdP-S2] M. van der Put and M.F.~Singer - {\bf Galois theory of
linear differential equations} - Grundlehren der mathematische Wissenschaften
{\bf 328}, Springer Verlag, 2003.\\
\noindent [R-S-Z] J.-P. Ramis, J. Sauloy and C. Zhang, - {\it La vari\'et\'e des classes 
analytiques d'\'equations aux $q$-diff\'erences dans une classe formelle } -
 C.R.Acad.Sci.Paris, Ser. I 338 (2004)\\
\noindent  [Sau1] J.~Sauloy - {\it Galois theory of fuchsian $q$-difference equations} - Ann. Sci. \'Ec. Norm. Sup. 
4$^{e}$ s\'erie {\bf 36}, no 6, 2003,   p. 925-968\\   
\noindent  [Sau2] J.~Sauloy - {\it Algebraic construction of the Stokes sheaf
for   irregular  linear $q$-difference equations} - Ast\'erisque {\bf 296},
2004,   p. 227-251\\  
\noindent [Sau3], J.~Sauloy -{\it  La filtration canonique par les
pentes d'un module aux $q$-diff\'erences et le gradu\'e associ\'e} - to appear
in Ann.Inst.Fourier. 2004.\\

\end{document}